\documentclass[final,3p,times, num,sort&compress]{elsarticle}

\usepackage{amssymb}
\usepackage{amsmath}
\usepackage{multibib}
\usepackage{mathtools}
\usepackage{caption}
\usepackage{lineno}
\usepackage{subcaption}
\usepackage{float}
\numberwithin{equation}{section}
\graphicspath{{Results/}}
\usepackage{soul}
\usepackage{xcolor}
\usepackage{hyperref}

 \usepackage{lineno}
\newcommand{\dx}{{\Delta{x}}}
\newcommand{\dt}{{\Delta{t}}}
\newcommand{\eps}{\epsilon}
\journal{Computational Physics}
\begin{document}
\begin{frontmatter}



\title{Numerical Artifacts in the Discontinuous Generalized Porous Medium Equation: How to Avoid Spurious Temporal Oscillations}

\author[label1]{Danielle C. Maddix}
\author[label2]{Luiz Sampaio}
\author[label1,label2]{Margot Gerritsen}
\address[label1]{Institute of Computational and Mathematical Engineering}
\address[label2]{Energy Resources Engineering \\Stanford University}

\begin{abstract} 
Numerical discretizations of the Generalized Porous Medium Equation (GPME) with discontinuous coefficients are analyzed with respect to the formation of numerical artifacts.  In addition to the degeneracy and self-sharpening of the GPME with continuous coefficients, detailed in \cite{maddix_pme}, increased numerical challenges occur in the discontinuous coefficients case.  These numerical challenges manifest themselves in spurious temporal oscillations in second order finite volume discretizations with both arithmetic and harmonic averaging.  The integral average, developed in \cite{vandermeer2016}, leads to improved solutions with monotone and reduced amplitude temporal oscillations.  In this paper, we propose a new method called the Shock-Based Averaging Method (SAM) that incorporates the shock position into the numerical scheme.  The shock position is numerically calculated by discretizing the theoretical speed of the front from the GPME theory.  The speed satisfies the jump condition for integral conservation laws.  SAM results in a non-oscillatory temporal profile, producing physically valid numerical results.   We use SAM to demonstrate that the choice of averaging alone is not the cause of the oscillations, and that the shock position must be a part of the numerical scheme to avoid the artifacts.
\end{abstract}


\end{frontmatter}


\section{Introduction}
\label{intro}

The purpose of the paper is to identify the cause of the numerical artifacts reported in the literature \cite{maddix_pme, lipnikov2016,vandermeer2016} for second order finite volume discretizations of the Generalized Porous Medium Equation (GPME) with discontinuous coefficients, and to suggest a numerical approach that does not have these problems.
The GPME, commonly known as the Filtration Equation, can be expressed in both 
conservative and integral forms as: 
\begin{equation}
\begin{aligned}
p_t &= \nabla \cdot (k(p) \nabla p) 
 \\ &= \Delta \Phi (p),
\hspace{.25cm} \text{where} \\
\vspace{1cm}
\Phi(p) &= \int_0^p k(\tilde p) d\tilde p, \hspace{.25cm}  k(p) = \Phi'(p).
 \end{aligned}
 \label{eq:GPME}
\end{equation}
In the discontinuous coefficients case, $k(p)$ is given as:
\begin{equation}
\begin{aligned}
	k(p)& = 
	\begin{dcases}
  		k_{\max},           & p \ge p^* \\
    		k_{\min},              &p^* > p, \hspace{.25cm} \text{where}
	\end{dcases}
\label{eq:discont_k}
\end{aligned}
\end{equation} 
$k_{\max}$, $k_{\min}$ and $p^*$ are real positive constants.

In this paper, we are interested in a subclass of the GPME known as the Stefan problem \cite{stefan91, rubinstein71, meirmanov92, vazquez2007}, for which
	\begin{equation}
				\Phi(p) = \begin{dcases}
  									c_1(p - c_3)_+,          & \text{if} \hspace{.1cm} p \ge 0, \\
    									c_2p ,             & \text{otherwise},
						\end{dcases}
					\label{eq:stef_gen}
	\end{equation}
for arbitrary $c_1, c_2, c_3 \in \mathbb{R}$.  We look at a particular Stefan problem, where $p \ge 0, c_1 = k_{\max}$ and $c_3 = p^*$.  
Then, $\Phi(p)$ in Eqn. \eqref{eq:GPME} can be expressed as the positive part function
	\begin{equation}
		\Phi(p) = k_{\max}(p-p^*)_+ = \begin{dcases}
  											k_{\max}(p-p^*),          & \text{if} \hspace{.1cm} p \ge p^*,  \\
    											0,              & \text{otherwise},
									\end{dcases}
		\label{eq:stef_sp}
	\end{equation}
	and $k(p)$ is given by Eqn. \eqref{eq:discont_k} with $k_{\min} = 0$.
	
 The Stefan problem can also be formulated in its classical form \cite{sethian88,osher97} as  
 \[
\begin{aligned}
  		\frac{\partial p}{\partial t}= k_{\max} \Delta p,   \hspace{0.25cm}        & p \ge p^*, \\
\end{aligned}
\]
\[
\begin{aligned}
    		\frac{\partial p}{\partial t} = k_{\min} \Delta p,      \hspace{0.25cm}        & p < p^*. \\       
\end{aligned}
\]
In the classical Stefan problem, the above two parabolic equations are defined on domains that are separated by a moving interface $x^*(t)$.
The Stefan condition is given on the interface as
\[
	(p_L -p_R)\frac{dx^*(t)}{dt} = -k_{\max}\frac{\partial p_L}{dx}+k_{\min}\frac{\partial p_R}{dx},
\]
where $p_L \equiv \lim_{x \rightarrow x^*(t)^-}p(x,t)$ and $p_R \equiv \lim_{x \rightarrow x^*(t)^+}p(x,t)$.
The Stefan condition can be derived by using the Rankine-Hugoniot jump condition for this conservation law \cite{myers15}.

The Stefan Problem has been key to both the numerical and theoretical developments of the GPME.   The Stefan Problem is used in modeling phase transitions, and was developed to study the evolution of a medium of two phases, water and ice \cite{vazquez2007}.  \citet{brattkus92} use the Stefan problem to model crystal growth.  \citet{sethian92} and \citet{osher97} develop a modified Stefan problem to model crystal growth as well as dendritic solidification.  Eqn. \eqref{eq:discont_k} is often used to illustrate the numerical challenges present in more complex porous media applications. 
For example, in \citet{vandermeer2016} a foam model prototype is developed with $k_{\max} = 1$, $k_{\min} = \eps \rightarrow 0$ and $p^* = 0.5$. 
 
  The continuous GPME already poses numerical challenges, caused by self-sharpening and degeneracy for near-zero $k(p)$ \cite{maddix_pme}.  
  Due to this self-sharpening and degeneracy in the continuous case, it is not just the discontinuity that poses the numerical challenges in the Stefan problem in Eqn. \eqref{eq:stef_sp}.  
The discontinuity in $k(p)$ does make the challenges more severe and because of this discontinuity, the Modified Equation Analysis approach in \cite{maddix_pme} for the continuous GPME is not applicable. Here, we look at the discontinuous GPME, and also refer to subclasses of the continuous GPME.  The Porous Medium Equation (PME) subclass, where 
$
k(p) = p^m \text{ and }  m \ge 1,
$
is used to model gas flow through a porous medium \cite{maddix_pme, vazquez2007, ngo2016}.  Another application in thermodynamics is the superslow diffusion equation \cite{vazquez2007, maddix_pme}, where
$
	k(p)  = \exp(-1/p)$.  Further applications of the GPME for continuous $k(p)$ are detailed in \cite{maddix_pme}.  



\subsection{Understanding the Behavior of the GPME} 

\begin{figure}[H]
		\center
		\includegraphics[width =0.49\textwidth]{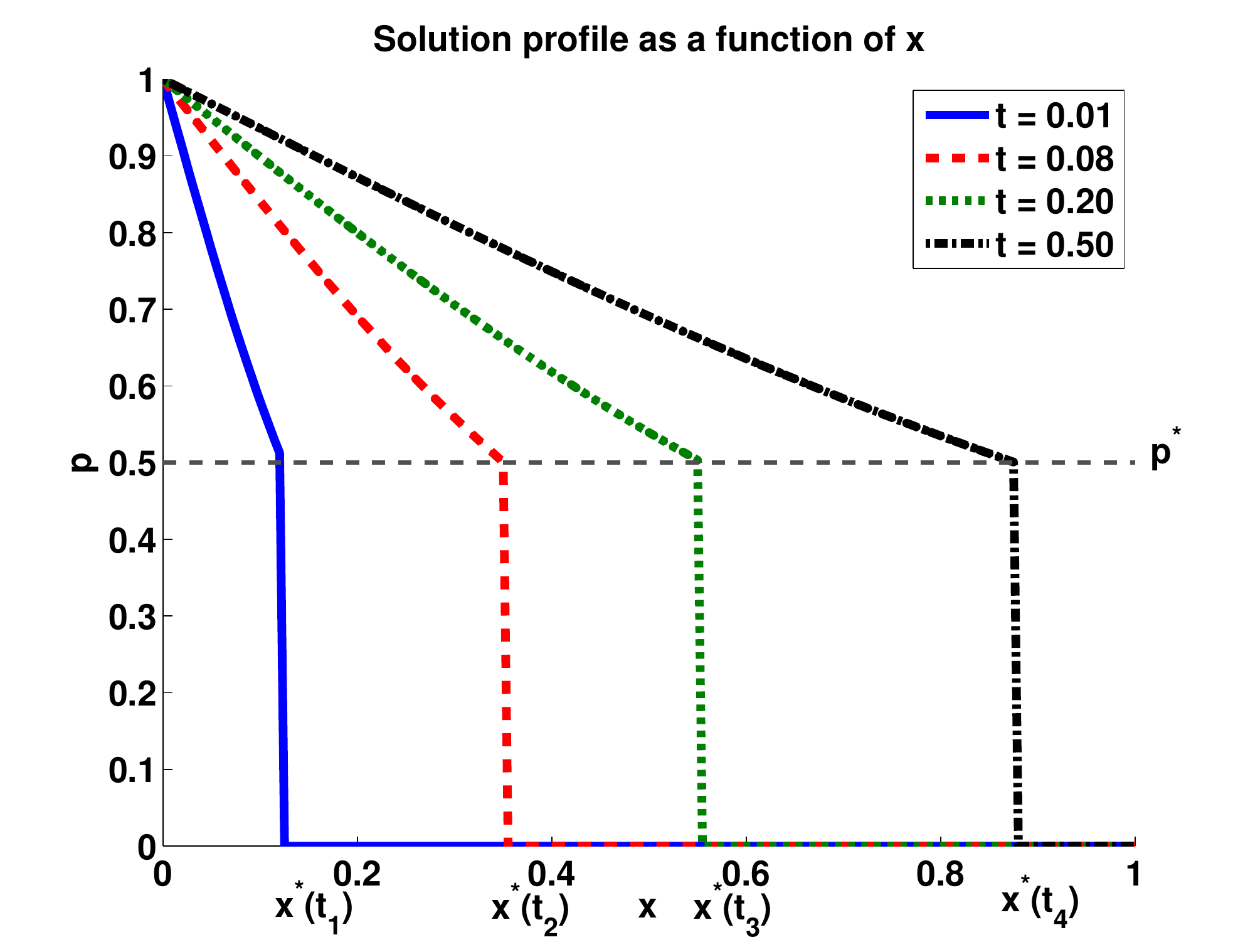}
		\caption{The exact solution of Eqn. \eqref{eq:GPME} with $k(p)$ given by Eqn. \eqref{eq:discont_k}, where $k_{\max} = 1$ and $k_{\min} = 0$, evolves as a rightward moving shock over time.  The moving shock position is given by $x^*(t)$ and the fixed $p$ value at the shock is given by $p^* = 0.5$.}
		\label{exactsol_time}
\end{figure}

The GPME in Eqn. \eqref{eq:GPME} at first appears like a heat equation.  Contrary to solutions of the heat equation, where the propagation speed is infinite, solutions of the GPME are known to have a finite speed of propagation \cite{vazquez2007}.  This results from the degeneracy of the GPME for compactly supported initial data, and is a property that distinguishes the GPME from classical parabolic theory.  This degeneracy leads to self-sharpening and moving interface solutions, as illustrated in Figure \ref{exactsol_time}.  In addition, for certain compactly supported smooth initial data, the waiting time phenomenon can occur.  This is discussed in the literature \cite{fischer15, angenent88} \cite[Chapter~3]{antontsev15} and is illustrated in Figure \ref{wait_time}.  The term waiting time refers to the fact that an interface moves only after sufficient sharpening.

\label{rh_cond}

The theoretical propagation speed for compactly supported initial data is known and is given by Darcy's Law \cite{vazquez2007} as
\[
	V = -\hspace{-.2cm}\lim_{x \rightarrow x^*(t)^-}\nabla v,
\] where
\begin{equation}	
	v = \int_0^p \frac{\Phi'(\tilde p)}{\tilde p} d\tilde p ,
	\label{eq:vel}
\end{equation} and $x^*(t)$ is the shock position.  Combining these expressions gives	\begin{equation}
		V = -\hspace{-.2cm}\lim_{x \rightarrow x^*(t)^-}\frac{\Phi'(p) \nabla p}{p} =  -\hspace{-.2cm}\lim_{x \rightarrow x^*(t)^-}\frac{k(p) \nabla p}{p}.
		\label{eq:vel_final}
	\end{equation}
	
Eqn. \eqref{eq:vel_final} holds for any $k(p)$.  
In the particular case of the Stefan problem in Eqn. \eqref{eq:stef_sp}, 
 $k(p) = k_{\max} = 1.0$ to the left of the shock.  Substituting the $k(p)$ limit into Eqn. \eqref{eq:vel_final} gives 
\begin{equation}
	V = -\hspace{-.2cm}\lim_{x \rightarrow x^*(t)^-}\frac{\nabla p}{p} = -\hspace{-.2cm}\lim_{x \rightarrow x^*(t)^-}\nabla \log(p).
	\label{eq:front_speed}
\end{equation}	

 The expression for the velocity $V$ can also be expressed in terms of fluxes.  The flux $F$ for the integral conservation law in Eqn. \eqref{eq:GPME} is given by 
\begin{equation}
	F(p) = -k(p) \nabla p =  -\nabla \Phi(p).
	\label{eq:anal_flux}
\end{equation}
Then, from Eqn. \eqref{eq:vel_final}, 
\begin{equation}
 	V = \frac{F(p_L)}{p_L},
	\label{eq:RH_pre}
\end{equation}
 where $p_L \equiv \lim_{x \rightarrow x^*(t)^-}p(x,t)$.  For a compactly supported initial condition and $k_{\min} = 0$, the flux and $p$ values to the right of the shock are zero.  
In this case, the velocity in Eqn. \eqref{eq:RH_pre} can be expressed in terms of the familiar jump condition for integral conservation laws as\begin{equation}
 	V = \frac{F(p_L) - F(p_R)}{p_L - p_R},
	\label{eq:RH}
\end{equation}
where $p_R \equiv \lim_{x \rightarrow x^*(t)^+}p(x,t)$ \cite{vazquez05, vazquez2007}.  The velocity in Eqn. \eqref{eq:RH} holds in general for any $k_{\min}, p_R \ge 0$, by the Rankine-Hugoniot condition \cite{lax57}.



\subsection{Numerical Methods to Approach this Problem} 
\label{lit}
\begin{figure}[H]
		\center
		\begin{subfigure}[H]{.33\textwidth}  
			\includegraphics[width =\textwidth]{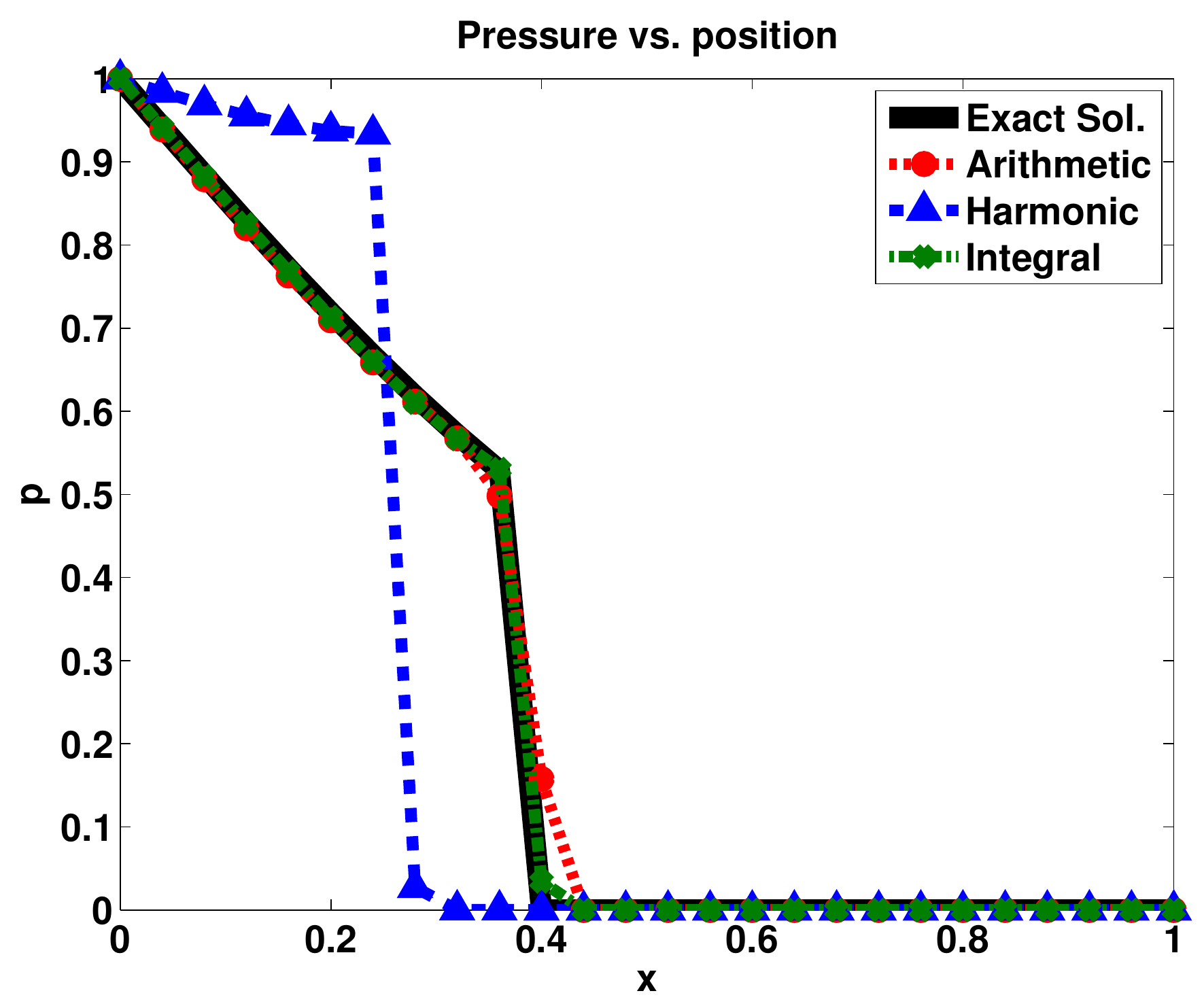}
			\caption{Spatial profiles at time $t = 0.05$.}
			\label{art_space}
		\end{subfigure}
		\begin{subfigure}[H]{.33\textwidth}  
			\includegraphics[width =\textwidth]{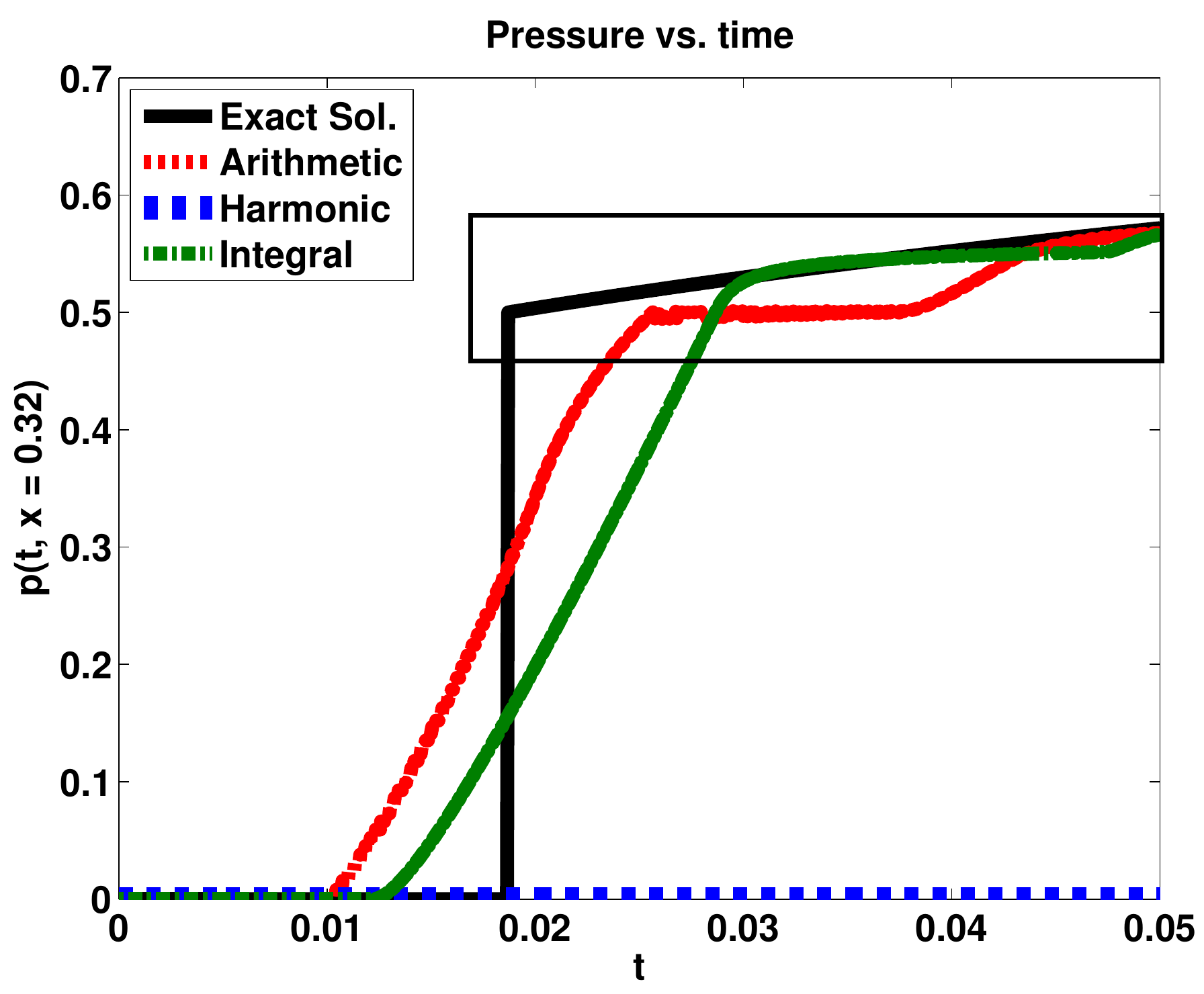}
			\caption{Temporal profiles at position $x = 0.32$.}
			\label{art_time}
		\end{subfigure}
		\begin{subfigure}[H]{.33\textwidth}  
			\includegraphics[width =\textwidth]{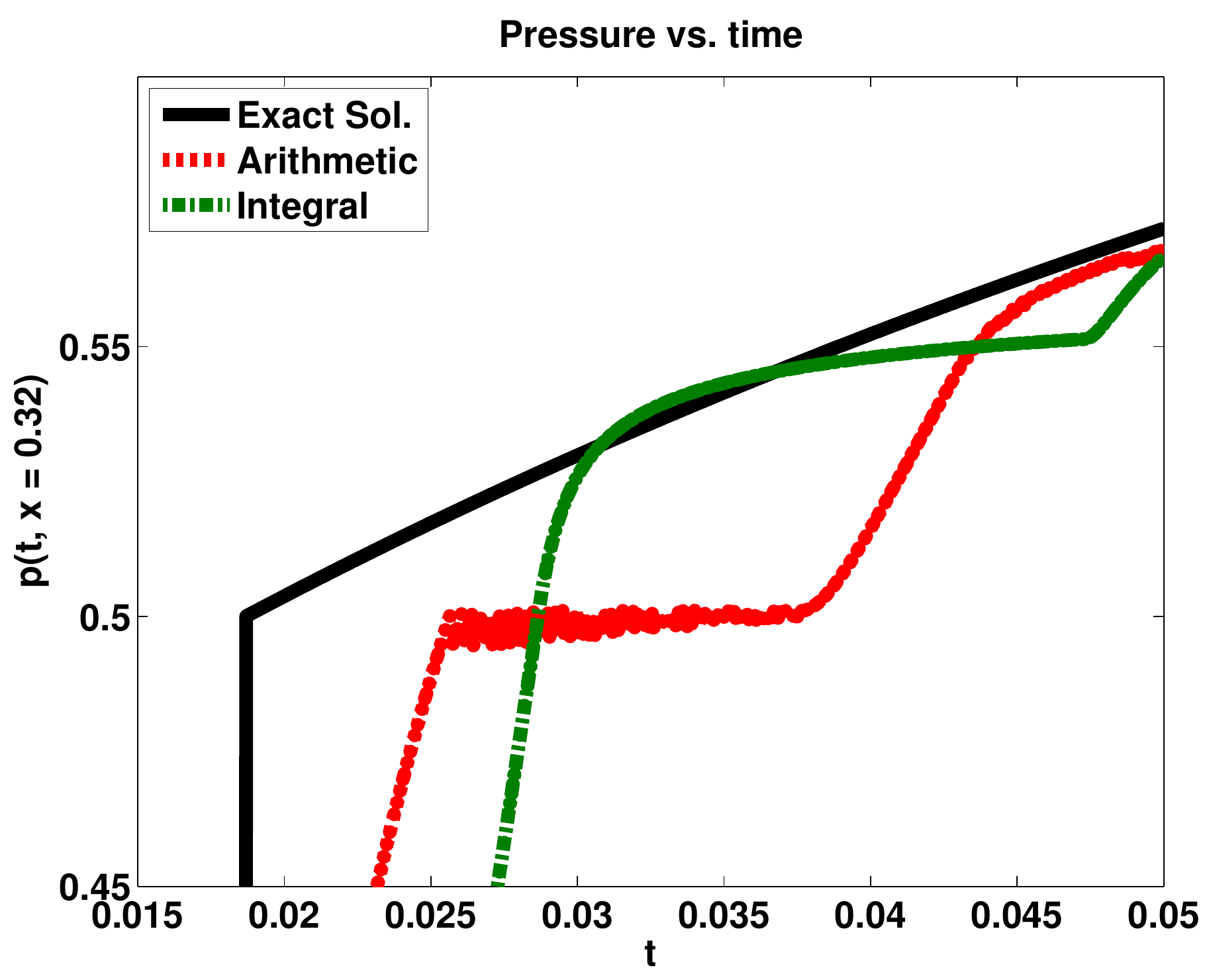}
			\caption{Zoomed in temporal profiles.}
			\label{art_timezoom}
		\end{subfigure}
		\caption{Comparison of various averages, where $k_{\max} = 1$, $k_{\min} = 0$ and $p^* = 0.5$.  The spatial step size $\dx = 0.04$ and the time step size $\dt = \dx^2/32$.}
		\label{fig:intro}
		\end{figure}

\begin{figure}[H]
		\center
		\begin{subfigure}[H]{0.33\textwidth}  
			\includegraphics[width =\textwidth]{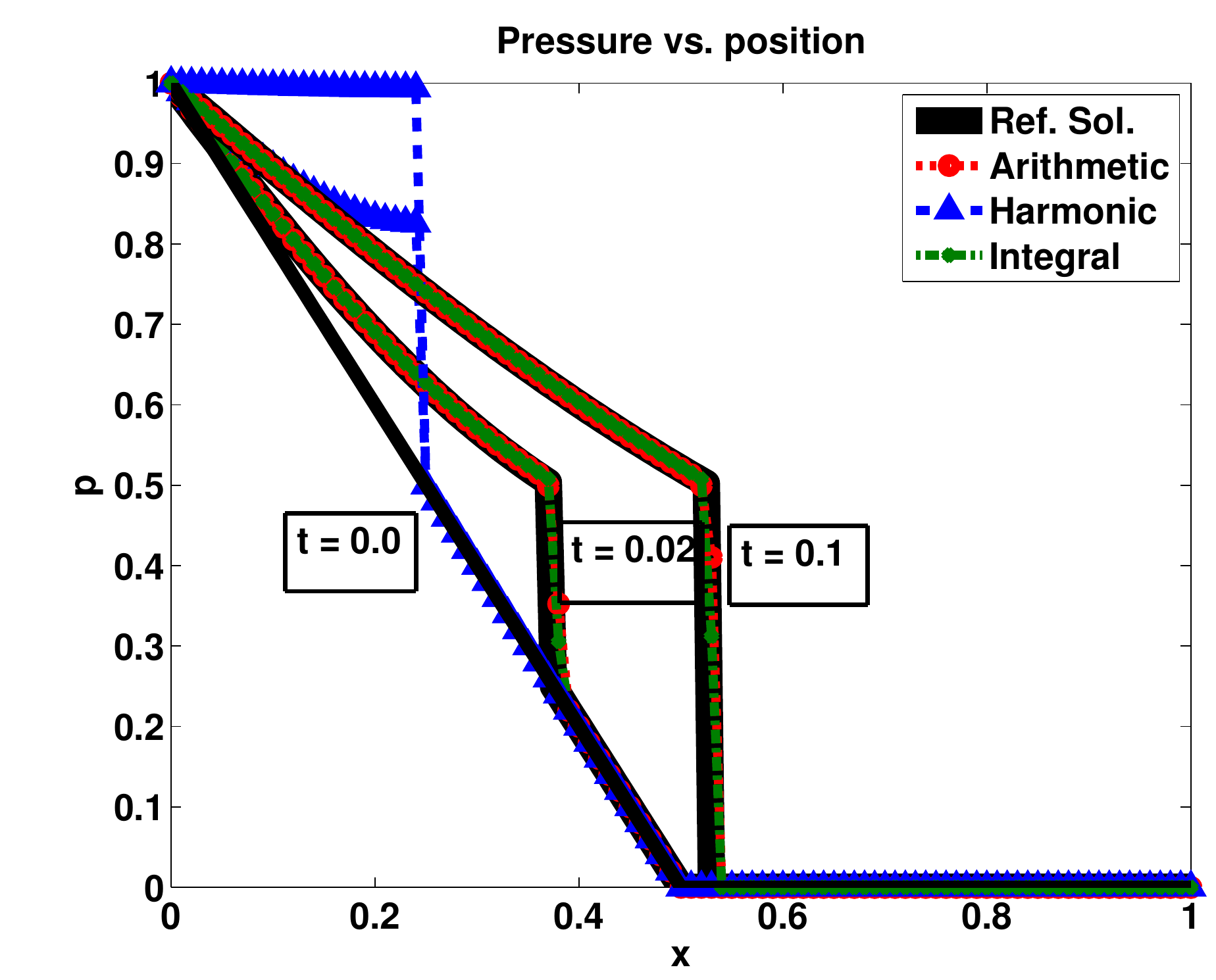}
			\caption{Spatial profiles}
		\end{subfigure}
		\begin{subfigure}[H]{0.33\textwidth}  
			\includegraphics[width =\textwidth]{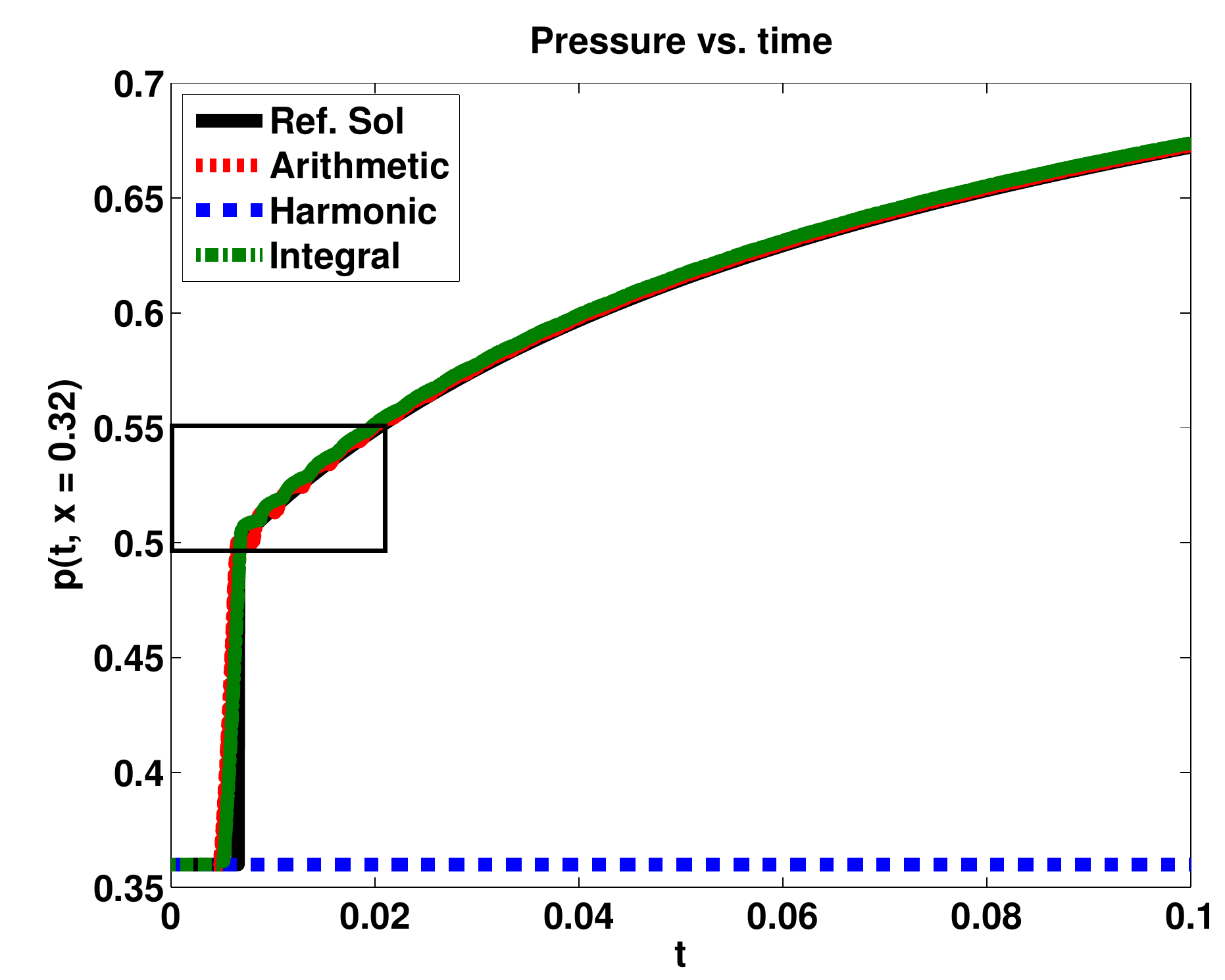}
			\caption{Temporal profiles at position $x = 0.32$.}
		\end{subfigure}
		\begin{subfigure}[H]{0.33\textwidth}  
			\includegraphics[width =\textwidth]{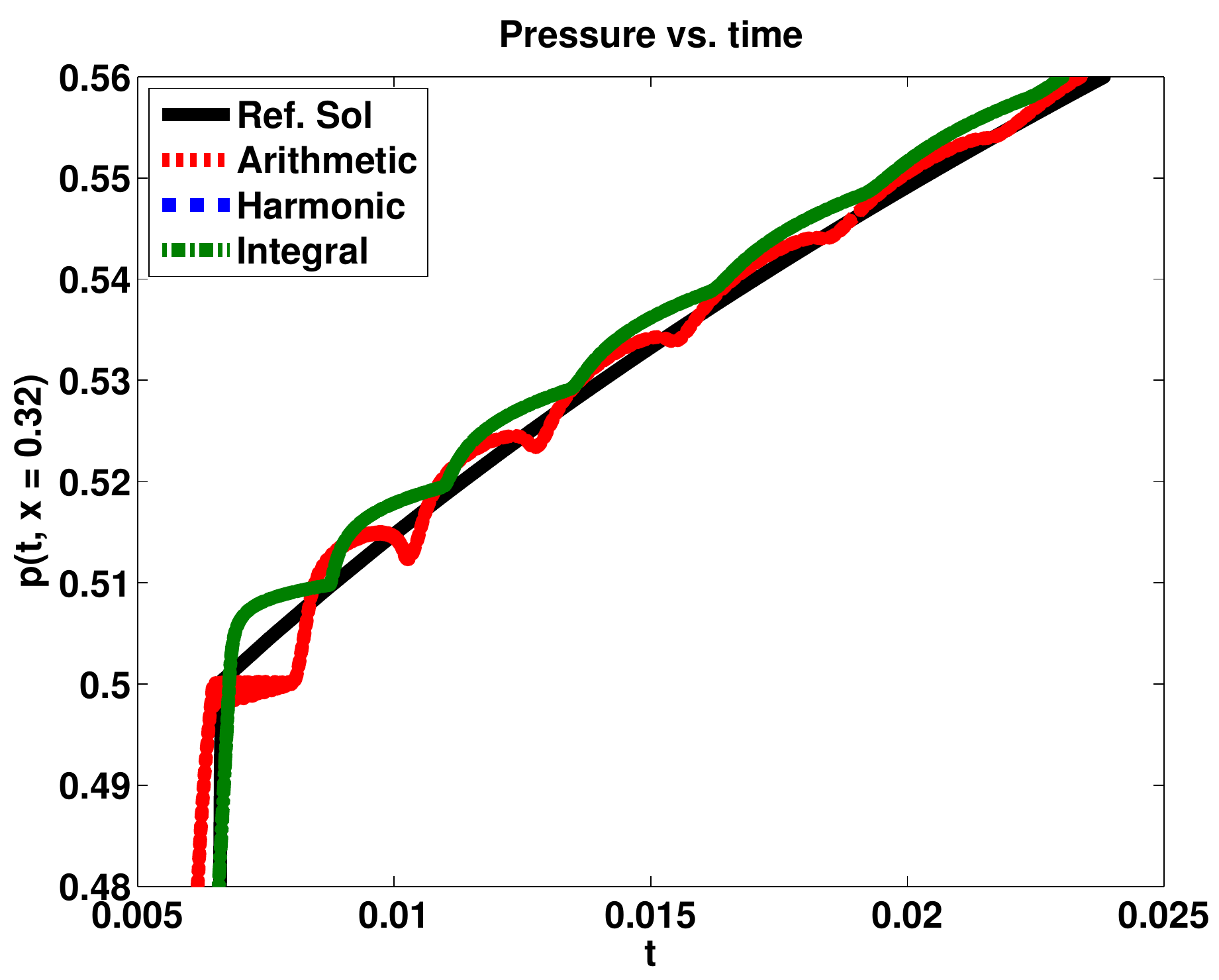}
			\caption{Zoomed in temporal profiles.}
		\end{subfigure}
			\caption{Comparison of various averages, where $k_{\max} = 1$, $k_{\min} = 0$ and $p^* = 0.5$ for a waiting time phenomenon example.  
The piecewise linear initial condition self-sharpens until it is sharp enough for the support interface position to move rightward at later times past the initial interface position at $x = 0.5$.  The spatial step size $\dx = 0.01$ and the time step size $\dt = \dx^2/32$.}
			\label{wait_time}
\end{figure}

The degeneracy, self-sharpening and nonlinearity of the GPME pose interesting numerical challenges.  The numerical approaches used by practitioners vary based on the field.  In the porous media communities, central flux-based finite volume methods are widely employed to solve the related variable coefficients problem
$p_t  = \nabla \cdot (k(x) \nabla p)$.
The coefficient $k(x)$ is defined at the cell-centers, and harmonic and arithmetic averaging are commonly used to compute the coefficient, which may be discontinuous, at the cell interface.  Harmonic averaging is often preferred in these problems because it leads to more physical solutions.  This common finite volume averaged-based approach has also been extended to the nonlinear GPME in Eqn. \eqref{eq:GPME}.  For our case, the selection of the $k(p)$ averaging is not as straightforward as in the variable coefficient case.     
After carefully comparing arithmetic to harmonic averaging, \citet{lipnikov2016} prefer arithmetic averaging over harmonic for the continuous PME with near-zero $k(p) = p^m$.  An integral average is developed in \citet{vandermeer2016} for the discontinuous GPME in Eqn. \eqref{eq:discont_k}.  With these average-based approaches, it is difficult to satisfy the jump condition in Eqn. \eqref{eq:RH}. Violation of this condition manifests itself in locking and spurious temporal oscillations, as illustrated in Figures \ref{fig:intro} and \ref{wait_time} and further discussed in Section \ref{static_avg}.   
In other fields, adaptive and moving mesh approaches have been applied to this problem also with increased refinement near the shock.  In \cite{ngo2016}, for example, a finite element moving mesh method for the PME is developed.  Adaptive Mesh Refinement (AMR) will be discussed further in Section \ref{static_avg}. 
 
 Variational particle schemes have been developed for the PME in \cite{wilkening10, otto01,villani03}, where a variational principle is defined and steepest descent is applied to optimize the corresponding energy function.  Variational principles have also been defined for the Stefan problem in \cite{rossi04} and a similar gradient flow method can be applied for the discontinuous GPME.
 
In the crystal application communities, numerical work has been done on solving the Stefan problem, as detailed in \cite{berger79, kharab86, Rizwan98, caldwell2002} and the references therein.  In \cite{rose90, rose93}, an enthalpy scheme for Stefan problems is introduced.  This method consists of a compact finite difference stencil with an implicit temporal scheme.  The shock position is not incorporated into the scheme, and the error is observed to be concentrated at the moving interface.  \citet{brattkus92} show that finite difference and finite element methods that do not take special care at the interface result in a significant error of $\mathcal{O}(\sqrt \dx)$.  They propose an effective method based on an integral equation formulation that 
 requires an additional integral evaluation.

In two-phase flow applications of the Stefan problem, level set methods \cite{sethian88} could also be useful in tracking the evolving interface.  Thus far, they have not been directly applied to Eqn. \eqref{eq:stef_sp}. Related problems with Stefan boundary conditions are discussed in \citet{sethian92} and \citet{osher97}.  More recent approaches using the level set method in the solidification community are detailed in \cite{kim00,tan07}.  An advantage of the level set approach over traditional front-tracking approaches \cite{rawahi02, li03} is that the interface is treated implicitly.  In doing so, the level set method can be more easily extended to higher dimensions, and is robust under complex topologies.  
Numerical methods combining the level set approach with the Extended Finite Element Method (X-FEM) have also been developed to solve these Stefan problems \cite{bernauer12, salvatori09}.
 A modified Stefan problem is also solved in \citet{zhao2016} using a phase-field method.  Phase-field methods converge as the interface thickness parameter tends to zero.  The use of phase-field approaches in solidification and the effect of the interface thickness parameter have also been studied in \cite{yang2015, han2013}.  While these papers in the solidification and phase change communities provide approaches to solve the Stefan problem, they do not give details on understanding why the numerical artifacts occur with finite volume average-based approaches.  This leads us to the main goal of this paper.

\subsection{Main Goal of this Paper}
The main goal of this paper is to shed light on the origin of the numerical problems reported in the literature for finite volume averaged-based approaches of the discontinuous GPME.   
Before we discuss the cause of the artifacts, we first introduce an alternate numerical method 
for Eqn. \eqref{eq:stef_sp} that does not exhibit the numerical errors observed in the literature.  This newly developed Shock-Based Averaging Method (SAM) is discussed in Section \ref{sam_exact}.  We discuss this first because it shows that incorporation of the shock position in the scheme is key, and also helps with the understanding of the artifacts in later sections.  Section \ref{exact_shock} discusses the derivation of SAM when the shock position is known.  The shock position can also be approximated using the jump condition, which is detailed in Section \ref{shock_speed}.  In Section \ref{static_avg}, we cast SAM in the finite volume framework to help identify what has been lacking in the other approaches.  The main issue is that these discretizations do not contain enough information about the shock.

\section{Proposed Numerical Method: Shock-Based Averaging Method (SAM)}
\label{sam_exact}
Due to the discontinuities in $k(p)$ and $p$, we adopt a finite volume approach to the integral form of the governing equation: 
\begin{equation}
	\begin{aligned}
	p_t &= \nabla \cdot (k(p) \nabla p), \hspace{.25cm} \text{where}
\\
	k(p)& = 
	\begin{dcases}
  		k_{\max},           & p \ge p^* \\
    		k_{\min},              &p^* > p,
	\end{dcases}
\label{eq:prob_def}
\end{aligned}
\end{equation}
where $p^*$ is the left limiting value.
We first define a finite volume grid with cell-centers $x_j$ for $j = 1,\dots,N + 1$. The boundaries of the domain are at $x_1$ and $x_{N+1}$, respectively and the corresponding unknowns $p$ are fixed by Dirichlet boundary conditions.  The remaining $N\!-\!1$ degrees of freedom $p_j$, and the corresponding coefficients $k_j$, are defined at the nodes $x_j$ for $j = 2,\dots,N$.  We define control volumes $CV_j$ with width $\dx_j$.  
The cell faces $x_{j+1/2}$ of each $CV_j$ are at a distance of $\dx_j/2$ from the cell-centers $x_j$.
We assume that the solution is monotone and non-increasing, and that the shock is located between $x_i$ and $x_{i+1}$, such that $ p_i \ge p^* \ge p_{i+1}$.  

The semi-discrete numerical discretization for $CV_j$ with volume $\dx_j$ is given by
\begin{equation}	
	\dx_j\frac{dp_j}{dt}= F_j^- - F_j^+,
	\label{eq:fv_discret}
\end{equation}
where $F_j^-$ 
represents the in-flux and $F_j^+$ 
represents the out-flux of $CV_j$.    
The numerical fluxes  
\begin{equation}
	 F_{j}^+ = -k_{j+1/2}\frac{p_{j+1}-p_j}{\dx}, \hspace{.25cm} j \ne i,
	\label{eq:stand_outflux}
\end{equation}
and
\begin{equation}
	F_{j}^- = -k_{j-1/2}\frac{p_{j}-p_{j-1}}{\dx}, \hspace{.25cm} j \ne i+1,
	\label{eq:stand_influx}
\end{equation}
are defined at all faces away from the shock cell with the standard two-point flux approximation, where $\dx \equiv \dx_j$.  The coefficient $k_{j+1/2}$ at the cell face represents a local average of its neighboring coefficients.  Then, for any two-neighbor average away from the shock cell, the coefficient is constant and given by
\begin{equation} 
k_{j+1/2} =
	\begin{dcases}
  		 k_{\max},           & 1 \le j \le i - 1,   \\
		  k_{\min}, & i+1 \le j \le N.
	\end{dcases}
	\label{eq:gen_coeff}
\end{equation}
Analogous definitions are used for $j\!-\!1/2$.
Because of the discontinuity, Eqns. \eqref{eq:stand_outflux}-\eqref{eq:stand_influx} cannot be used for the flux at the cell interface $x_{i+1/2}$.

\subsection{Formulation of the Fluxes near the Shock}
\label{aux_grid}
\begin{figure}[H]
		\center 
		\includegraphics[width =.49\textwidth]{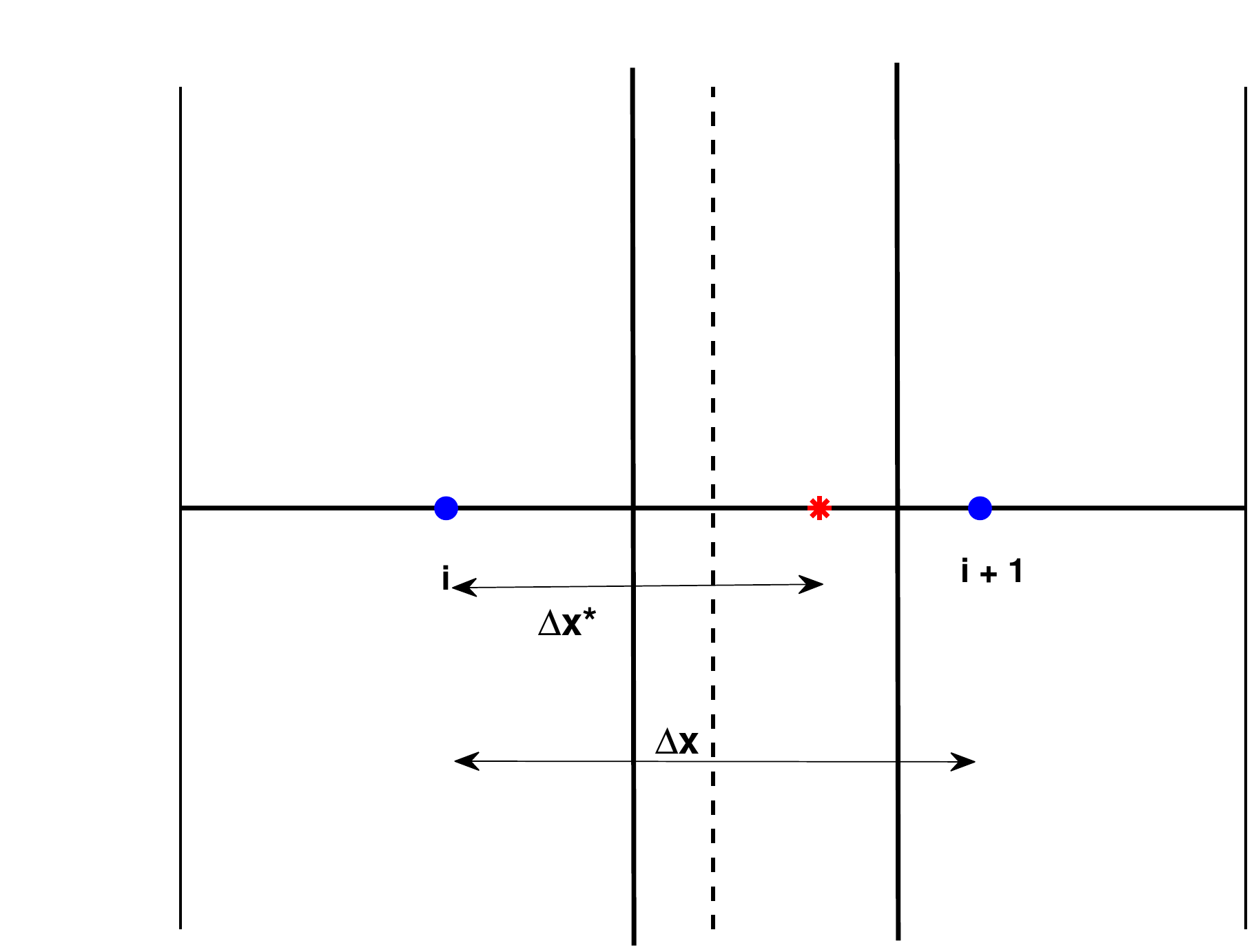}
		\caption{Illustration of the additional control volume $CV_*$ around the shock position $x^*(t).$  The grid points $x_i$ and $x_{i+1}$ are depicted by the blue dots, and $x^*(t)$ is depicted by the red star.   The dashed line represents the deleted cell face between $CV_i$ and $CV_{i+1}$.  The new cell faces are located at a distance $\dx^*(t)/2$ from $x_i$ and $[\dx - \dx^*(t)]/2$ from $x_{i+1}$.}
		\label{fig:FV}
\end{figure}

To estimate the fluxes out of $CV_i$ and into $CV_{i+1}$, we borrow ideas from hyperbolic systems \cite{rao09} and place a control volume around the discontinuity.  
Figure \ref{fig:FV} shows the auxiliary finite volume grid with the additional control volume $CV_*$ around the physical shock position $x^*(t)$, where $p(x^*(t),t) \equiv p^*$ is defined.      
Figure \ref{fig:FV} also shows $\dx^*(t)$ as the distance between $x_i$ and $x^*(t)$, where $0 \le \dx^*(t) \le \dx$.  We then remove the cell face at position $x_i + \dx/2$ and add the cell faces of $CV_*$.  
The new cell faces are chosen to be centrally located at a distance of 
\[
	\frac{|x^*(t) - x_i| }{2} \equiv \frac{\dx^*}{2},
\] from $x_i$ and $x^*(t)$, and a distance of 
\[
	\frac{|x^*(t) - x_{i+1}|} { 2} \equiv \frac{\dx - \dx^*}{2},
\]
 from $x^*(t)$ and $x_{i+1}$.  
 We have effectively added an additional degree of freedom $p^*$ at $x^*(t)$.  In Eqn. \eqref{eq:prob_def}, $p^*$ is known, and it does not need to be computed.  We present the approach in its general form for future extensions.

 With the additional control volume, we can now define the necessary out-flux $F_i^+$ of $CV_i$ and in-flux $F_{i+1}^-$ of $CV_{i+1}$. 
 The coefficient at each of the new cell faces is known.  By the monotonicity of $p$, for any $p(x,t)$ to the left of the interface, $k(p(x,t)) = k_{\max}$.  Similarly, for any $p(x,t)$ to the right of the interface, $k(p(x,t)) = k_{\min}$.  
For this piecewise constant coefficient, any two-neighbor average on this auxiliary grid gives the same coefficient value at the faces.  Based on these known coefficient values at the faces and the face-centered grid, the second-order fluxes are given by 
\begin{equation}
	F_{i}^+ = -k_{\max}\frac{p^*- p_i}{\dx^*}, 
	\label{eq:fluxi}
\end{equation}
at the cell face between $x_i$ and $x^*$
and
\begin{equation}
F_{i+1}^- = -k_{\min}\frac{p_{i+1}-p^*}{\dx-\dx^*}, 
\label{eq:fluxiplus}
\end{equation} 
 at the cell face between $x^*$ and $x_{i+1}$. 
 The resulting scheme is conservative because $F_i^+$ is equal to the in-flux $F_*^-$ of $CV_*$, by definition and similarly $F_{i+1}^-$ is equal to the out-flux $F_*^+$ of $CV_*$.  The fluxes are then substituted into the semi-discrete equation \eqref{eq:fv_discret}, where the cell volumes are now given by $\dx_i = (\dx+ \dx^*)/2$ and $\dx_{i+1} = \dx- \dx^*/2$, respectively.
 
 
For $k_{\min} = 0$, we can relate the expression for $F_{i}^+$ in Eqn. \eqref{eq:fluxi} to the expression for $\Phi(p) = k_{\max}(p-p^*)_+$ in Eqn. \eqref{eq:stef_sp} for the Stefan problem.  The analytical flux $F(p_L)$ is given by $-\nabla \Phi(p_L)$ in Eqn. \eqref{eq:anal_flux}.  The numerical flux $F_{i}^+$ can be interpreted as approximating this gradient with an upwind discretization as 
\[
	-\frac{\Phi(p^*)-\Phi(p_i)}{\dx^*} = -k_{\max}\frac{p^*-p_i}{\dx^*} = F_{i}^+.
\]
We see that the above formula is a first order approximation to the flux at $x^*(t)$.  From the jump condition in Eqn. \eqref{eq:RH_pre}, the velocity can then be computed.

Round-off errors can arise when the shock location approaches the grid points $x_i$ and $x_{i+1}$, and the denominators in Eqn. \eqref{eq:fluxi} and Eqn. \eqref{eq:fluxiplus} approach zero.  To avoid problems, we specify a tolerance of $\eps$ that is proportional to $\dx$.  If $\dx^* \le \eps$, we neglect the leftmost portion near $x_i$, and set $k_{i+1/2} = k_{\min}$ for $CV_i$.  Similarly, if $\dx - \dx^* \le \eps$, we ignore the rightmost portion near $x_{i+1}$, and set $k_{i+1/2} = k_{\max}$ for $CV_{i+1}$.  We then use these coefficients in the standard two-point flux approximation in Eqns. \eqref{eq:stand_outflux} and \eqref{eq:stand_influx}.  

We will refer to the above approach as the Shock Based-Averaging Method (SAM).

\section{SAM Numerical Results: Exact Shock Location}
\label{exact_shock}
\begin{figure}[H]
		\center
		\includegraphics[width =0.49\textwidth]{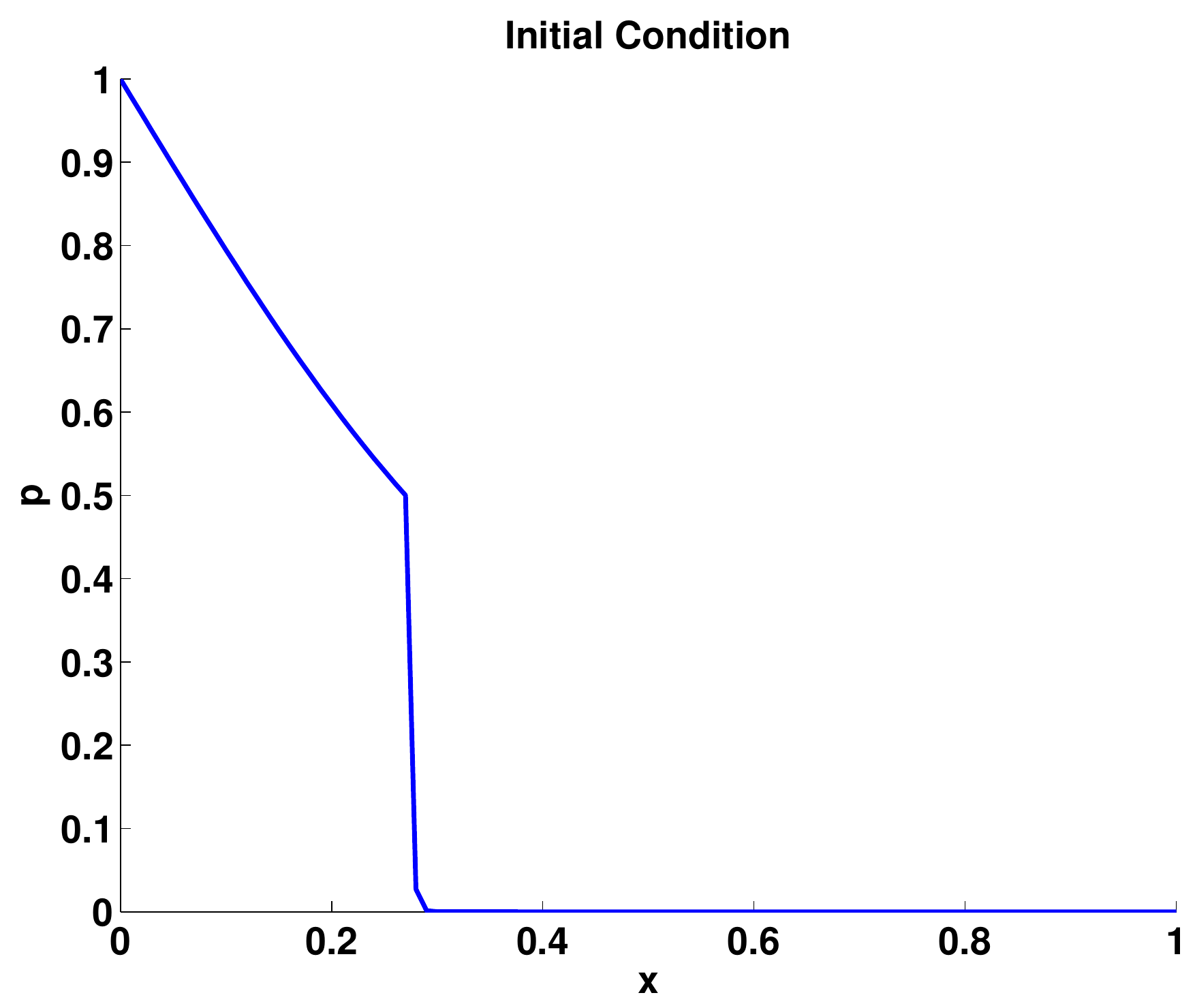}
		\caption{Spatial profile at time $t = 0$.}
		\label{IC}
\end{figure}

In the numerical results presented throughout the paper, the test problem is given by Eqn. \eqref{eq:prob_def} with 
$k_{\max}=1.0$, $k_{\min}=0.0$ and $p^* = 0.5$.  The choice of $k_{\max}$ and $p^*$ is arbitrary, while $k_{\min}$ is set to zero to test the behavior of the numerical method in the degenerate case.  The algorithm also works for arbitrarily small $k_{\min}$. 
The fixed Dirichlet boundary conditions are given by
	\begin{equation}
		\begin{aligned}
			p(0,t) &= 1.0, \hspace{.2cm} \forall t \ge 0, \\
			p(1,t) &= 0.0,\hspace{.2cm} \forall t \ge 0.
		\end{aligned}
		\label{eq:bc}
	\end{equation}
We first present results for the initial condition displayed in Figure \ref{IC}.  This initial condition is formed as the exact solution to Eqn. \eqref{eq:prob_def} at a later time.  The same parameters for $k_{\max}$ and $p^*$, as given above for the test problem, are used to generate the initial condition, whereas $k_{\min} = 0.01$ for some smoothness at the bottom of the front.  The coefficient value $k_{\min} = 0$ throughout the simulations.
For the representative test problem, an analytical solution exists,   
and is used to verify the implementation through detailed convergence studies. 
The derivation of the exact solution is generalized for arbitrary $k_{\max}$ in \ref{exact}.  In this section,  we use the exact shock position to compute $\dx^*(t) \equiv |x^*(t) - x_i|$, where 
\[
	x^*(t) = \alpha \sqrt{t}, \hspace{.5cm} \alpha  = 2\sqrt{k_{\max}}z_1,
\]
and $z_1$ is the solution to the nonlinear equation in Eqn. \eqref{eq:alpha}.  In the following section, we will show how the shock location can be approximated, if it is not available.

\begin{figure}[H]
		\center
		\begin{subfigure}[H]{.49\textwidth}  
			\includegraphics[width =\textwidth]{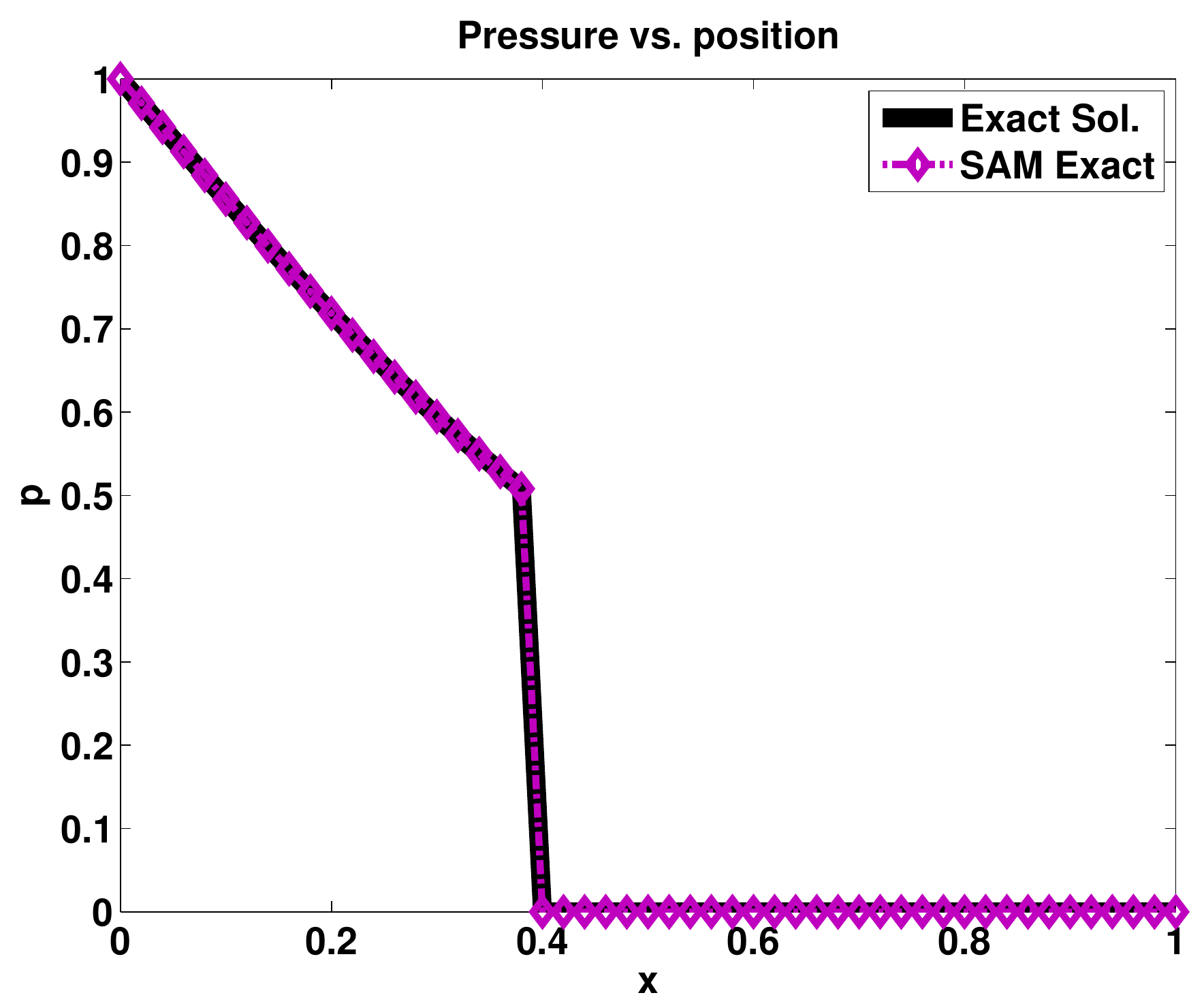}
			\caption{$N = 50$ grid points}
		\end{subfigure}
		\begin{subfigure}[H]{.49\textwidth}  
			\includegraphics[width =\textwidth]{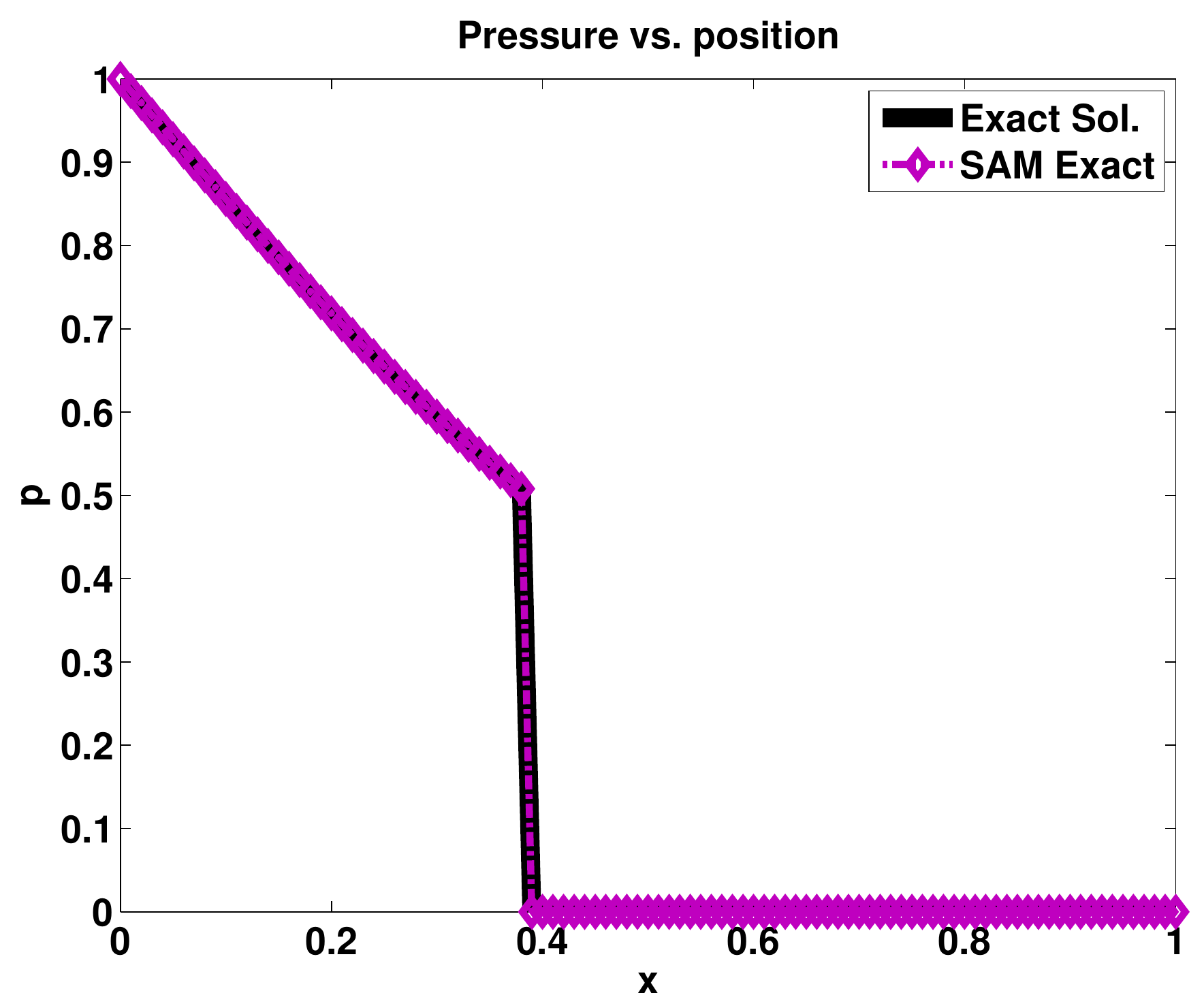}
			\caption{$N = 100$ grid points}
		\end{subfigure}
		\caption{Spatial profiles at time $t = 0.05$ with $\dt  = \dx^2/32$.}
		\label{fig:sam_exact_space}
\end{figure}

After numerical stability and accuracy tests, the time step is selected to be $\dt  = \dx^2/32$ for the explicit Forward Euler method.  The discontinuity in the coefficient and the degeneracy make the time step criterion be more restrictive than it is for classical parabolic equations.  If stability is the only interest, and not accuracy, the coefficient in the time step can be increased, but the time step must still be on the same order of $\mathcal{O}(\dx^2)$.
\begin{figure}[H]
	\center
	\begin{subfigure}[H]{0.49\textwidth}  
			\includegraphics[width =\textwidth]{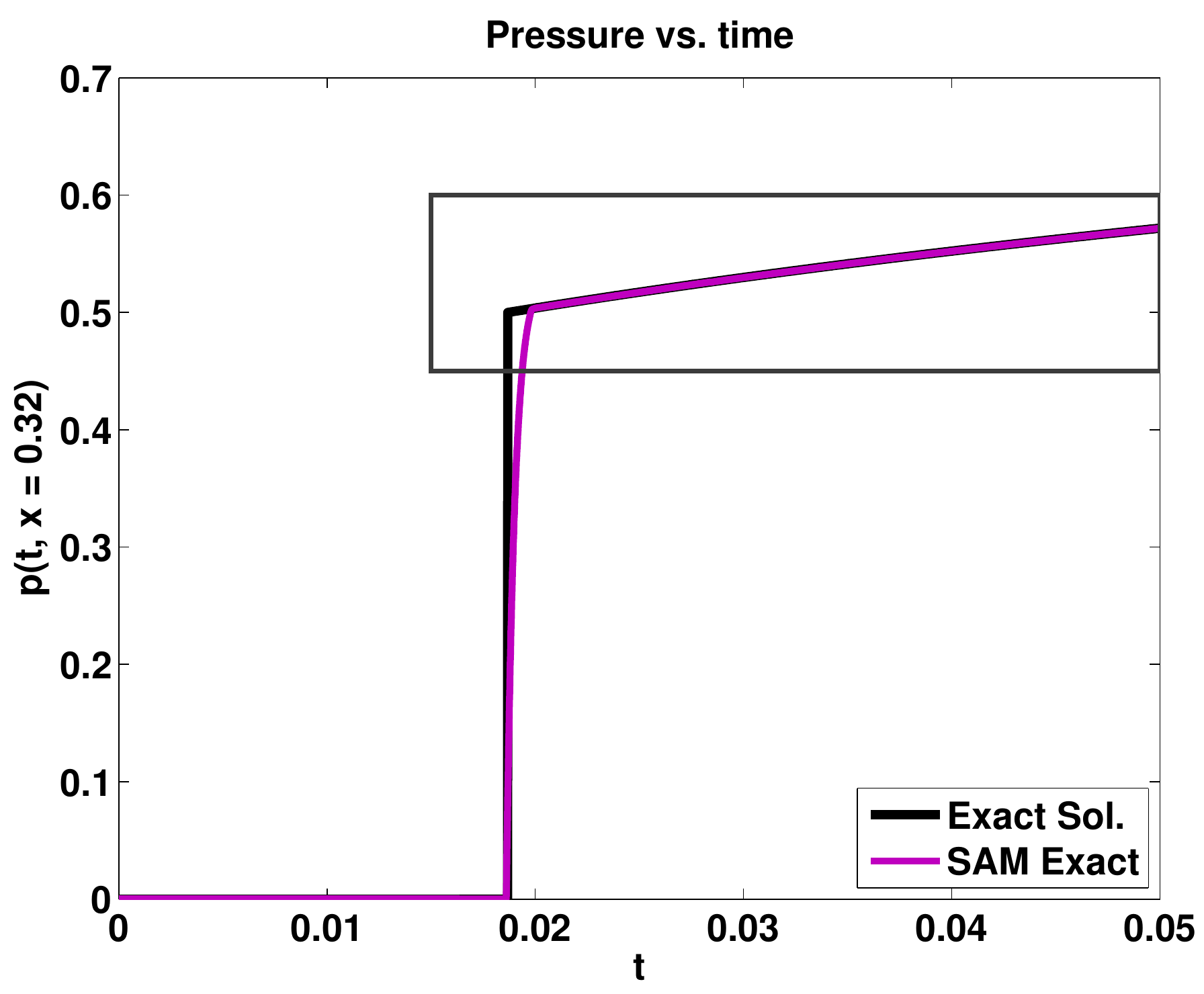}
		\end{subfigure}
		\begin{subfigure}[H]{0.49\textwidth}  
			\includegraphics[width =\textwidth]{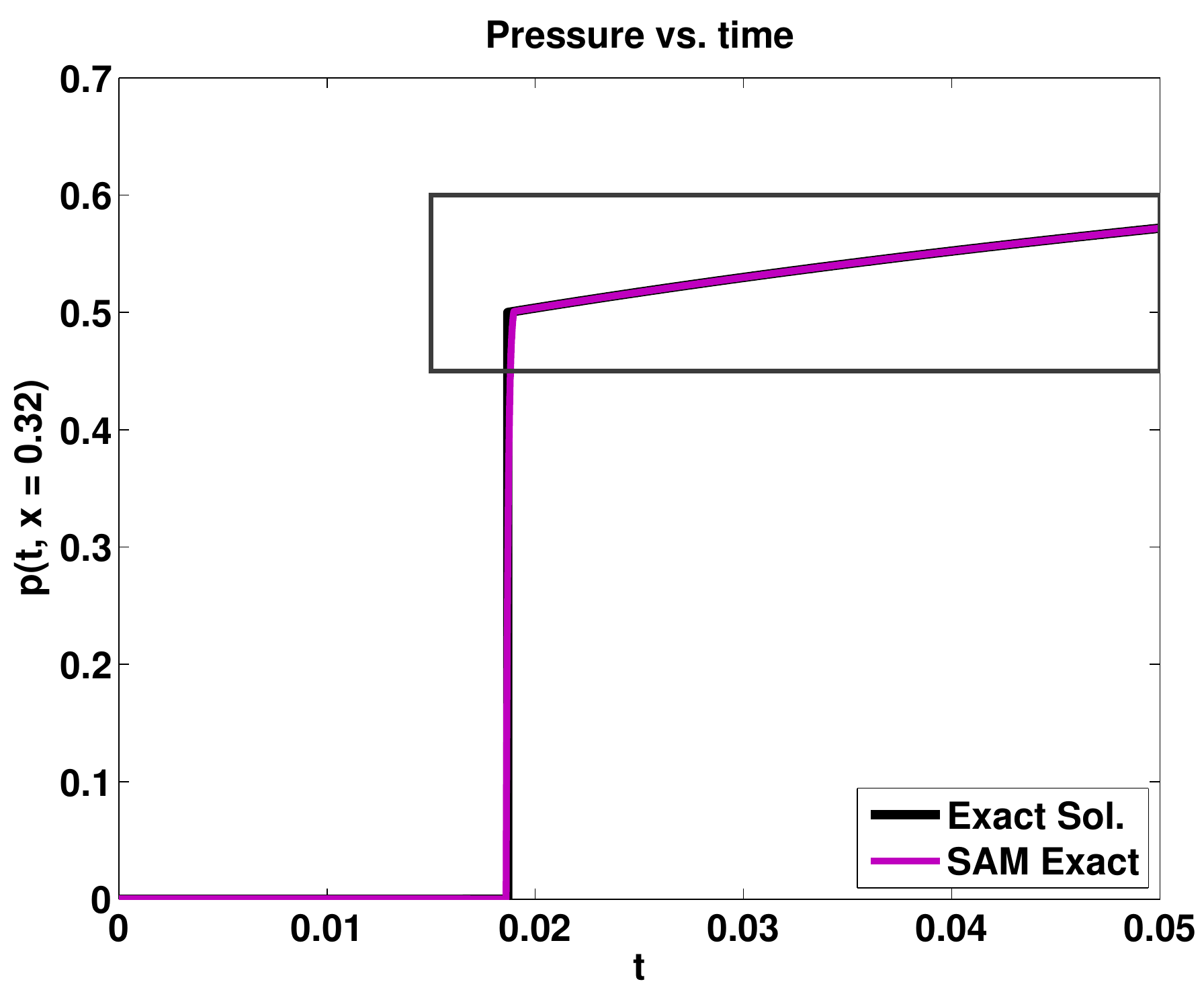}
		\end{subfigure}
		\begin{subfigure}[H]{0.49\textwidth}  
			\includegraphics[width =\textwidth]{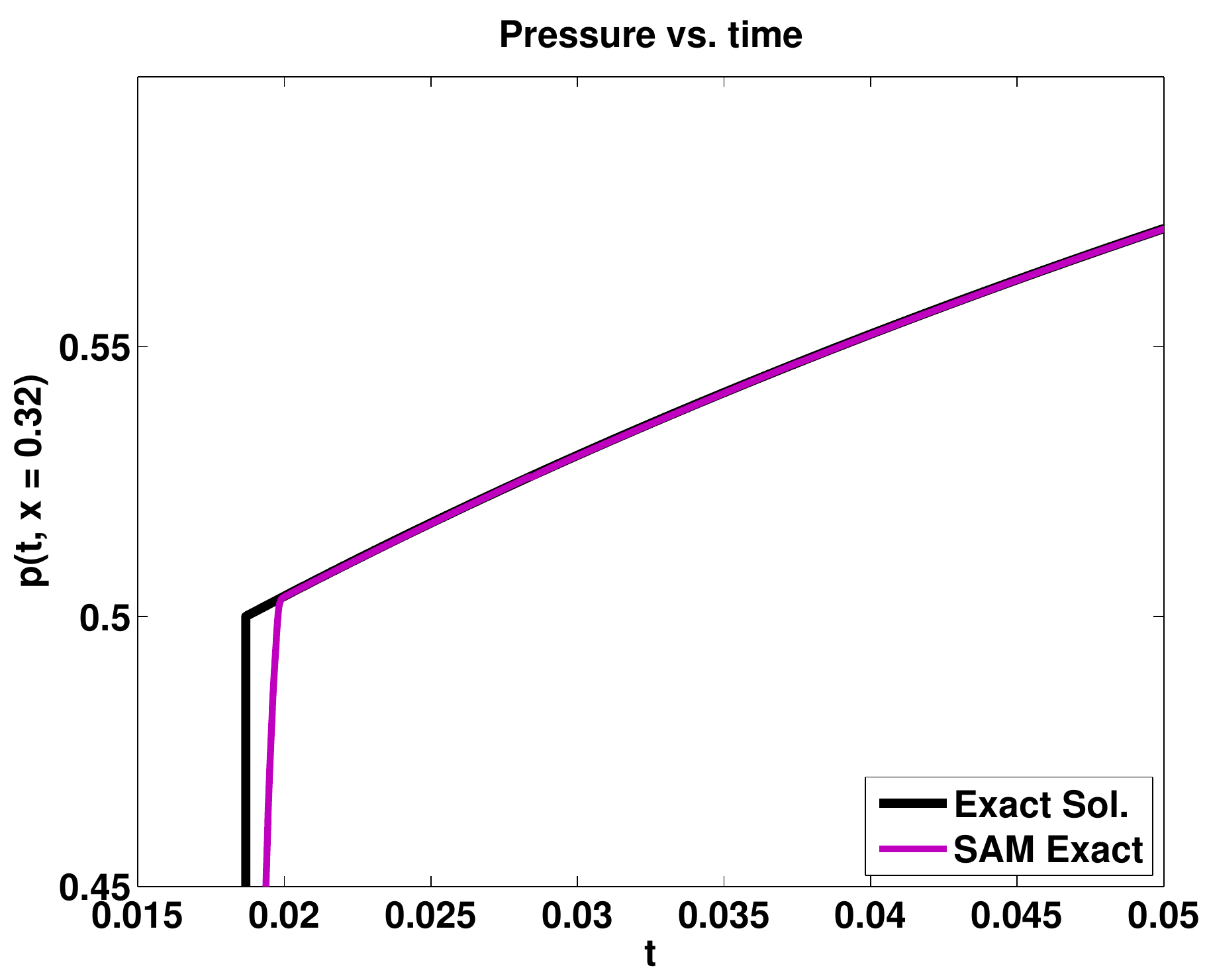}
			\caption{$N = 50$ grid points}
		\end{subfigure}
		\begin{subfigure}[H]{0.49\textwidth}  
			\includegraphics[width =\textwidth]{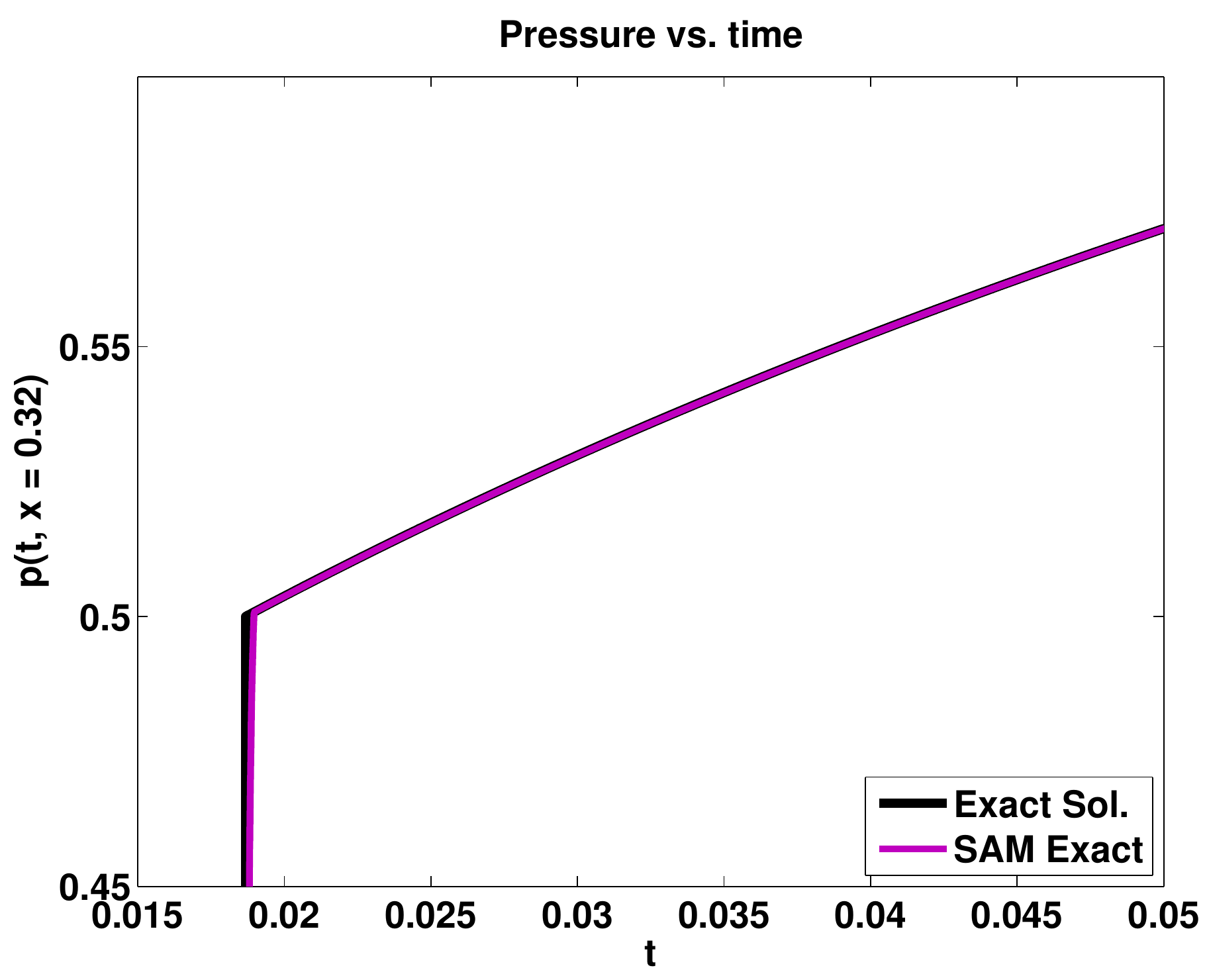}
			\caption{$N = 100$ grid points}
		\end{subfigure}
		\caption {Temporal profiles at position $x = 0.32$ with $\dt  = \dx^2/32$.}
		\label{fig:sam_exact_time}
\end{figure}
In Figure \ref{fig:sam_exact_time}, the numerical solution is plotted at the arbitrary position $x = 0.32$ over time.  By self-similarity of the solution \cite{vazquez2007, maddix_pme, ngo2016},  the temporal plots have the same profile for any $x$-coordinate that the front has passed through.  Figure \ref{fig:sam_exact_time} illustrates that the numerical solution with SAM has an accurate and non-oscillatory temporal profile, even on a coarse grid.  Figure \ref{fig:sam_exact_space} reveals that the shock location is accurate in space for this moving front problem.  The self-sharpening nature of the GPME \cite{maddix_pme, ngo2016, vazquez2007} is depicted, by the smooth lower corner in the initial condition in Figure \ref{IC} evolving into a sharp corner in Figure \ref{fig:sam_exact_space}.  Figure \ref{fig:sam_exact_space} also illustrates the sharp capture of the shock.  


\section{SAM Numerical Results: Approximate Shock Location}
\label{shock_speed}
The method in the prior section can be extended to more general GPME problems, where the exact shock position is not known.
The finite speed of propagation property and theoretical speed of the front for the GPME, discussed in Section \ref{rh_cond}, can be utilized to numerically approximate the shock location. 
Again the fluxes in Eqns. \eqref{eq:fluxi} and \eqref{eq:fluxiplus} are used with the difference that $\dx^*(t)$ is numerically calculated.  In this section, we show that this approximation does not introduce numerical artifacts, such as temporal oscillations.  The convergence results are provided in \ref{app:conv}. 

 We discretize $V$ in Eqn. \eqref{eq:front_speed} with upwinding for the derivative as 
\begin{equation}
	\hat{V} = -\frac{p_i - p_{i-1}}{\dx p_i}, 
	\label{eq:num_speed}
\end{equation}
 where $i$ is the index, such that $p_{i} \ge p^*  \ge p_{i+1}$.  Eqn. \eqref{eq:num_speed} can be interpreted as an approximation of the jump condition in Eqn. \eqref{eq:RH_pre}, where $p_L$ is approximated by $p_i$ and $p_R = 0$.  For problems where $p_R$ is not initially zero, such as the waiting time problem in Section \ref{lit}, we can also use a discrete approximation to Eqn. \eqref{eq:RH}.  Care must be taken numerically when approximating $p_R$ in Eqn. \eqref{eq:RH}, since $p_{i+1}$ is not guaranteed to be zero: all that is guaranteed is that $0 \le p_{i+1} \le p^*$.  In this case, we use $p_{i+2}$ to approximate $p_R$.
 
To obtain the numerical shock position, a simple time integration is implemented.  We let $\xi^n$ represent the approximate $x^*(t)$ at time step $n \ge 0$ and substitute 
\begin{equation}
	\dx^*(t) \approx \xi^n - x_i,
	\label{eq:deltax_star}
\end{equation}
into the expressions for the fluxes in Eqns. \eqref{eq:fluxi} and \eqref{eq:fluxiplus}.  We integrate the approximate shock speed in Eqn. \eqref{eq:num_speed} in time, using the following update
\begin{equation}
			\xi^{n+1} = \xi^n + \dt \hat{V}.
			\label{eq:xi}
\end{equation}
We assume that the initial position $\xi^0$ is known from the problem definition or can also be approximated. 

We can also approximate $\dx^*(t)$ using a level set method \cite{sethian88}. 
In the level set method, the interface $x^*(t)$ is represented by the zero level set $\{x \hspace{.1cm}	 | \, \phi(x,t) = 0\}$ of a signed distance function $\phi(x,t)$.  
Using the velocity $V$ of the interface in Eqn. \eqref{eq:front_speed}, the evolution equation can be written in terms of $\phi(x,t)$ as
\begin{equation}
	\phi_t + V \cdot \nabla \phi = 0.
	\label{eq:levelset}
\end{equation}
The velocity can be extended to the entire domain or in a neighborhood of the interface by using velocity extension methods based on the Fast Marching Method \cite{chopp09}.  The level set equation \eqref{eq:levelset} can be solved numerically using Forward Euler in time and upwinding for the gradient, where $V$ is discretized as $\hat{V}$ in Eqn. \eqref{eq:num_speed}.  Since $\dx^*(t)$ is defined as the distance from the shock or interface $x^*(t)$ to the grid point $x_i$, $\dx^*(t)$ can be calculated by $|\phi(x_i)|$.   Solving Eqn. \eqref{eq:levelset} for $\phi_i^n$ and setting $\dx^*(t) \approx |\phi_i^n|$ gives the same results as solving Eqns. \eqref{eq:deltax_star}-\eqref{eq:xi}.  

\begin{figure}[H]
		\center
		\begin{subfigure}[H]{.49\textwidth}  
			\includegraphics[width =\textwidth]{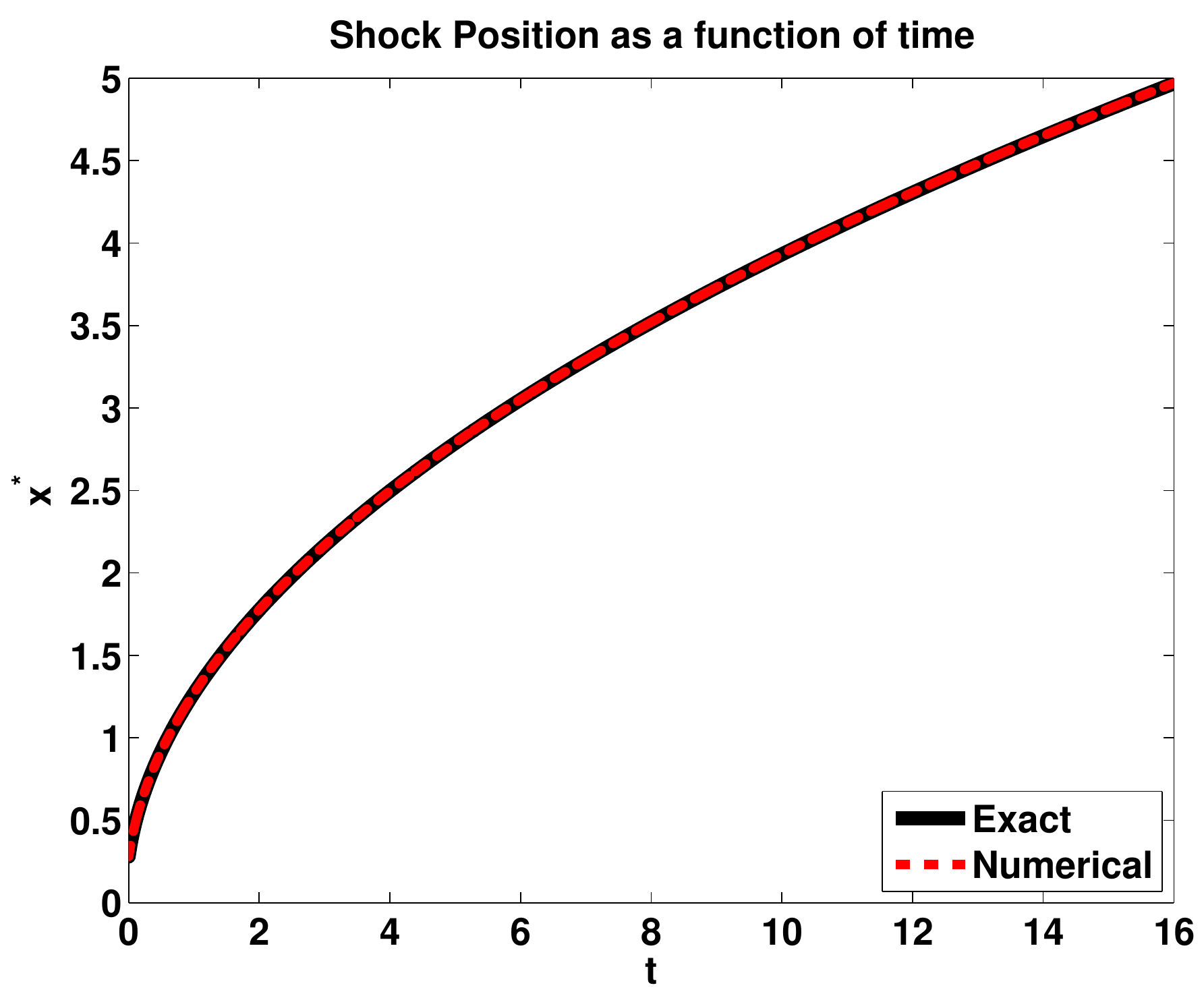}
			\caption{Shock position}
			\label{fig:shock_pos}
		\end{subfigure}
		\begin{subfigure}[H]{.49\textwidth}  
			\includegraphics[width =\textwidth]{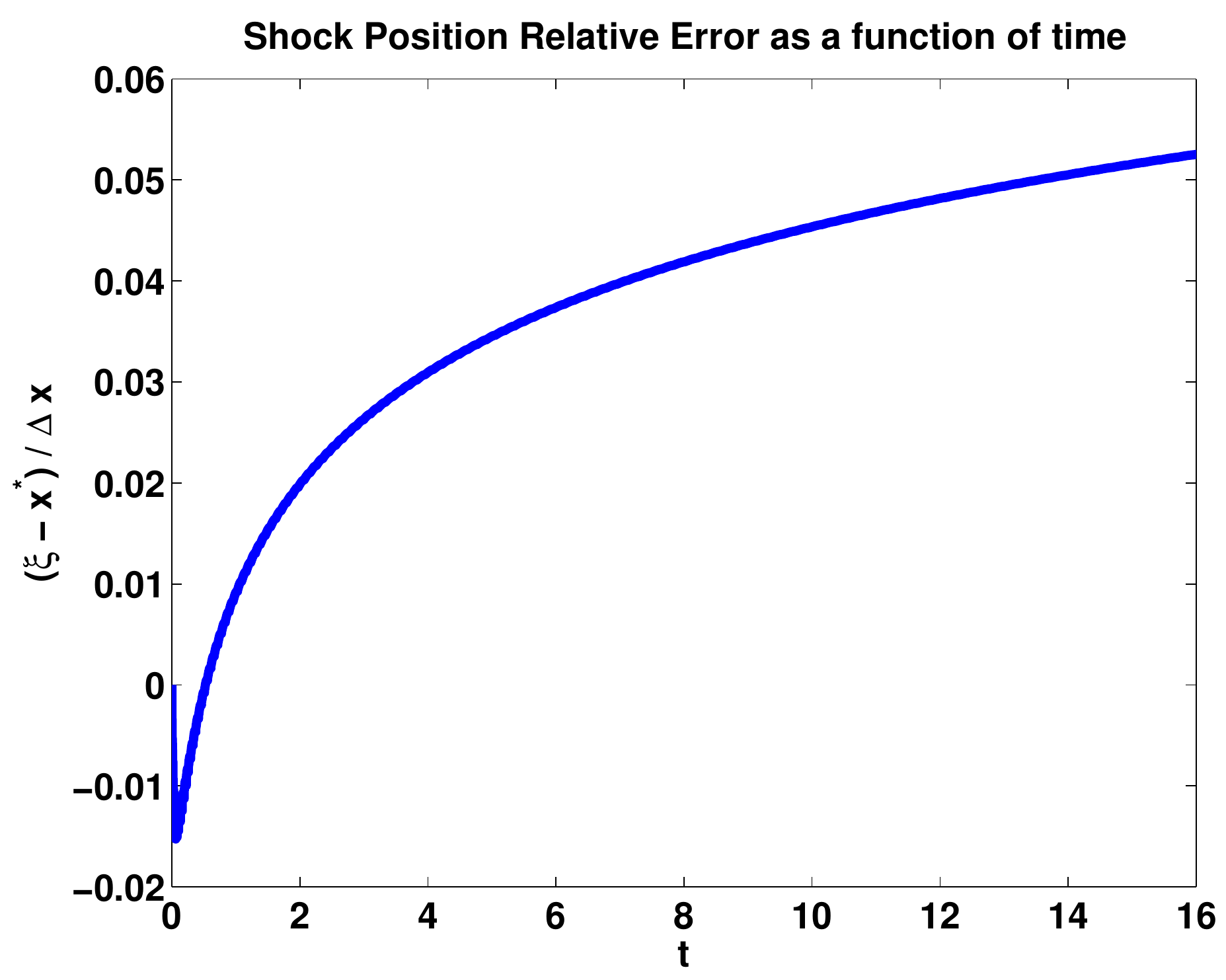}
			\caption{Shock position error relative to $\dx$}
			\label{fig:shock_pos_err}
		\end{subfigure}
		\caption{Comparison of the temporal profiles of the shock position and relative error for $\dx = 0.04$ and $\dt  = \dx^2/32$.}
		\end{figure}

Figures \ref{fig:shock_pos} and  \ref{fig:shock_pos_err} illustrate the long-time behavior of the numerical shock position evolution.  Figure \ref{fig:shock_pos} displays that the numerical shock position aligns with the exact shock position as a function of time, and that the shock position evolution is accurately captured without oscillations.  Both plots verify the numerical implementation, and show that there is no significant accuracy loss in estimating the shock position using Eqn. \eqref{eq:xi}.

\begin{figure}[H]
		\center
		\begin{subfigure}[H]{.49\textwidth}  
			\includegraphics[width =\textwidth]{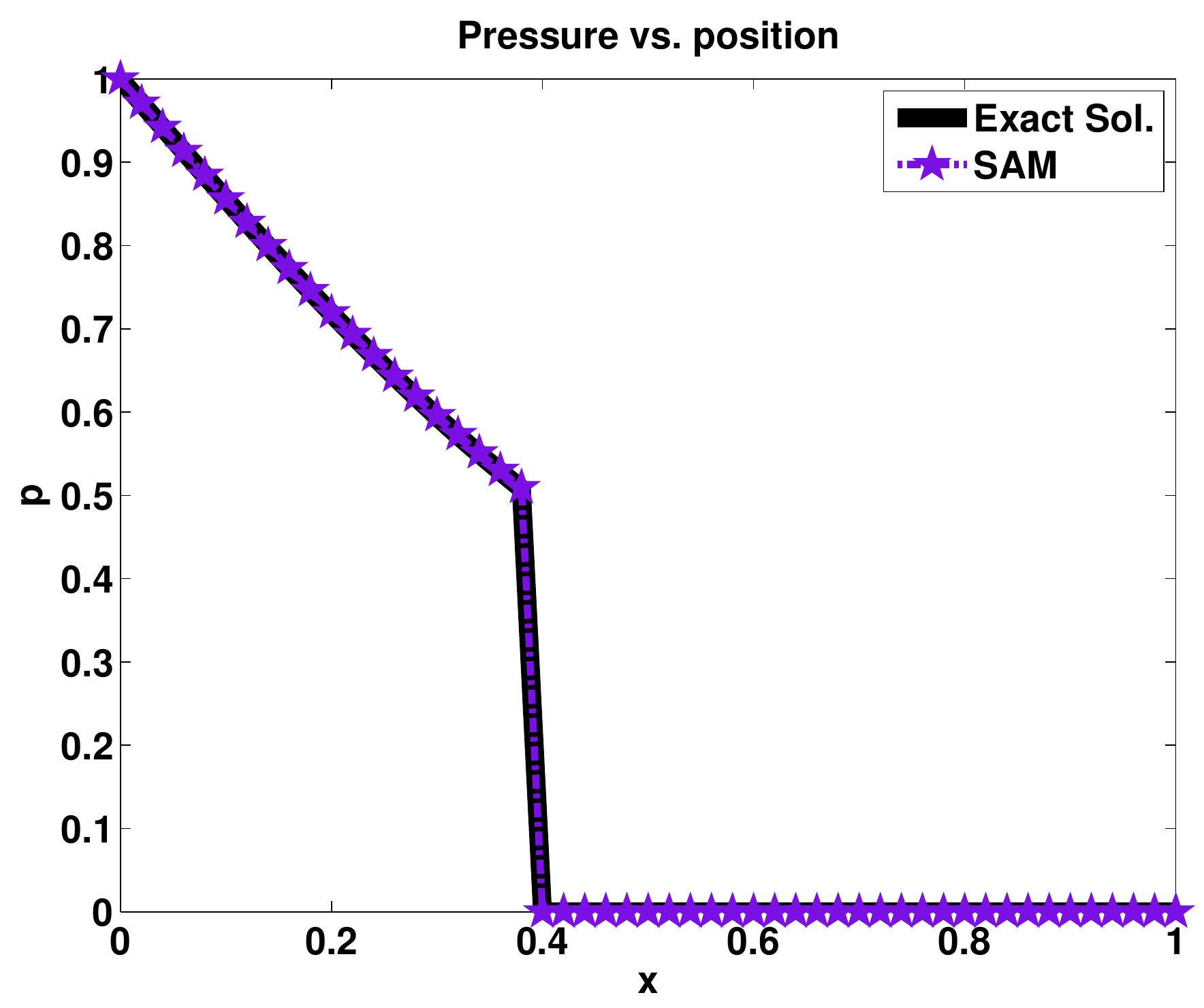}
			\caption{$N = 50$ grid points}
		\end{subfigure}
		\begin{subfigure}[H]{.49\textwidth}  
			\includegraphics[width =\textwidth]{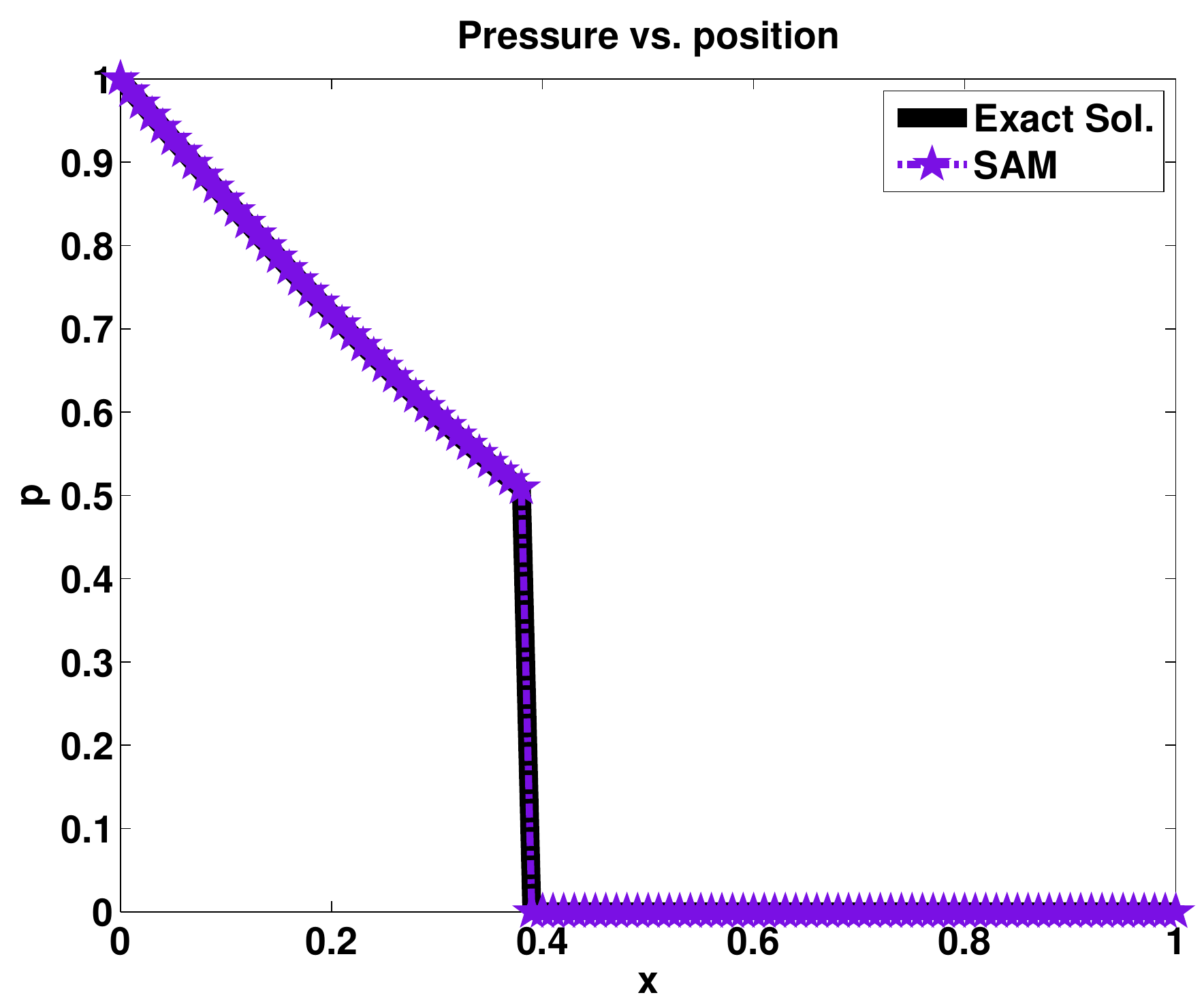}
			\caption{$N = 100$ grid points}
		\end{subfigure}
		\caption{Spatial profiles at time $t = 0.05$ with $\dt  = \dx^2/32$.}
		\label{fig:sam_space}
		\end{figure}
\begin{figure}[H]
		\center
		\begin{subfigure}[H]{0.49\textwidth}  
			\includegraphics[width =\textwidth]{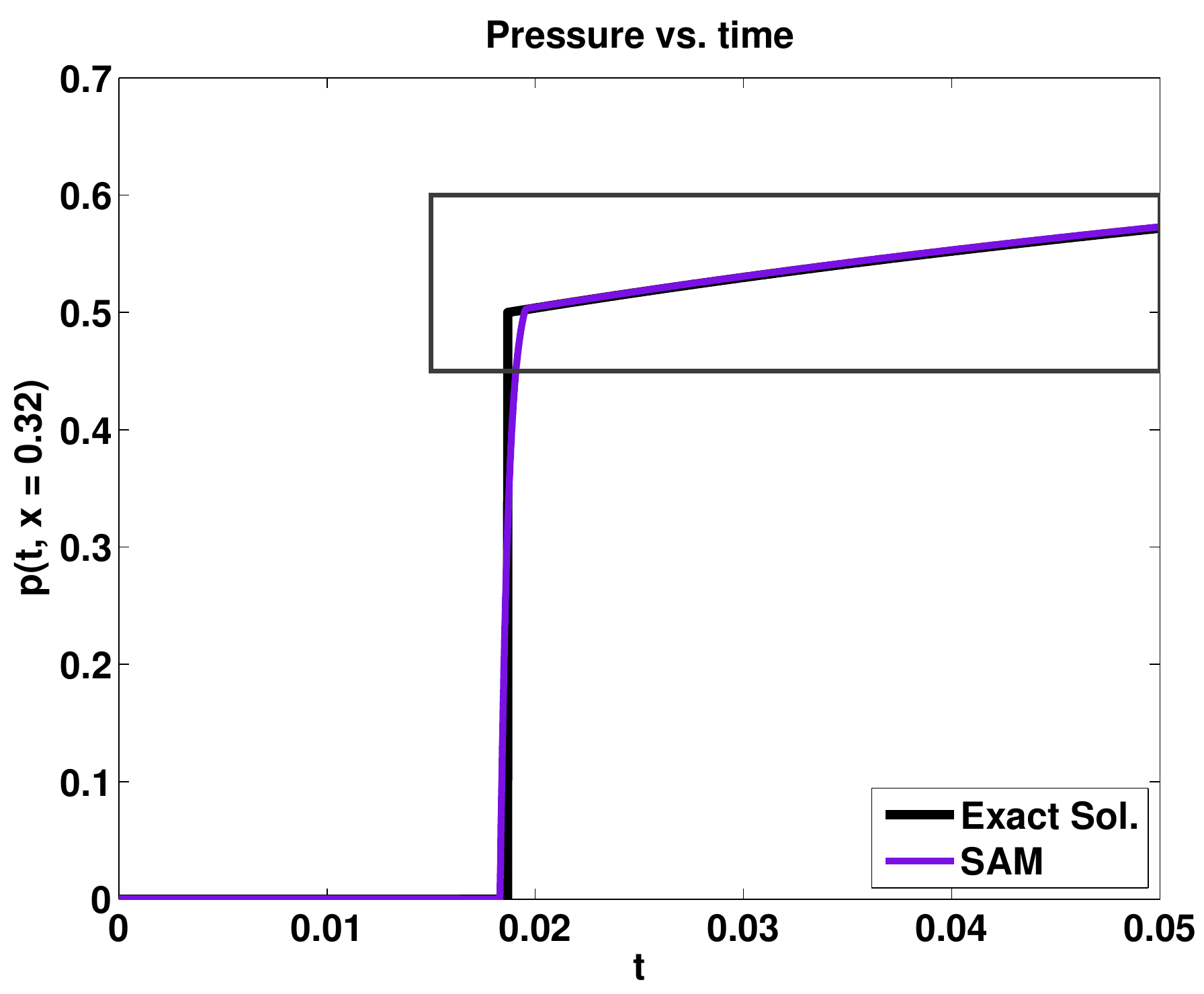}
		\end{subfigure}
		\begin{subfigure}[H]{0.49\textwidth}  
			\includegraphics[width =\textwidth]{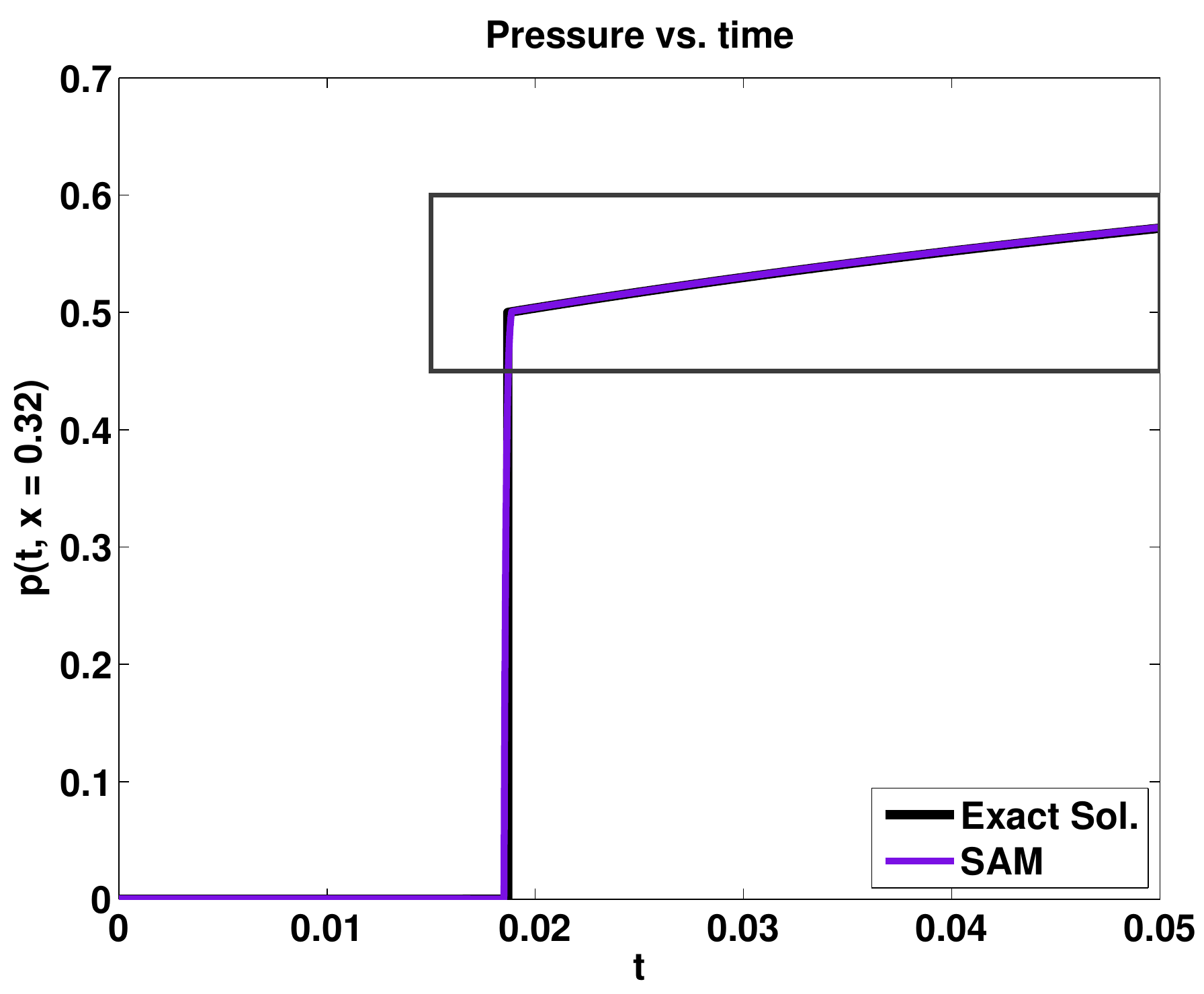}
		\end{subfigure}
		\begin{subfigure}[H]{0.49\textwidth}  
			\includegraphics[width =\textwidth]{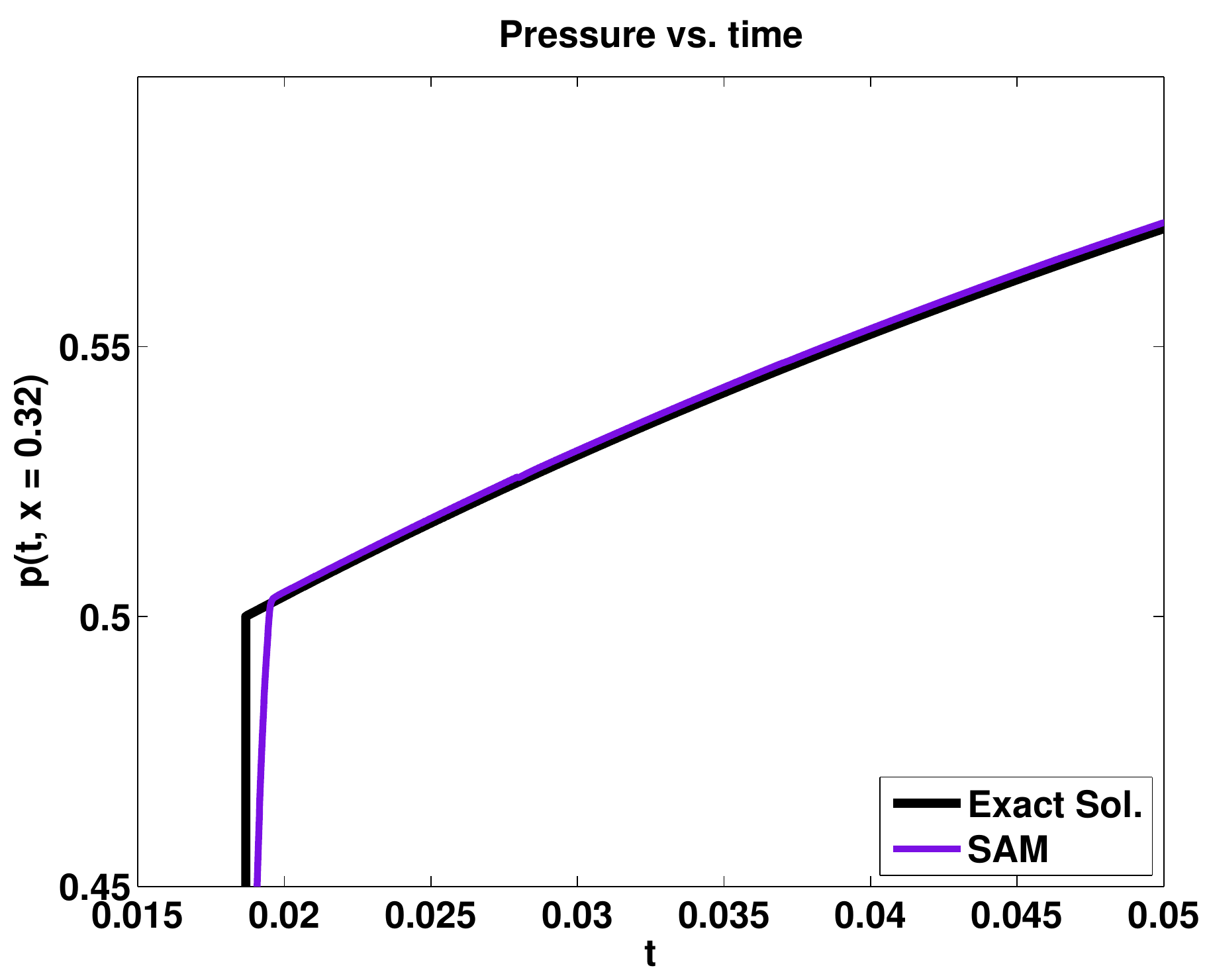}
			\caption{$N = 50$ grid points}
		\end{subfigure}
		\begin{subfigure}[H]{0.49\textwidth}  
			\includegraphics[width =\textwidth]{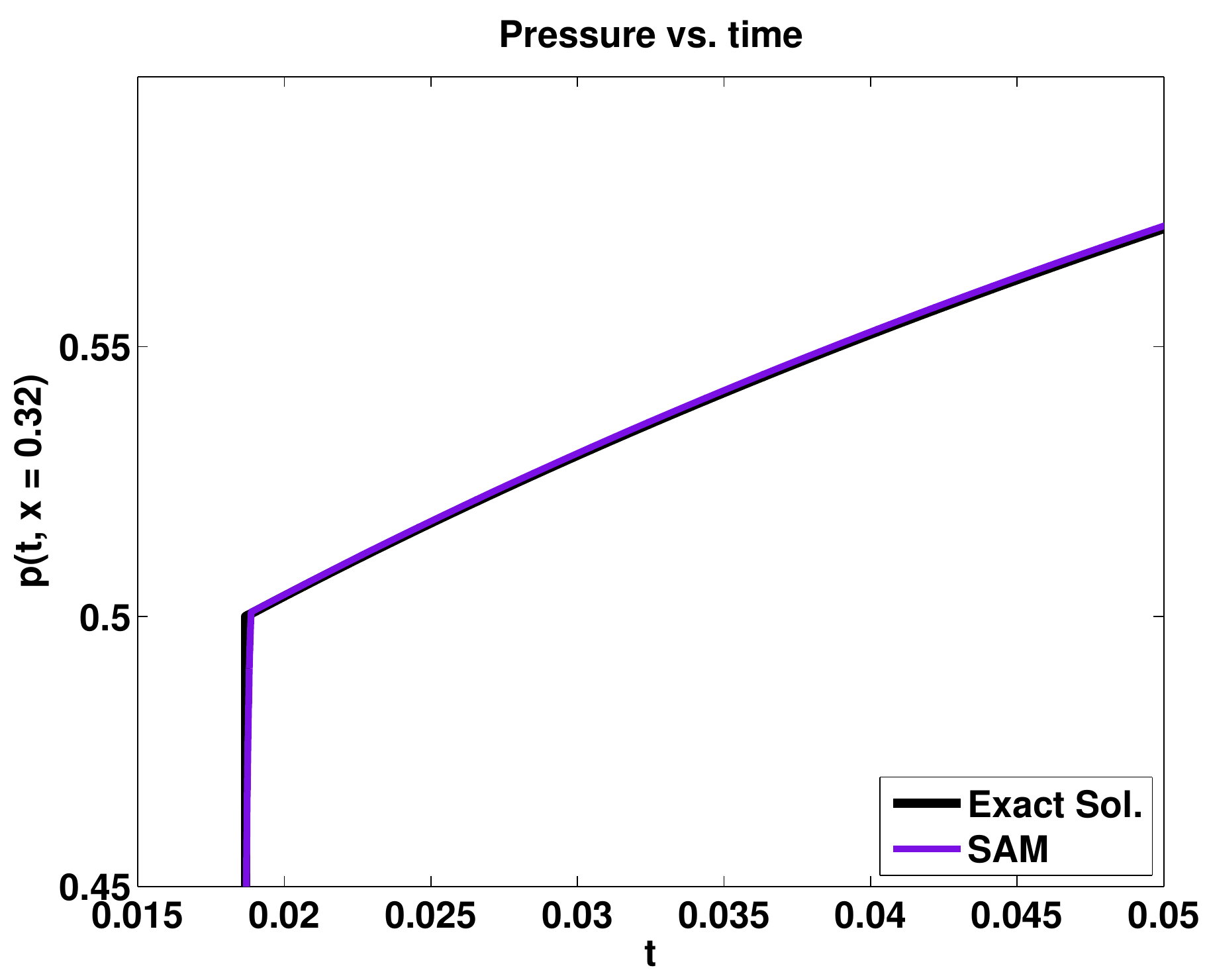}
			\caption{$N = 100$ grid points}
		\end{subfigure}
		\caption {Temporal profiles at position $x = 0.32$.  with $\dt  = \dx^2/32$. }
		\label{fig:sam_time}
\end{figure}

  The results in Figures \ref{fig:sam_space} and \ref{fig:sam_time} on the same test case from the prior section show that there is no significant change in the behavior of the SAM solution with the approximate shock speed than with the exact shock speed in the prior section.  

\begin{figure}[H]
		\center
		\begin{subfigure}[H]{0.33\textwidth}  
			\includegraphics[width =\textwidth]{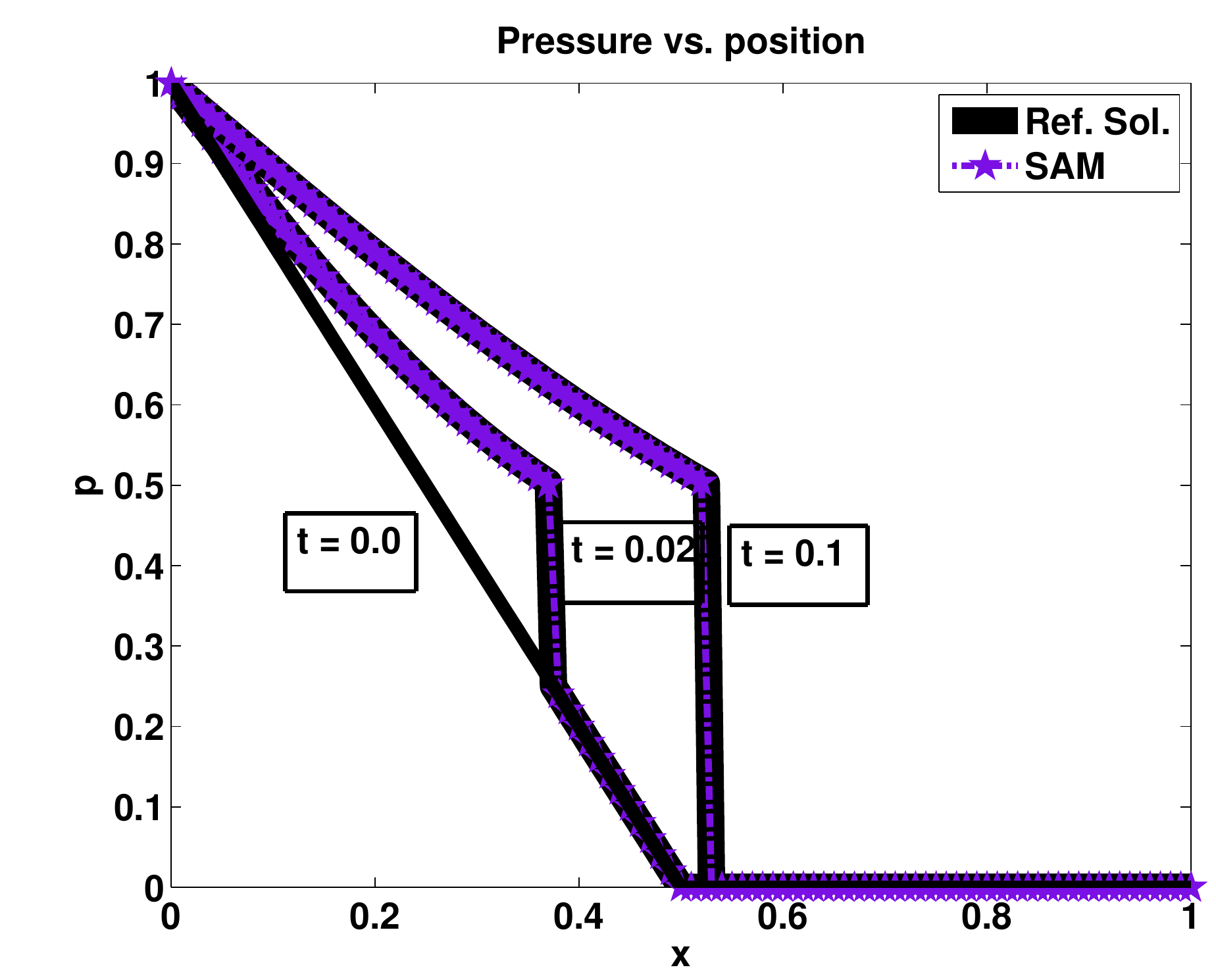}
			\caption{Spatial profiles}
			\label{wait_space}
		\end{subfigure}
		\begin{subfigure}[H]{0.33\textwidth}  
			\includegraphics[width =\textwidth]{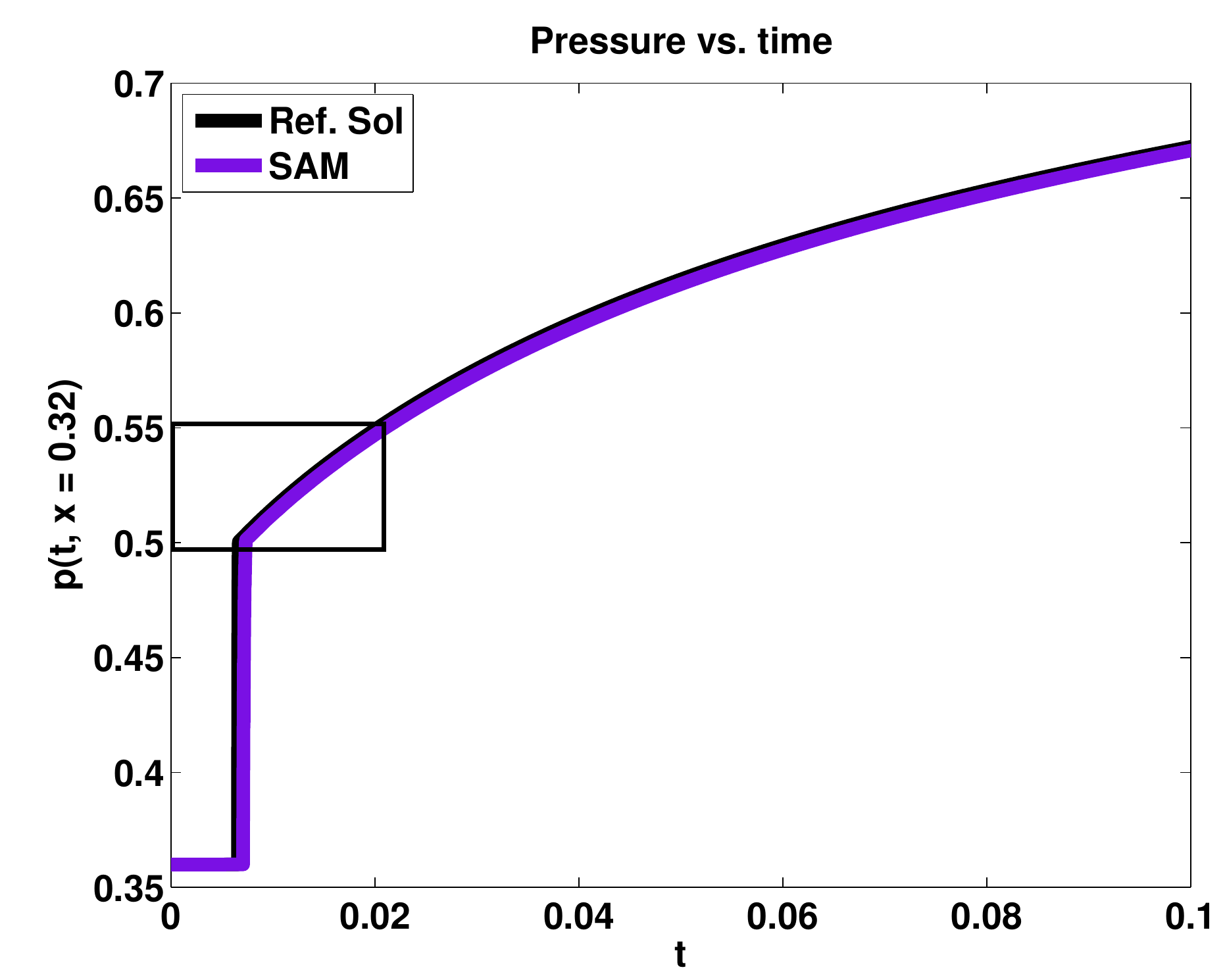}
			\caption{Temporal profiles at position $x = 0.32$.}
			\label{wait_time2}
		\end{subfigure}
		\begin{subfigure}[H]{0.33\textwidth}  
			\includegraphics[width =\textwidth]{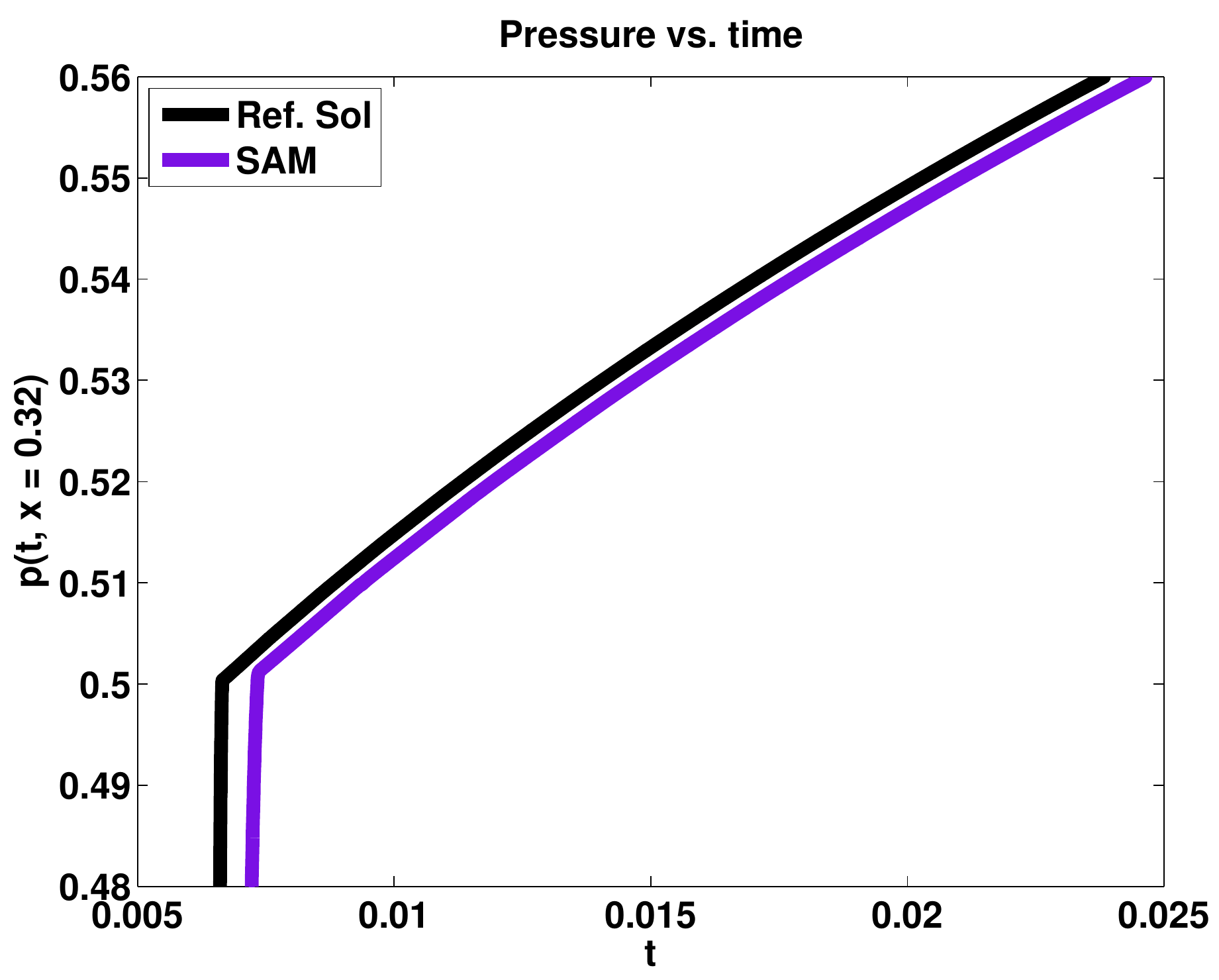}
			\caption{Zoomed in temporal profiles.}
			\label{wait_timezoom}
		\end{subfigure} 
		\caption{SAM numerical solution with $N = 100$ grid points and $\dt = \dx^2/32$ for the waiting time phenomenon example from Figure \ref{wait_time}.}
\end{figure}

Figures \ref{wait_space}-\ref{wait_timezoom} show results for SAM for the waiting phenomenon problem discussed in Section \ref{rh_cond}.  For this initial condition, the exact solution is not known.  The figures show that the SAM solution is accurate and does not possess the numerical artifacts shown in Figure \ref{wait_time}.      
    
\section{Average-based approaches} 
\label{static_avg} 
In the prior sections, we showed that we developed an accurate method, SAM, for the Stefan problem.  
The main goal of the paper is to understand why the artifacts are occurring with finite volume averaged-based methods from the literature.  We will now recast SAM in this framework to shed light on what is happening in these methods that are used frequently in the porous media community.

We implement the finite volume average-based methods using the Forward in Time, Central in Space (FTCS) discretization on a uniform Cartesian grid, as is commonly done in the porous media literature.  The numerical fluxes are given by 
\begin{equation}
	 F_{j}^+ = -k_{j+1/2}\frac{p_{j+1}-p_j}{\dx},
	 \label{eq:out_flux_arith}
\end{equation}
and
\begin{equation}
	F_{j}^- = -k_{j-1/2}\frac{p_{j}-p_{j-1}}{\dx},
	\label{eq:in_flux_arith}
\end{equation}
for an arbitrary average $k_{j+1/2}$ of the neighboring coefficients at the cell face $x_{j+1/2}$ for all $j = 1, \dots, N$.  For verification, we successfully repeated the test case presented for the arithmetic and integral average in \citet{vandermeer2016}.  



\subsection{Arithmetic and Harmonic Averages}
In the FTCS finite volume scheme with arithmetic 
 \begin{equation}
 	k^A_{j+1/2} = \frac{k_j + k_{j+1}}{2}, 
	\label{eq:arith}
\end{equation}
 and harmonic
\begin{equation}
	k^H_{j+1/2} = \frac{2k_jk_{j+1}}{k_j + k_{j+1}},
	\label{eq:harm}
\end{equation}
averaging, 
the fluxes in and out of the control volumes surrounding the shock are discretized using Eqns. \eqref{eq:out_flux_arith}-\eqref{eq:in_flux_arith}. 
The corresponding arithmetic fluxes are given by
\begin{equation}
	{F^+_{i}}^A =  {F^-_{i+1}}^A = -k^A_{i+1/2} \frac{p_{i+1} - p_i}{\dx}.
	\label{eq:arith_flux}
\end{equation}
The analogous expression holds for the harmonic fluxes ${F^+_{i}}^H$ and ${F^-_{i+1}}^H$, where the arithmetic average $k^A_{i+1/2}$ in Eqn. \eqref{eq:arith_flux} is replaced with the harmonic average $k^H_{i+1/2}$.
 
\begin{figure}[H]
	\center
	\begin{subfigure}[H]{0.49\textwidth}  
			\includegraphics[width =\textwidth]{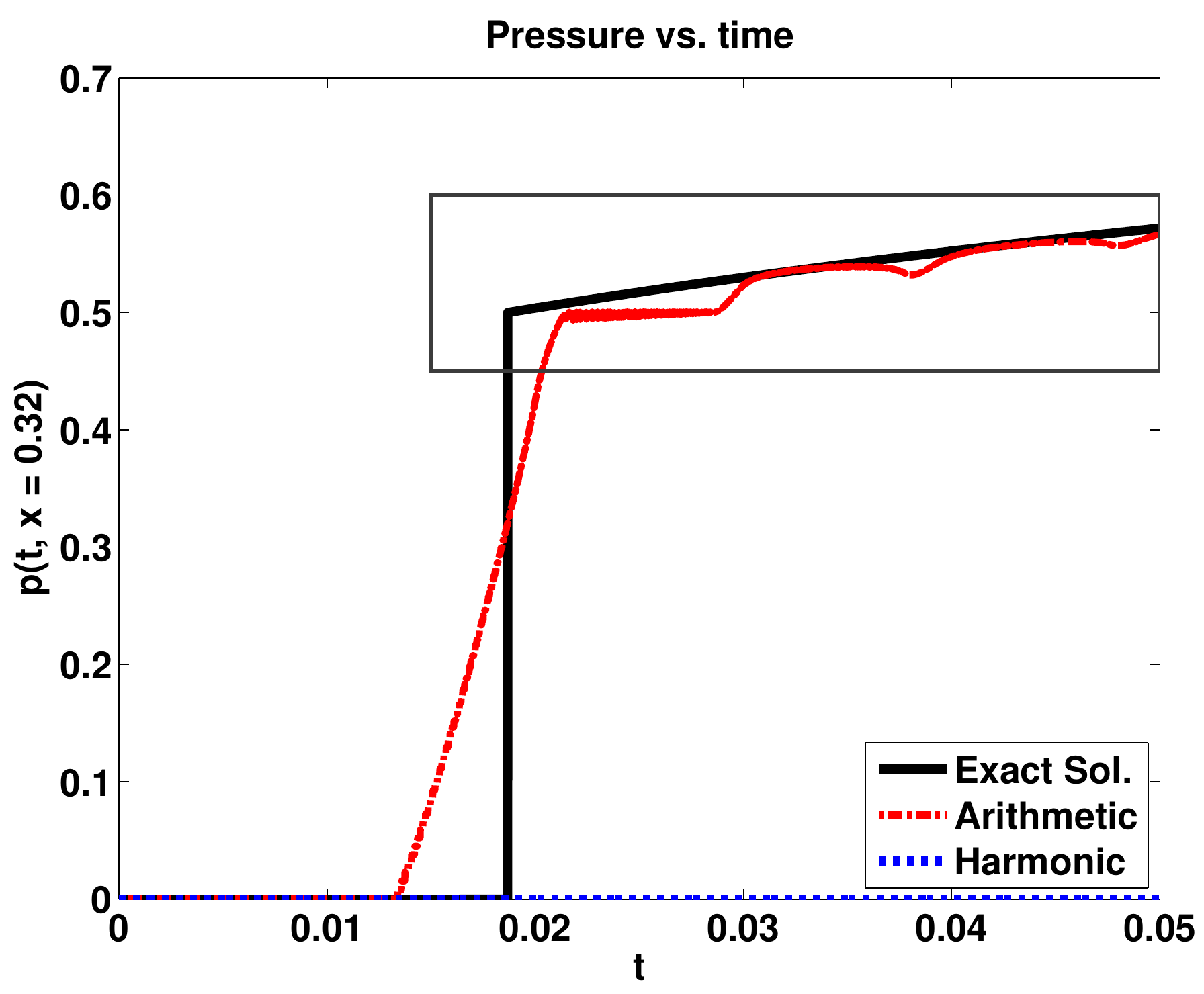}
		\end{subfigure}
		\begin{subfigure}[H]{0.49\textwidth}  
			\includegraphics[width =\textwidth]{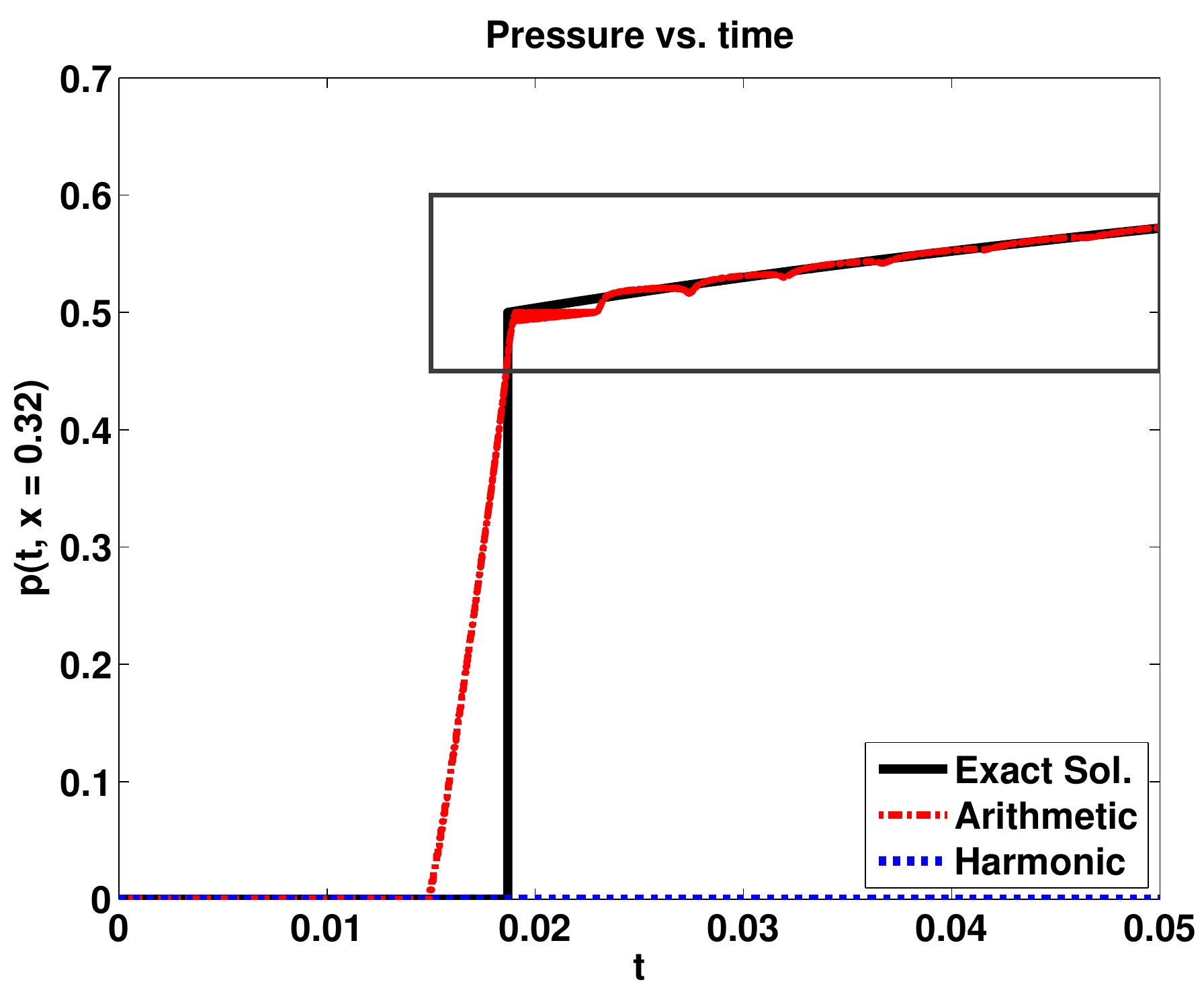}
		\end{subfigure}
		\begin{subfigure}[H]{0.49\textwidth}  
			\includegraphics[width =\textwidth]{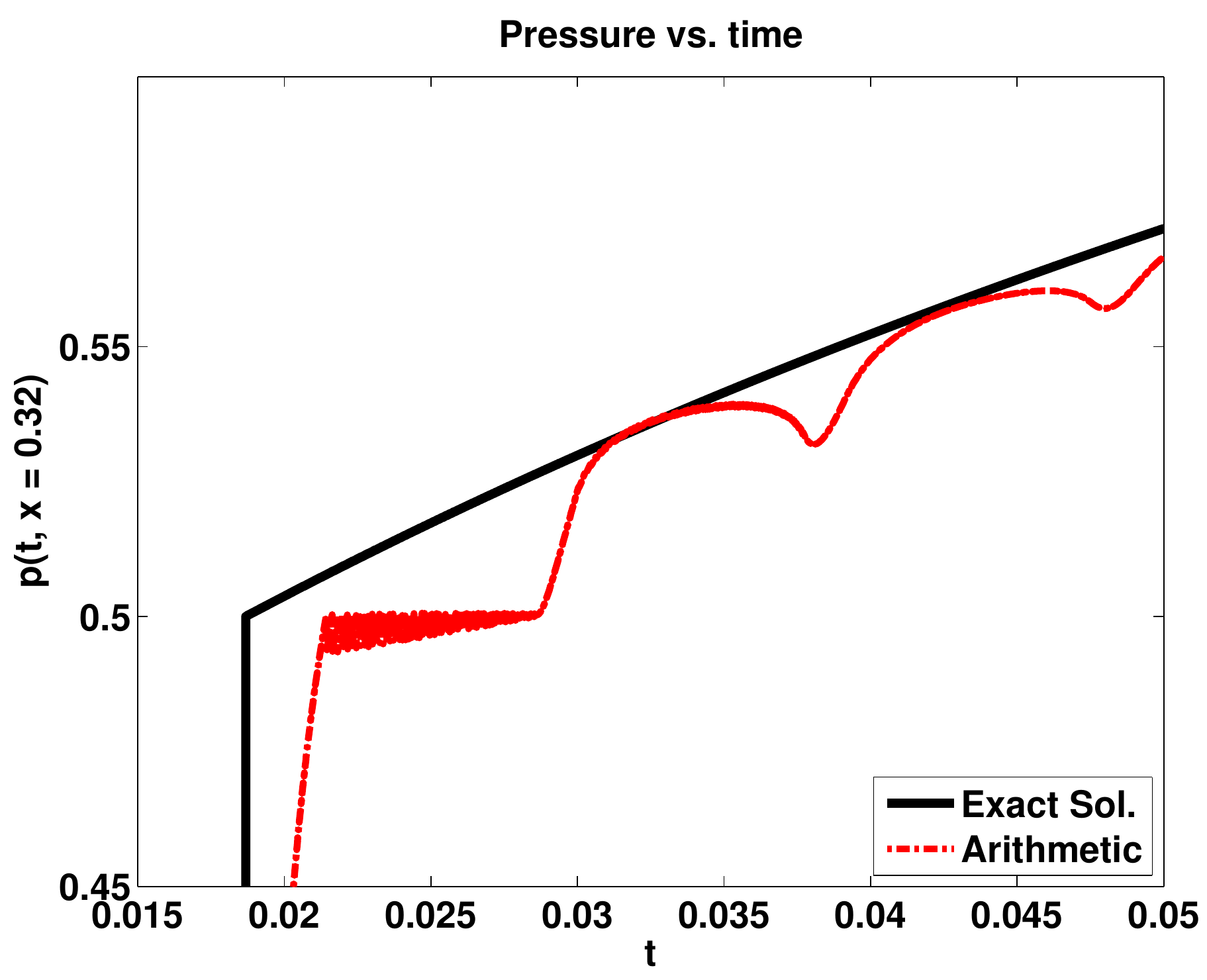}
			\caption{$N = 50$ grid points}
			\label{fig:arith_res50}
		\end{subfigure}
		\begin{subfigure}[H]{0.49\textwidth}  
			\includegraphics[width =\textwidth]{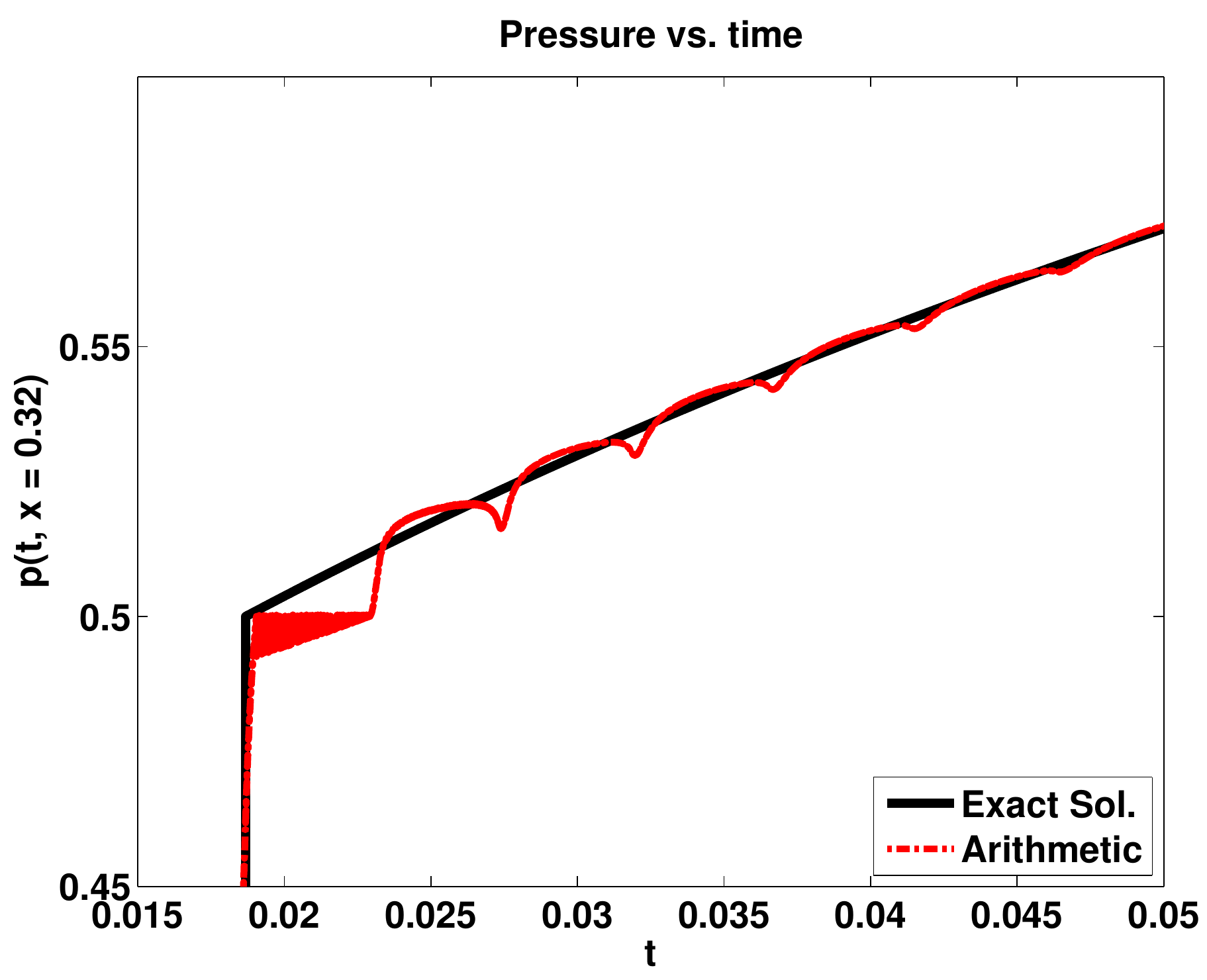}
			\caption{$N = 100$ grid points}
			\label{fig:arith_res100}
		\end{subfigure}
		\caption{Temporal profiles at position $x = 0.32$  with $\dt  = \dx^2/32$.} 
		\label{fig:arith_time}
\end{figure}

\begin{figure}[H]
		\center
		\begin{subfigure}[H]{.49\textwidth}  
			\includegraphics[width =\textwidth]{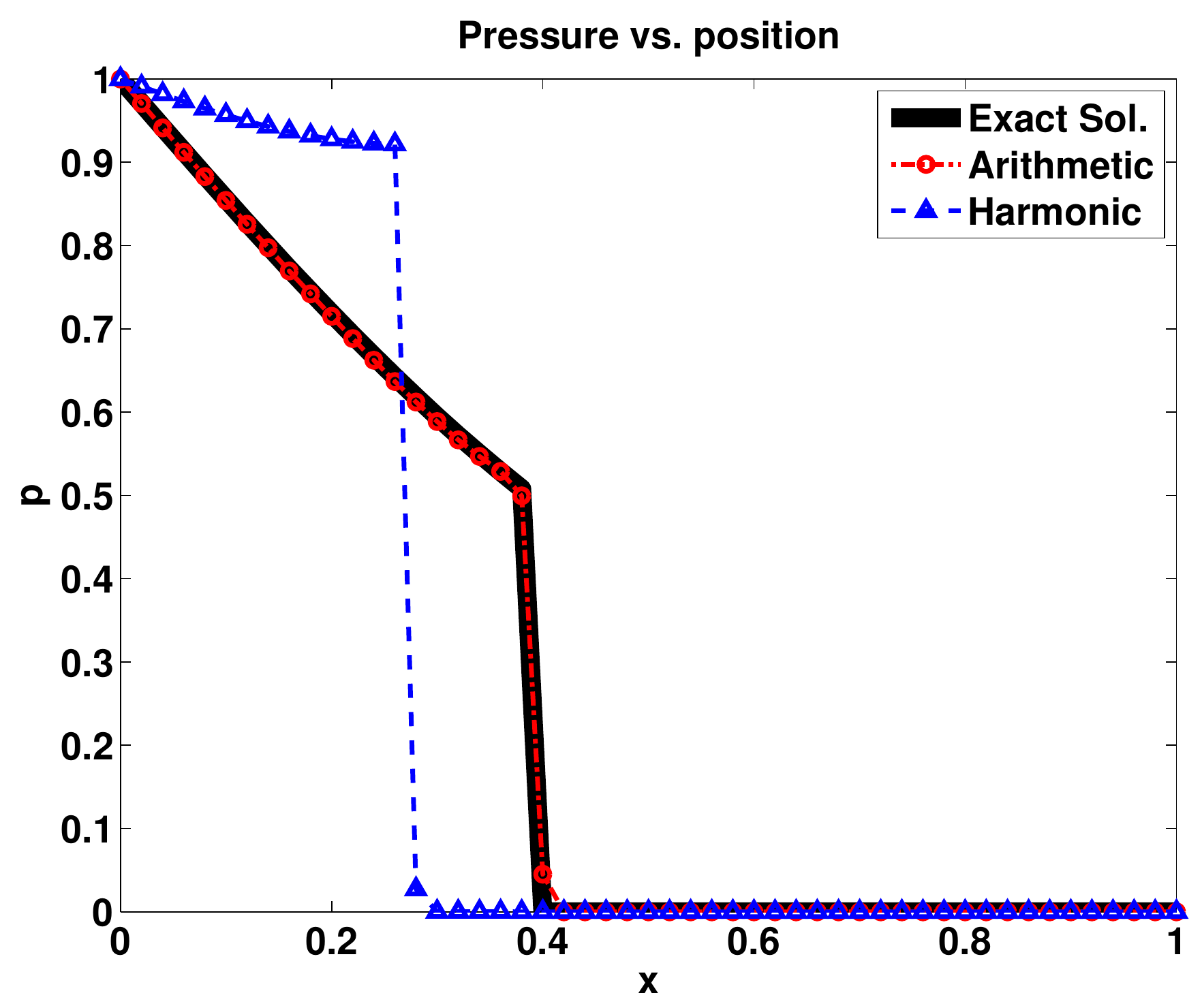}
			\caption{$N = 50$ grid points} 
		\end{subfigure}
		\begin{subfigure}[H]{.49\textwidth}  
			\includegraphics[width =\textwidth]{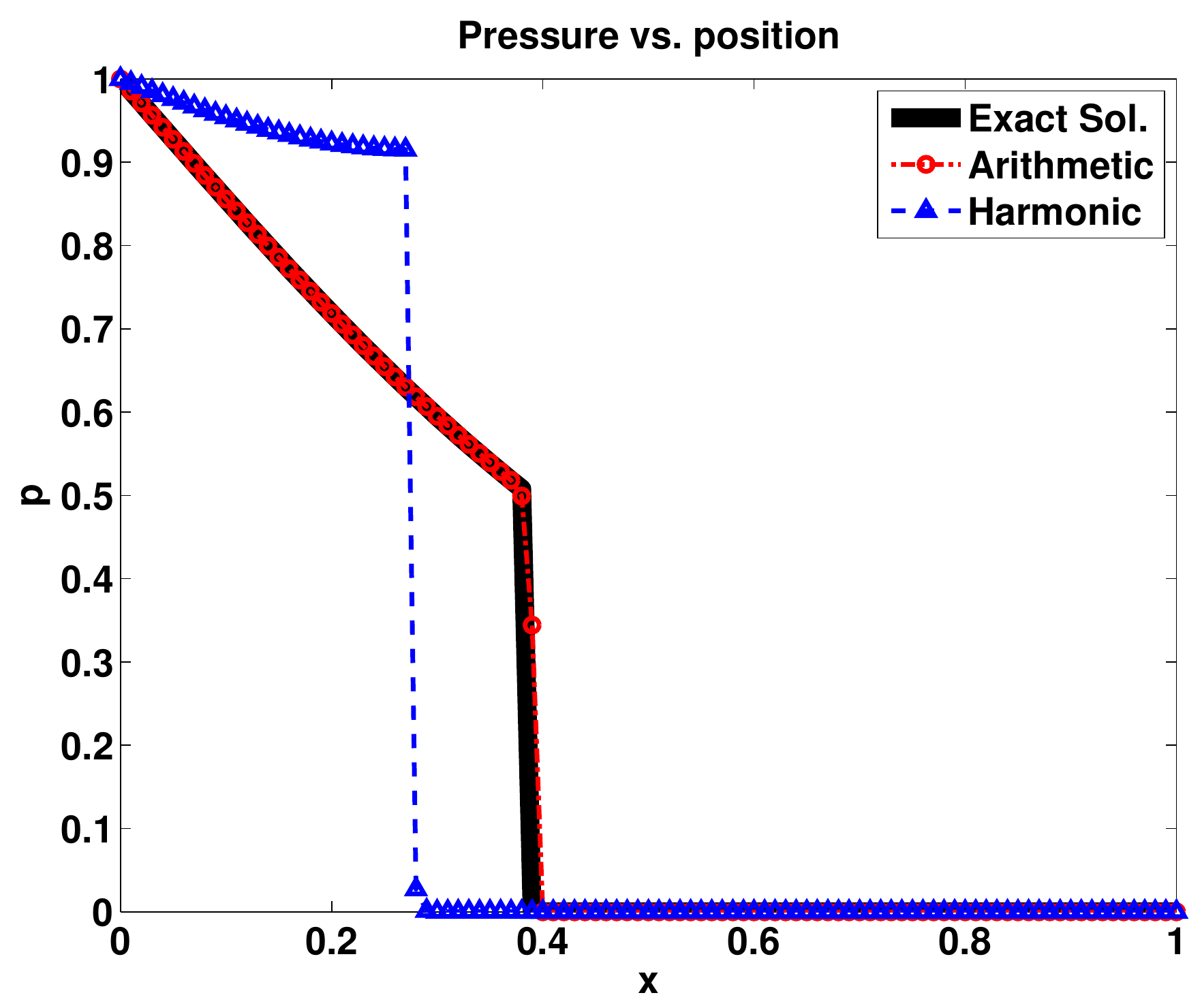}
			\caption{$N = 100$ grid points} 
		\end{subfigure}
				\caption{Spatial profiles at time $t = 0.05$ with $\dt  = \dx^2/32.$} 
		\label{fig:space} 
\end{figure}

\begin{figure}[H]
		\center
		\begin{subfigure}[H]{.32\textwidth}  
			\includegraphics[width =\textwidth]{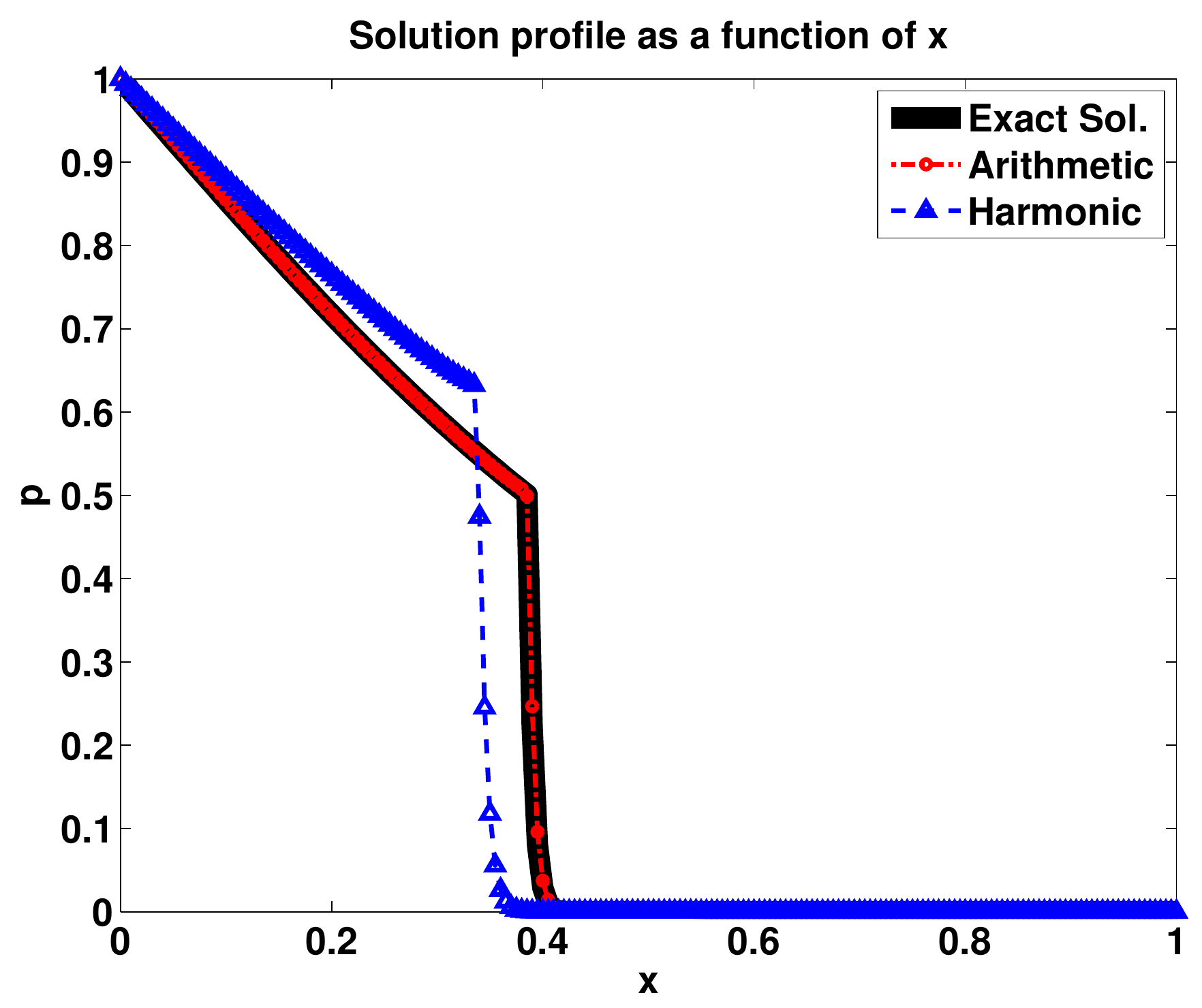}
			\caption{Spatial profiles at $t = 0.05$}
			\label{spatial_01}
		\end{subfigure}
		\begin{subfigure}[H]{.32\textwidth}  
			\includegraphics[width =\textwidth]{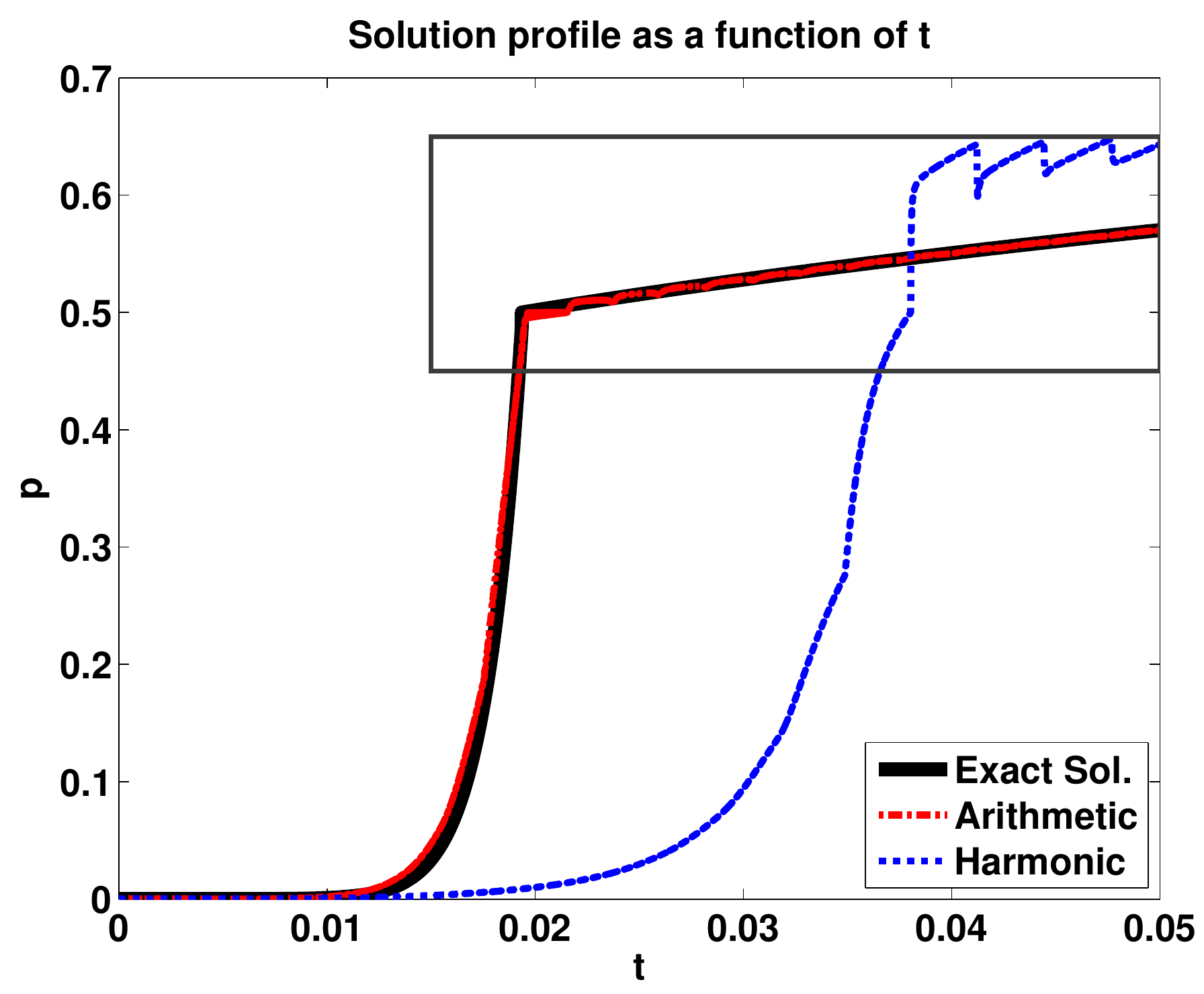}
			\caption{Temporal profile at $x = 0.32$}
			\label{time_01}
		\end{subfigure}
		\begin{subfigure}[H]{0.32\textwidth}  
			\includegraphics[width =\textwidth]{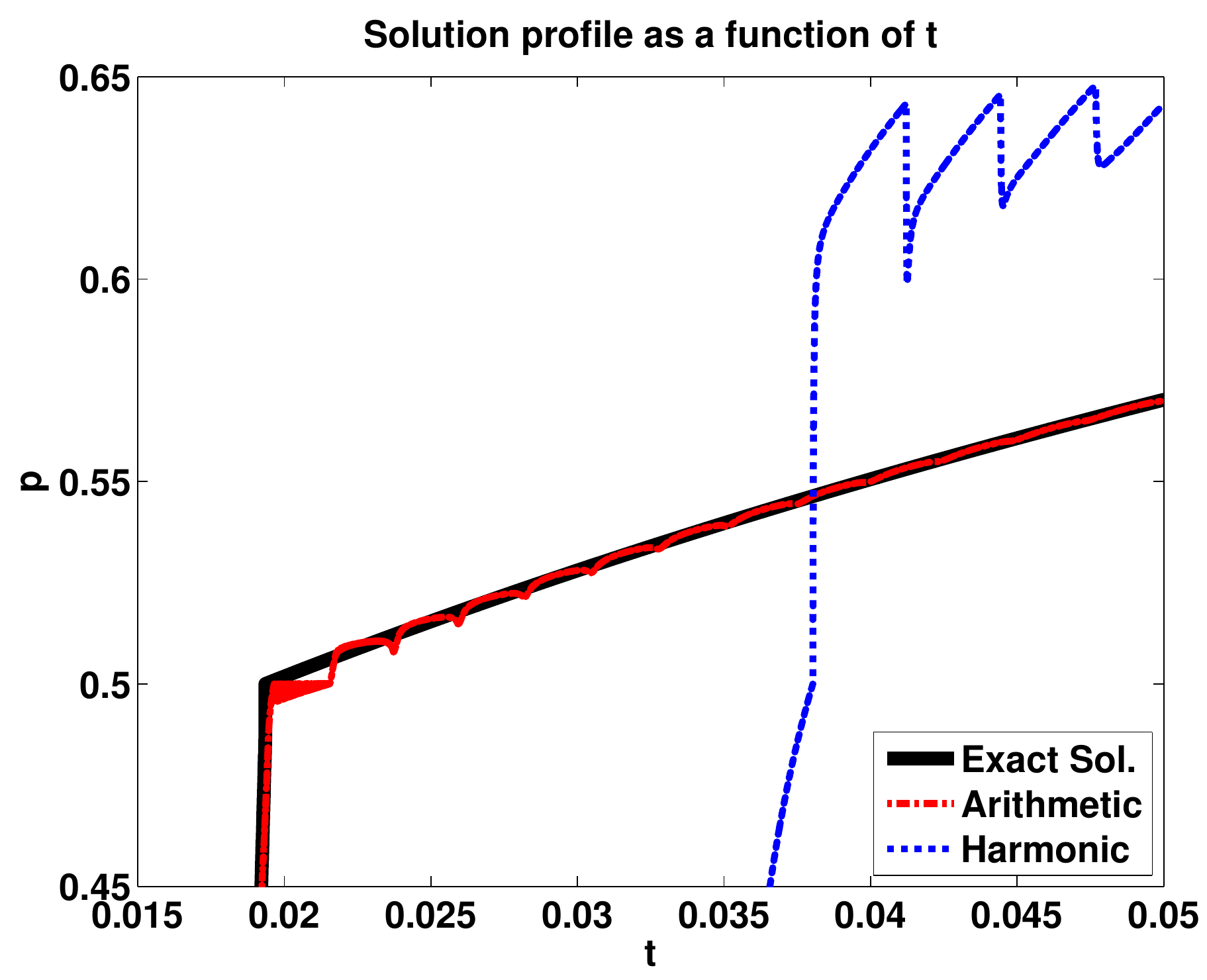}
			\caption{Zoomed in temporal profiles}
			\label{time_01_zoom}
		\end{subfigure}
	\caption{Solution profiles with $N = 200$ grid point and $\dt  = \dx^2/32$.}
	\label{fig:time_01}
\end{figure}

 The temporal plots at $x = 0.32$ in Figure \ref{fig:arith_time} reveal that the numerical interface with harmonic averaging is locked, that is, does not move at all, and so does not advance to this position.  Figure \ref{fig:space} shows the spatial evolution and we again see the locking with harmonic averaging.  This numerical solution in Figure \ref{fig:space} evolves to a step function with values at 1 and 0.  Eqn. \eqref{eq:harm} explains this locking numerical artifact, since $k^H_{j+1/2} = 0$ when either coefficient is zero.  
 Although here locking is seen as a drawback, there are situations, where it can be desirable, such as in variable coefficient problems, $p_t  = \nabla \cdot (k(x) \nabla p)$, when the interface separates a permeable and impermeable material.  The arithmetic average is known to cause leakage across such interfaces.  We could relax $k_{\min} = 0$ to $k_{\min} = \eps$ for some small $\eps > 0$.  Figure \ref{spatial_01} shows that for $\eps = 0.01$, the numerical solution with harmonic averaging no longer completely locks, but is still lagging behind the true front location.  From Eqn. \eqref{eq:harm}, we see that $k^H_{i+1/2}$ favors the smaller coefficient and is approximately equal to $2\eps$, as the shock moves through the interval $[x_i,x_{i+1}]$.  The constant and small averaged $k$ value at the cell face results in a solution whose numerical speed is too slow.  Grid dependent oscillations with much larger amplitudes than those in the solution with arithmetic averaging are present in its temporal profile in Figures \ref{time_01}-\ref{time_01_zoom}.  In \cite{vandermeer2016}, the harmonic average is not considered, since it is known in the literature to behave poorly for the GPME with near-zero coefficients.  In the rest of this subsection, we focus on the behavior of the solution with arithmetic averaging.
The numerical solutions depicted in Figures \ref{fig:arith_time} and \ref{fig:space} show results for the arithmetic average.  We see that the numerical solution with arithmetic averaging does not sharply capture the shock.  Figure \ref{fig:arith_time} shows spurious oscillations of low and high frequency.  The low frequency oscillations are grid dependent.  The high frequency oscillations are apparent in the zoomed in region around the solution shock value of $p^* = 0.5$.  Figure \ref{fig:space} shows that the solution with arithmetic averaging does not have spatial oscillations.  In the temporal view, oscillations are present.  

There are two main differences in the artifacts for the continuous and discontinuous GPME.  Here, for the discontinuous case, there are temporal oscillations of low and high frequencies for the arithmetic average.  For the continuous GPME, temporal oscillations were only observed in solutions with harmonic averaging  \cite{maddix_pme}. 
The lower frequency grid dependent oscillations were analyzed in \cite{maddix_pme} for discretizations with harmonic averaging.  The high frequency oscillations are an additional numerical artifact with arithmetic averaging that were not observed for the harmonic average in \cite{maddix_pme}.  

The low frequency temporal oscillations have been studied in the literature \cite{maddix_pme, vandermeer2016, zanganeh14, rossen99}.  It has been shown in these works that the temporal oscillations occur as the front crosses a grid cell.  The spurious oscillations then do not vanish with grid refinement.  The frequency increases and the amplitude decreases, as the number of grid points increases.  The differences in frequency and amplitude are illustrated by the solutions for $N = 50$ in Figure \ref{fig:arith_res50} and $N = 100$ in Figure \ref{fig:arith_res100}.  Another observed characteristic of these oscillations in \cite{vandermeer2016, zanganeh14, rossen99} is that the amplitude decreases as shock moves further away from the solution probe point.  This is because the possible error in the shock position relative to the distance from the probe point to the shock decreases.  

\begin{figure}[H]
		\center
		\begin{subfigure}[H]{0.49\textwidth}  
			\includegraphics[width =\textwidth]{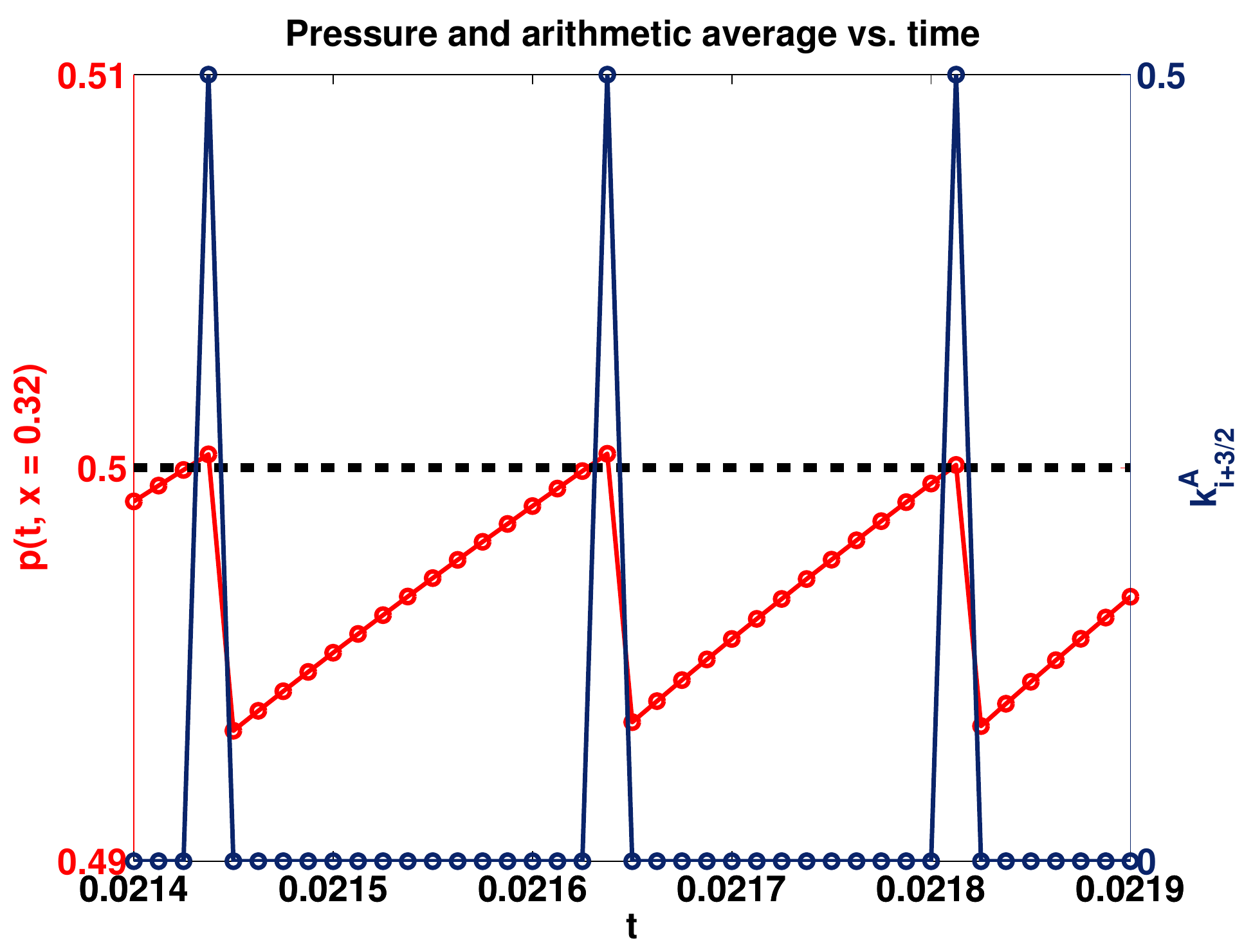}
			\caption{$\dt  = \dx^2/32$}
			\label{fig:arith_oscill_50}
		\end{subfigure}
		\begin{subfigure}[H]{0.49\textwidth}  
			\includegraphics[width =\textwidth]{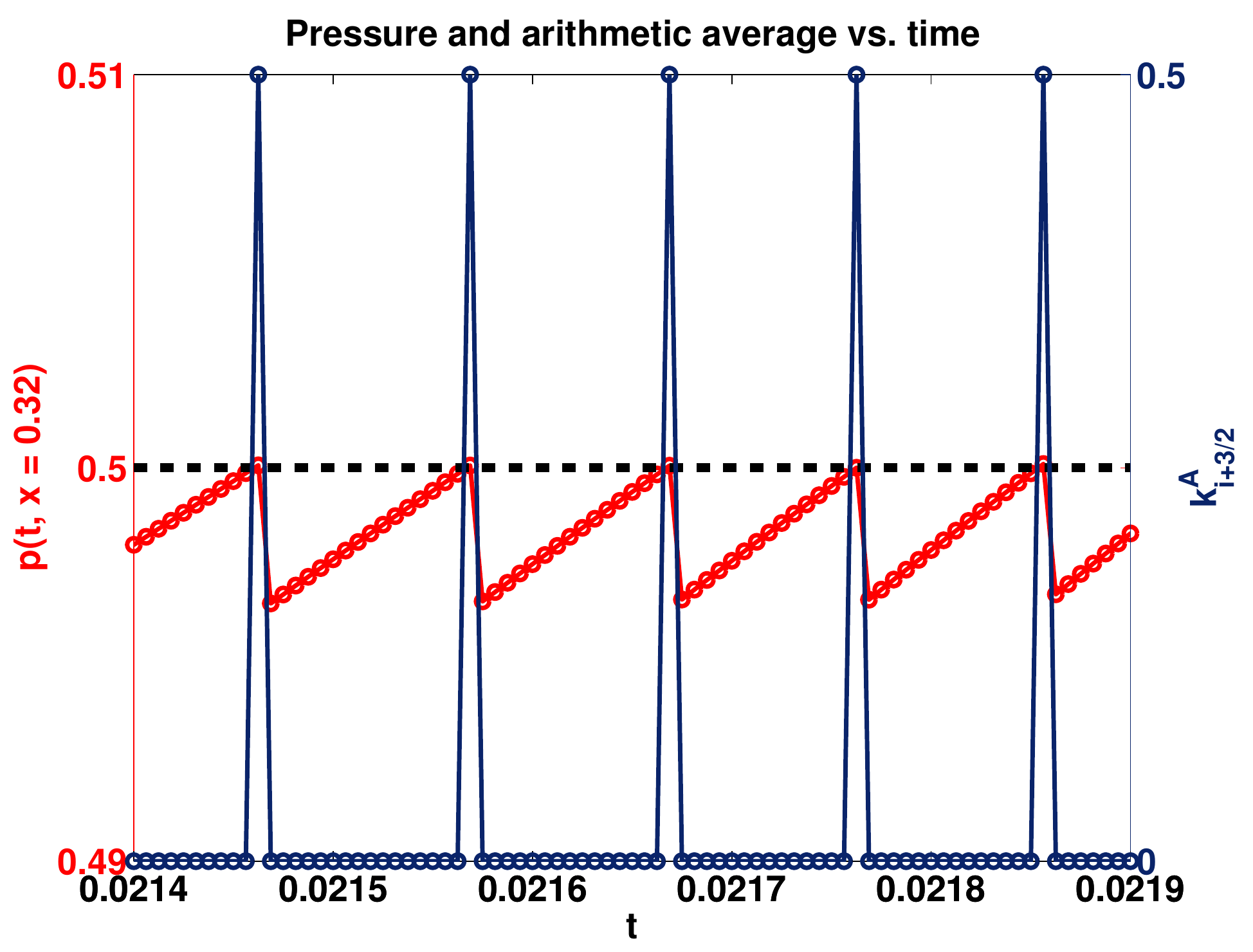}
			\caption{$\dt  = \dx^2/64$}
			\label{fig:arith_oscill_100}
		\end{subfigure}
						\caption{Zoomed in region of the high frequency temporal oscillations near $p^* = 0.5$ in the left cell $p_{i+1}$ for the arithmetic average $k^A_{i+3/2}$ with $N = 50$ grid points at position $x = 0.32$.} 
\end{figure}

The high frequency oscillations were also observed in \cite{vandermeer2016}.  Figures \ref{fig:arith_oscill_50}-\ref{fig:arith_oscill_100} illustrate the temporal profile on a shorter time interval, where the high frequency oscillations occur for $p$ values near $p^* = 0.5$.  We see that $p_{i+1}$ is slowly increasing, until it crosses the threshold at $p^* = 0.5$.  The corresponding $k_{i+1}$ then jumps from $k_{\min} = 0$ to $k_{\max} = 1$, according to the model for $k(p)$ in Eqn. \eqref{eq:discont_k}.  The arithmetic average $k^A_{i+3/2}$ in Eqn. \eqref{eq:arith} jumps from 0 to 0.5, as illustrated in Figures \ref{fig:arith_oscill_50}-\ref{fig:arith_oscill_100}.  The increase in $k^A_{i+3/2}$ causes $p$ to drop below $p^* = 0.5$ and $k^A_{i+3/2}$ to jump back down to 0 at the next time step.  The cycle then repeats itself.

\begin{figure}[H]
		\center
		\begin{subfigure}[H]{0.49\textwidth}  
			\includegraphics[width =\textwidth]{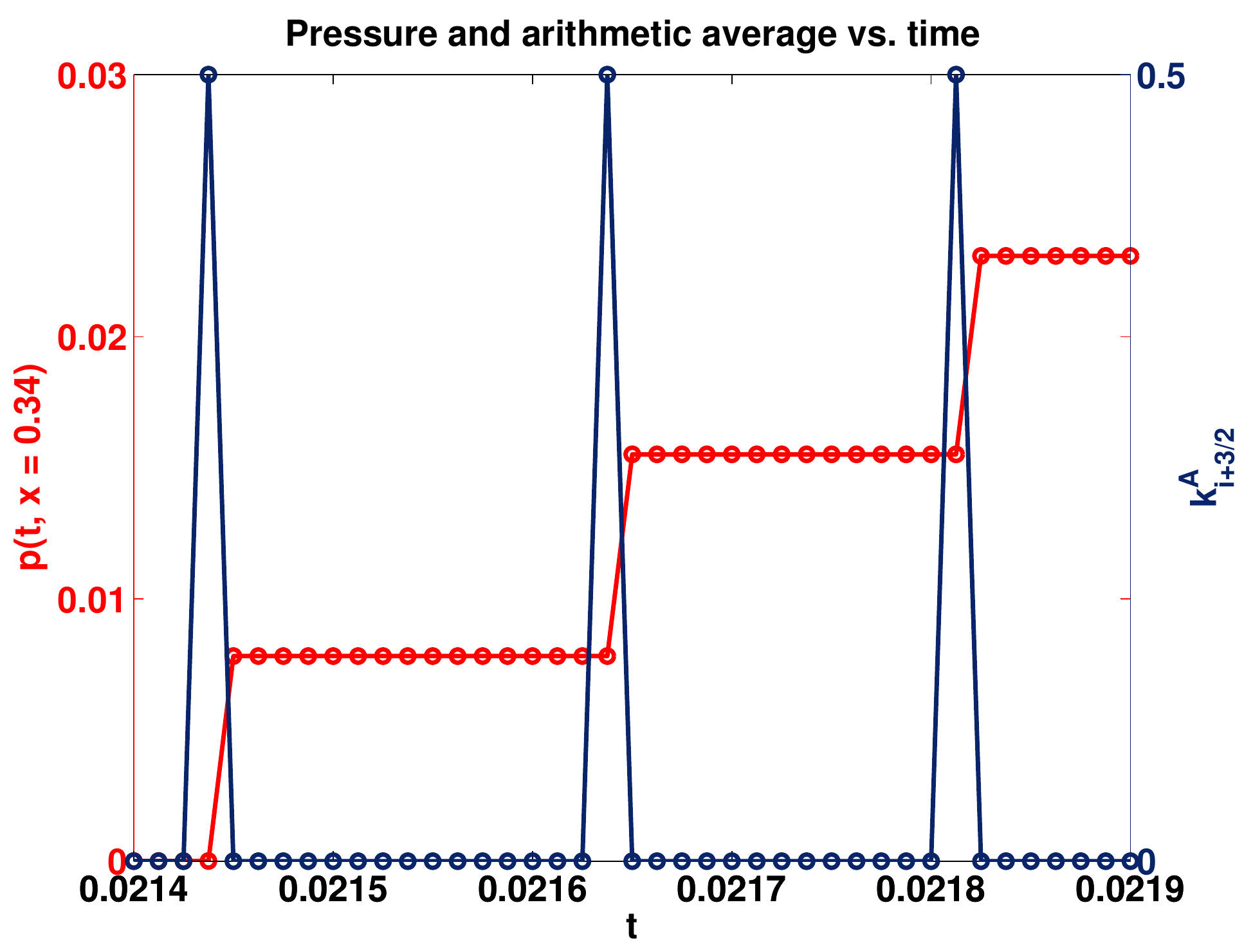}
			\caption{$\dt  = \dx^2/32$}
			\label{fig:arith_oscill_50left}
		\end{subfigure}
		\begin{subfigure}[H]{0.49\textwidth}  
			\includegraphics[width =\textwidth]{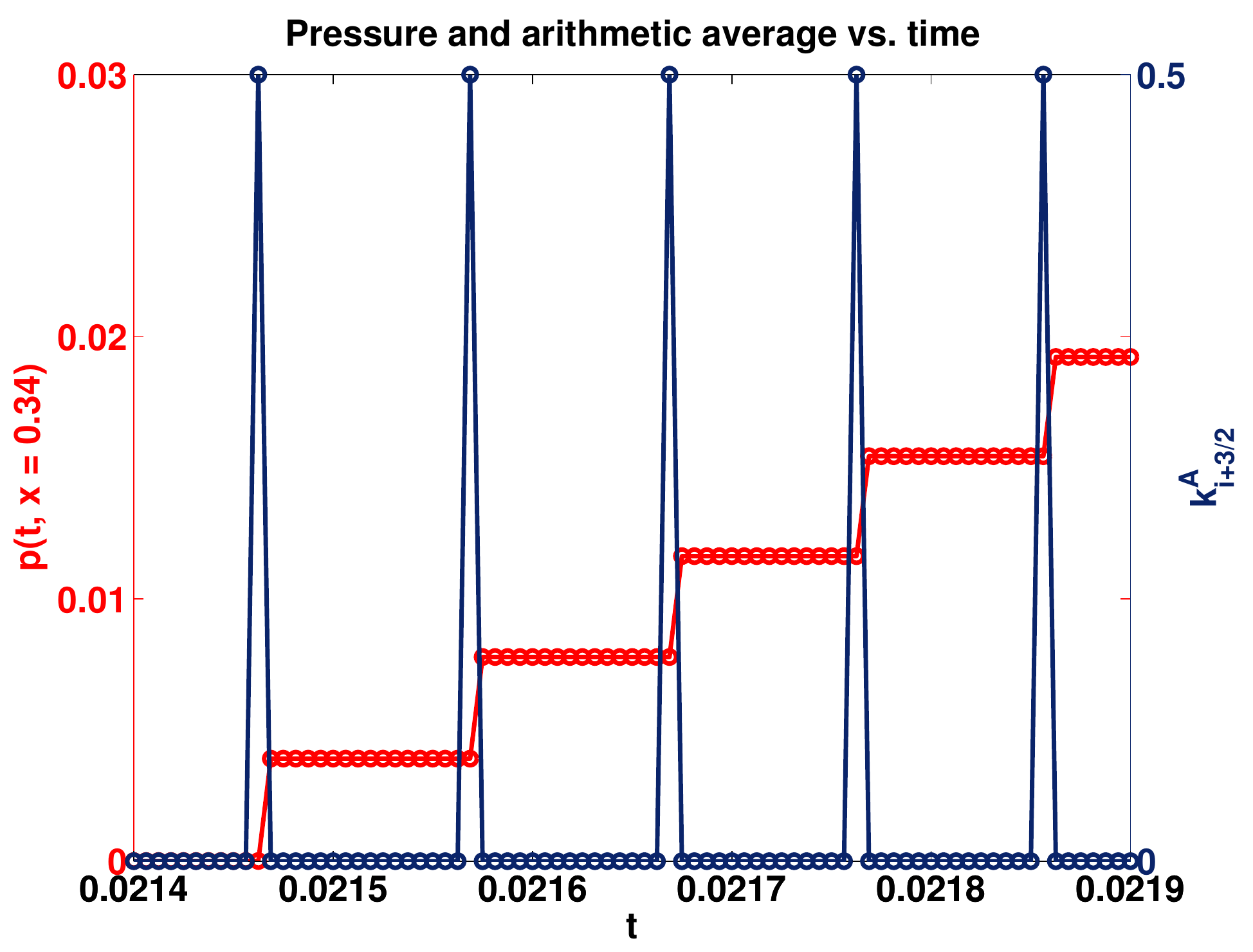}
			\caption{$\dt  = \dx^2/64$}
			\label{fig:arith_oscill_100left}
		\end{subfigure}
						\caption{Zoomed in region of the solution in the right cell $p_{i+2}$ for the arithmetic average $k^A_{i+3/2}$ with $N = 50$ grid points at position $x = 0.34$.} 
\end{figure}

We use the FTCS discretization to explain the effect of the jump in $k^A_{i+3/2}$ on the $p$ value at the next time step.  We discretize the semi-discrete equation \eqref{eq:fv_discret} with Forward Euler in time and substitute in the arithmetic fluxes to obtain
 \begin{equation}
 	\begin{aligned}
	p^{n+1}_{i+1} = p^n_{i+1} + \frac{\dt}{\dx}\bigg(-k^A_{i+1/2}\Big[\frac{p^n_{i+1}-p^n_{i}}{\dx}\Big] + k^A_{i+3/2}\Big[\frac{p^n_{i+2}-p^n_{i+1}}{\dx}\Big]\bigg),
		\label{eq:scheme}
	\end{aligned}
\end{equation}
where $n$ is a time step before the solution drop and $p^n_{i+1} \ge 0.5$.
The first term in Eqn. \eqref{eq:scheme} is small in magnitude, since $p^n_i \approx 0.5$.  The second term in Eqn. \eqref{eq:scheme} has a larger gradient of $p$, since $p^n_{i+2} \approx 0$.  It is also negative, and so the increase in $k^A_{i+3/2}$ results in the drop in $p^{n+1}_{i+1}$.  At the following time steps, $k^A_{i+3/2} = 0$, until the solution crosses the threshold again.  

As $p_{i+1}$ is increasing and is below the threshold, $k^A_{i+1/2}$ is fixed at $0.5$, regardless of the shock position in the cell.  
From Eqn. \eqref{eq:scheme}, we see that the positive quantity $\dt/(2\dx^2)[p_i^n - p_{i+1}^n]$ is added to the current $p$ value at each time step.  The constant arithmetic average at the interface $x_{i+1/2}$ allows the solution to artificially increase above $p^* = 0.5$, resulting in the high frequency oscillations. 

Figures \ref{fig:arith_oscill_50}-\ref{fig:arith_oscill_100} also show that the high frequency oscillations are dependent on the time step size.  As $\dt$ is decreased by half, the number of oscillations in Figure \ref{fig:arith_oscill_100} doubles from those in Figure \ref{fig:arith_oscill_50}.  The amplitude of the high frequency oscillations also decreases, as $\dt$ decreases.  The amplitude decrease is expected from Eqn. \eqref{eq:scheme}, since the additional term is proportional to $\dt$.  In \cite{vandermeer2016}, it is observed that as the time step is decreased, the arithmetic average solution in the high frequency region reaches a constant state at $p^* = 0.5$, rather then converging to the true solution. 

Figures \ref{fig:arith_oscill_50}-\ref{fig:arith_oscill_100} display the profile of $p_{i+1}$ over a time window when the shock is in the interval $[x_i,x_{i+1}]$.  In Figures  \ref{fig:arith_oscill_50left}-\ref{fig:arith_oscill_100left}, we now look at the profile of $p_{i+2}$ over the same time window to see the effect of the jump in $k^A_{i+3/2}$ on the neighboring cell.  Figures \ref{fig:arith_oscill_50left}-\ref{fig:arith_oscill_100left} show that the profile of $p_{i+2}$ is an increasing piecewise constant.  Since $p_{i+2}$ and $p_{i+3}$ are both less than $p^*$, $k^A_{i+5/2} = 0$ and the solution update for $p^{n+1}_{i+2}$ reduces to
 \begin{equation}
 	\begin{aligned}
p_{i+2}^{n+1} &=  
	p^n_{i+2} + k^A_{i+3/2}\bigg(p^n_{i+1}-p^n_{i+2}\bigg)\frac{\dt}{\dx^2}.
		\label{eq:scheme2}
	\end{aligned}
\end{equation}
Eqn. \eqref{eq:scheme2} shows that at the times when $p^n_{i+1} < p^* = 0.5$, $k^A_{i+3/2} = 0$, and so $p_{i+2}^{n+1} = p_{i+2}^{n}$.  Otherwise, at the times when $k^A_{i+3/2}$ jumps to 0.5, the solution increases proportional to $\dt$.  As expected, the increase in out-flux causes the solution to decrease in the left cell (Figures \ref{fig:arith_oscill_50}-\ref{fig:arith_oscill_100}) and the same increase in in-flux causes the solution in the right cell to increase (Figures \ref{fig:arith_oscill_50left}-\ref{fig:arith_oscill_100left}). 

\subsection{Integral Average}
\citet{vandermeer2016} developed the integral average 
\begin{equation}
	k^{I}_{j+1/2} =  \frac{ \int_{p_j}^{p_{j+1}} k(\tilde{p})d\tilde{p}}{p_{j+1}-p_j},
	\label{eq:int_avg}
\end{equation}
which is effective in reducing the numerical artifacts.    The integral average is derived by expressing the coefficient $k(p) = \Phi_p$ in Eqn. \eqref{eq:GPME} to obtain 
\[
	p_t =  \nabla \cdot (\Phi_p \nabla p).
	\label{eq:int_form}
\] 
Discretizing $\Phi_p$ directly with central differences at the cell face $x_{j+1/2}$ gives 
\begin{equation}
	k^{I}_{j+1/2} = \frac{\Phi(p_{j+1}) - \Phi(p_j)}{p_{j+1} - p_{j}}.
	\label{eq:int_central}
\end{equation} 
Using the definition of $\Phi(p) = \int_0^p k(\tilde{p})d\tilde{p}$, Eqn. \eqref{eq:int_central} simplifies to Eqn. \eqref{eq:int_avg}.

\begin{figure}[H]
\center
	\begin{subfigure}[H]{0.49\textwidth}  
			\includegraphics[width =\textwidth]{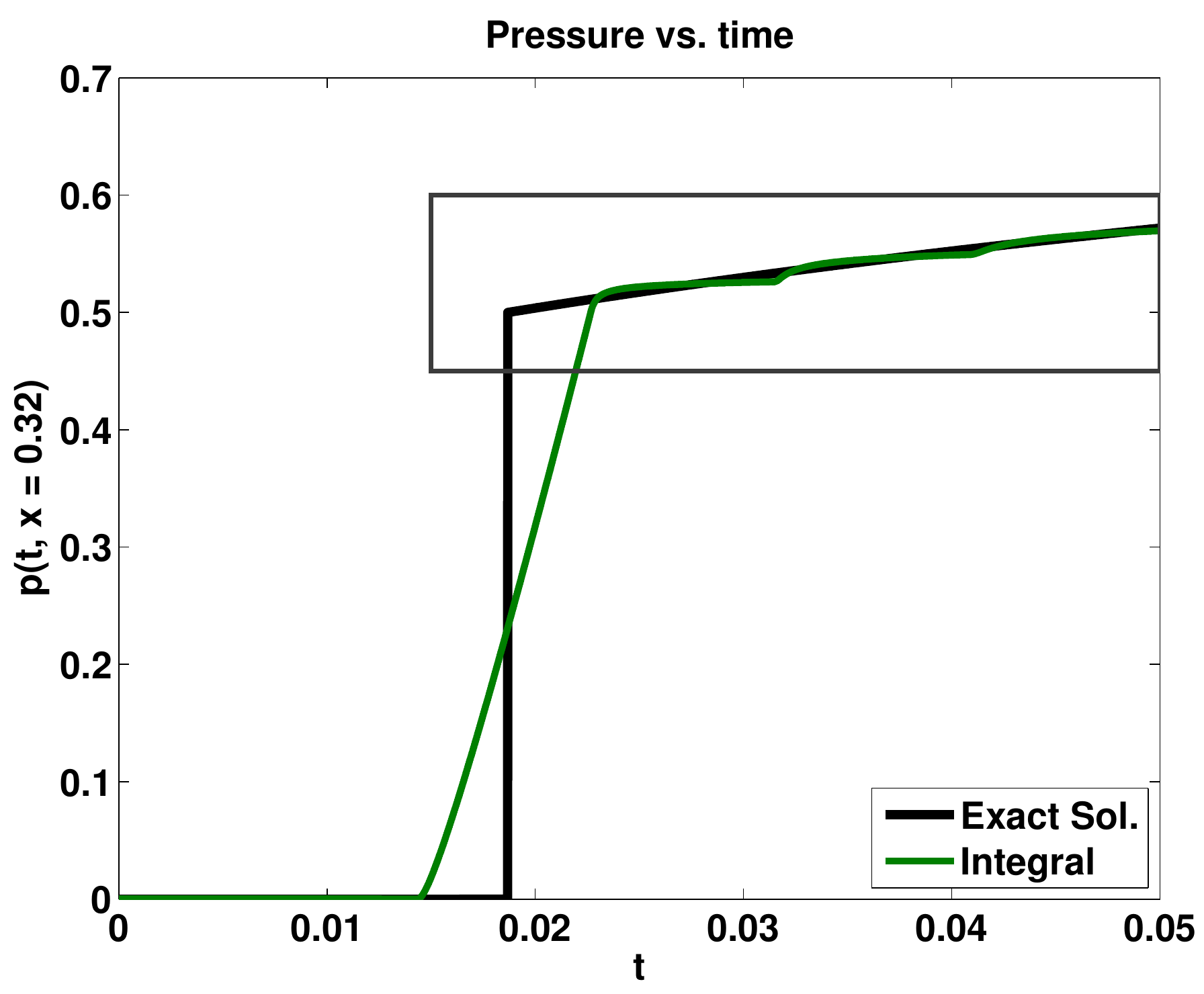}
		\end{subfigure}
		\begin{subfigure}[H]{0.49\textwidth}  
			\includegraphics[width =\textwidth]{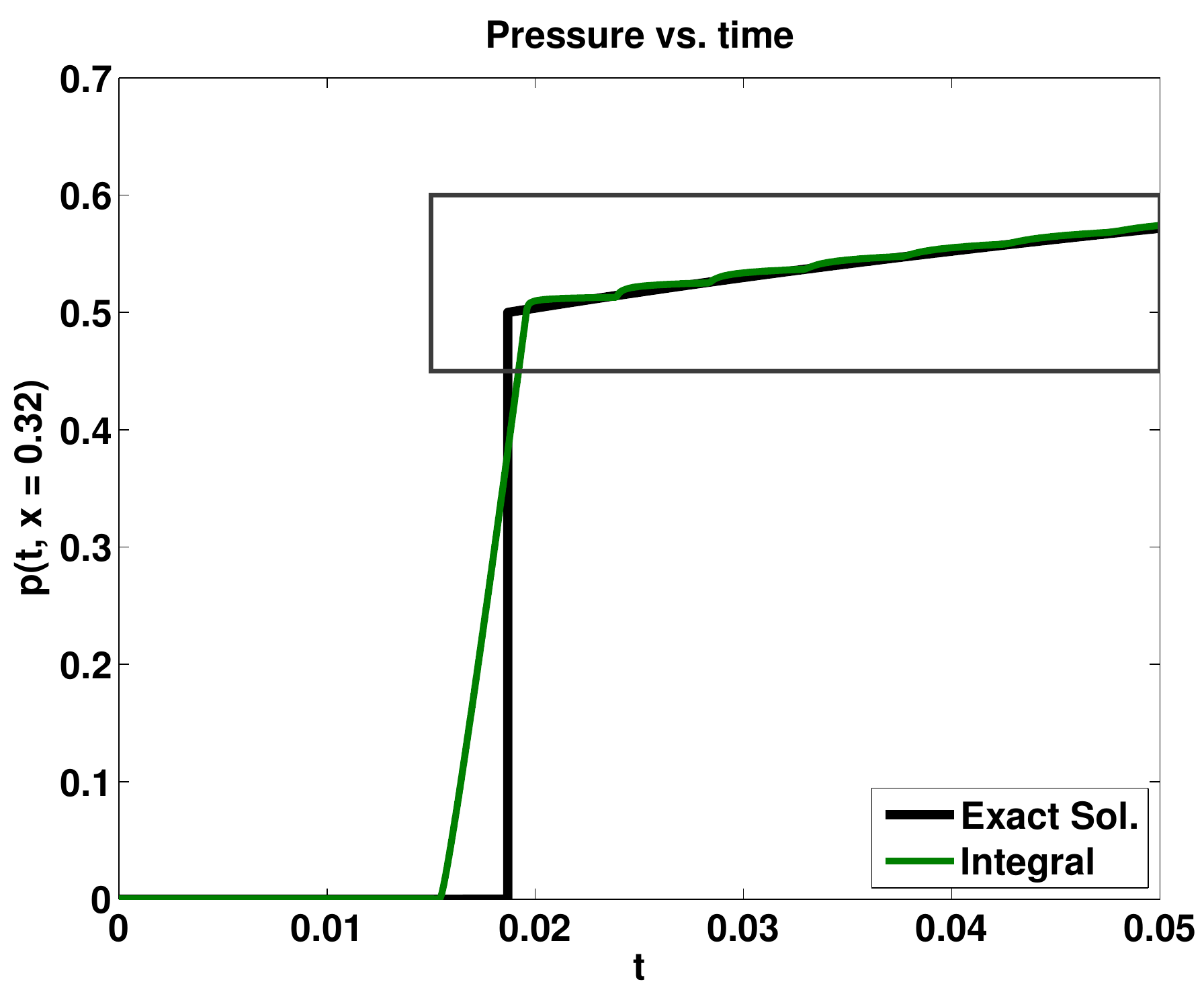}
		\end{subfigure}
		\begin{subfigure}[H]{0.49\textwidth}  
			\includegraphics[width =\textwidth]{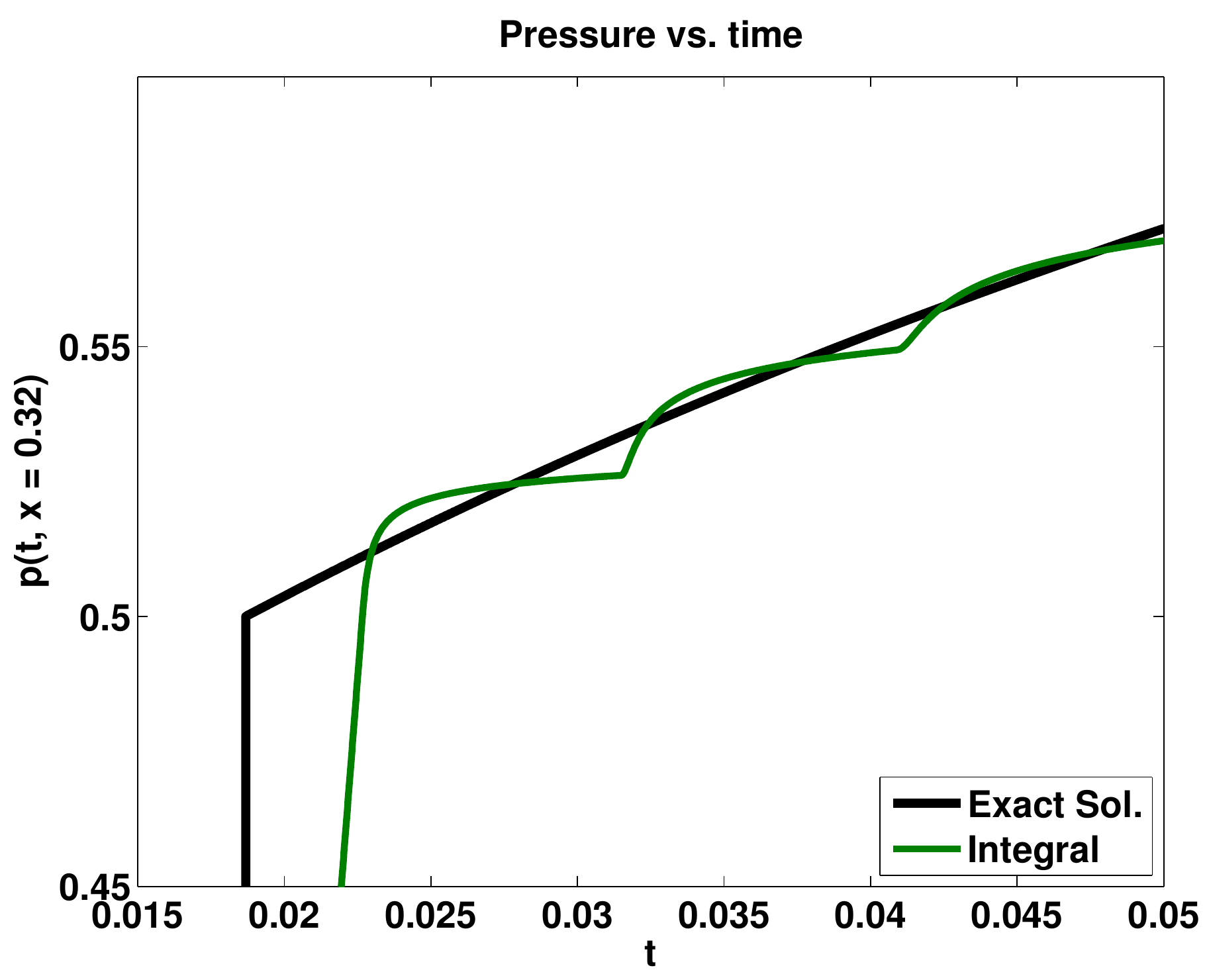}
			\caption{$N = 50$ grid points}
		\end{subfigure}
		\begin{subfigure}[H]{0.49\textwidth}  
			\includegraphics[width =\textwidth]{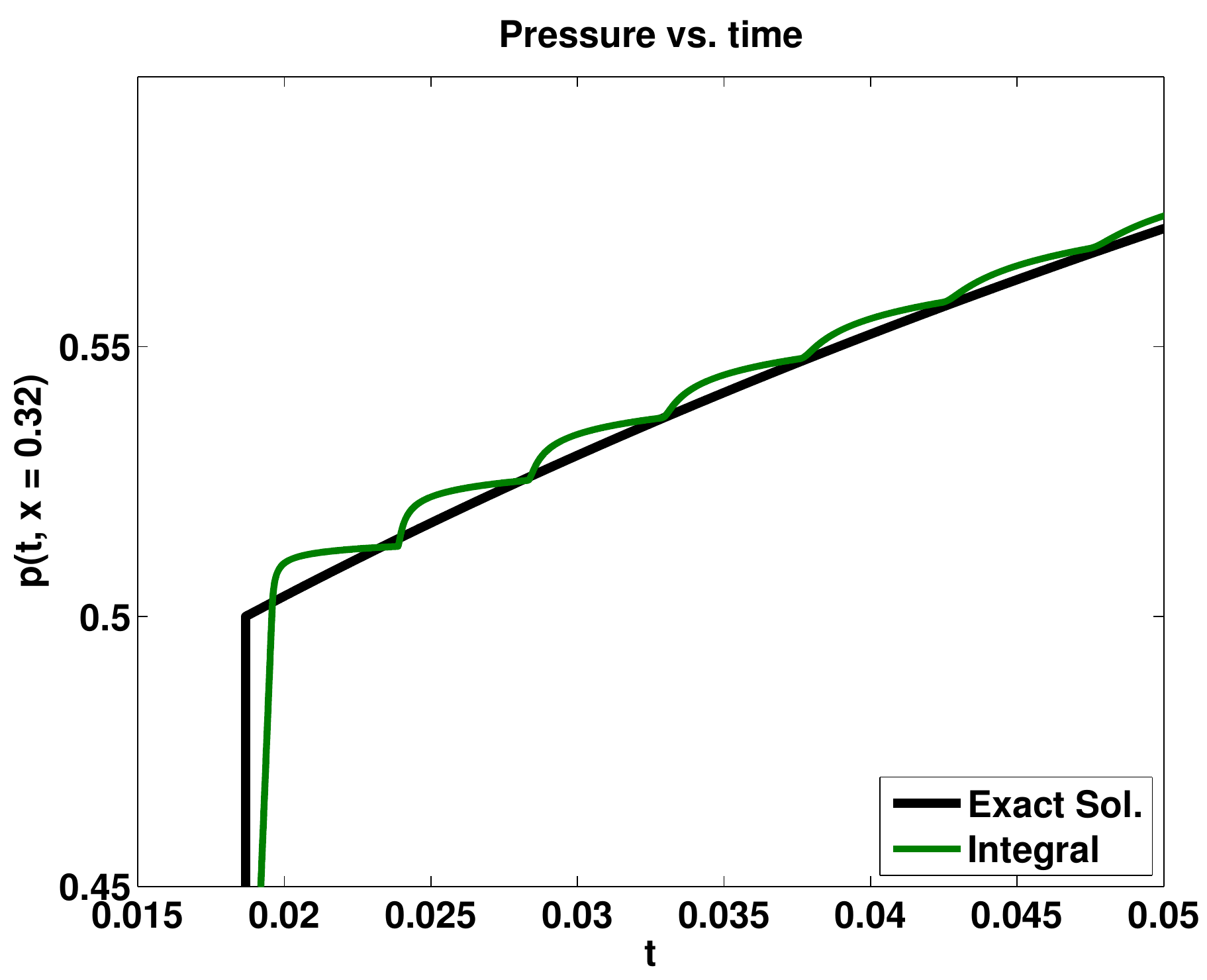}
			\caption{$N = 100$ grid points}
		\end{subfigure}
		\caption{Temporal profiles at position $x = 0.32$  with $\dt  = \dx^2/32$.}
		\label{fig:int_time}
\end{figure}
\begin{figure}[H]
		\center
		\begin{subfigure}[H]{.49\textwidth}  
			\includegraphics[width =\textwidth]{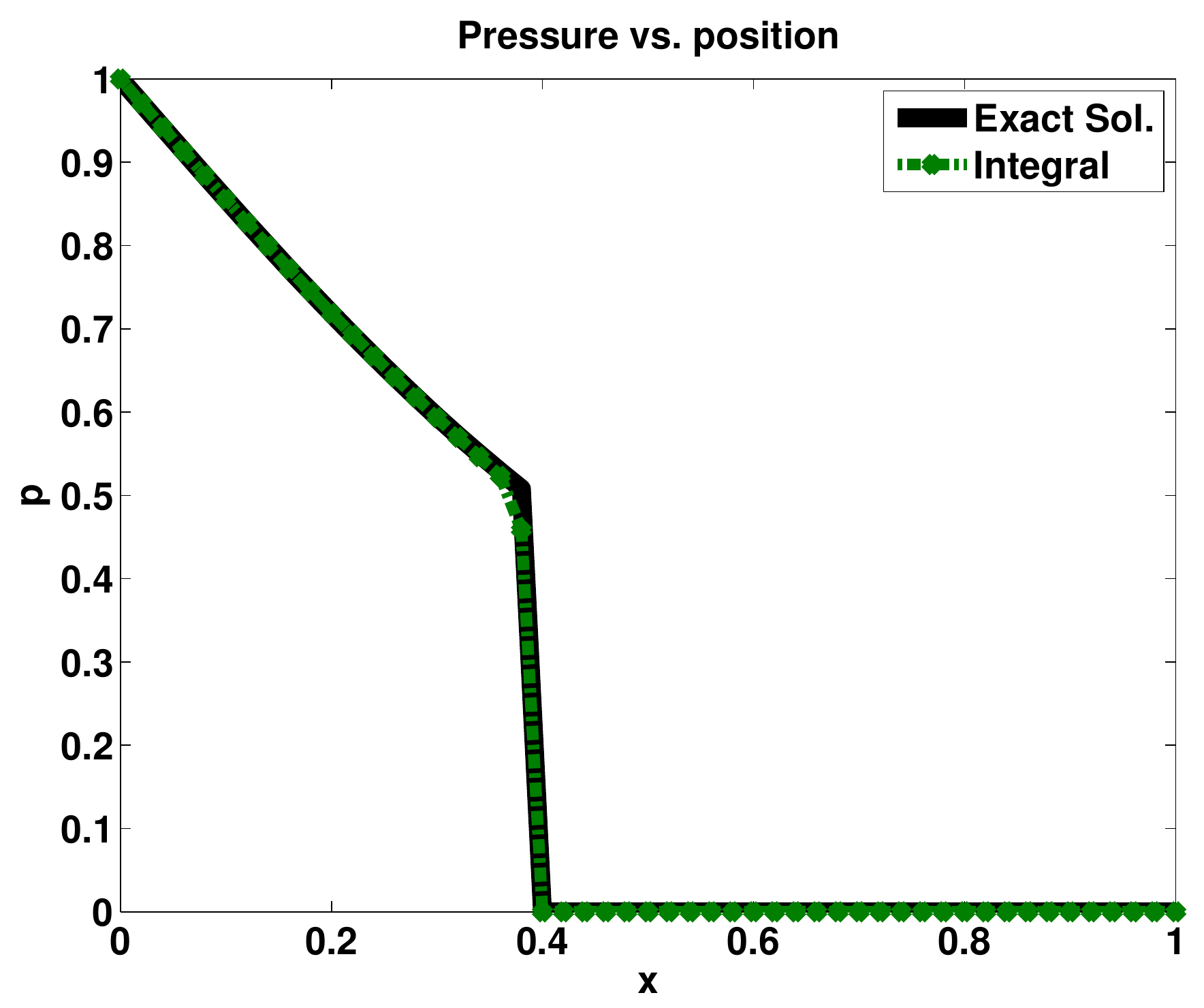}
			\caption{$N = 50$ grid points}
		\end{subfigure}
		\begin{subfigure}[H]{.49\textwidth}  
			\includegraphics[width =\textwidth]{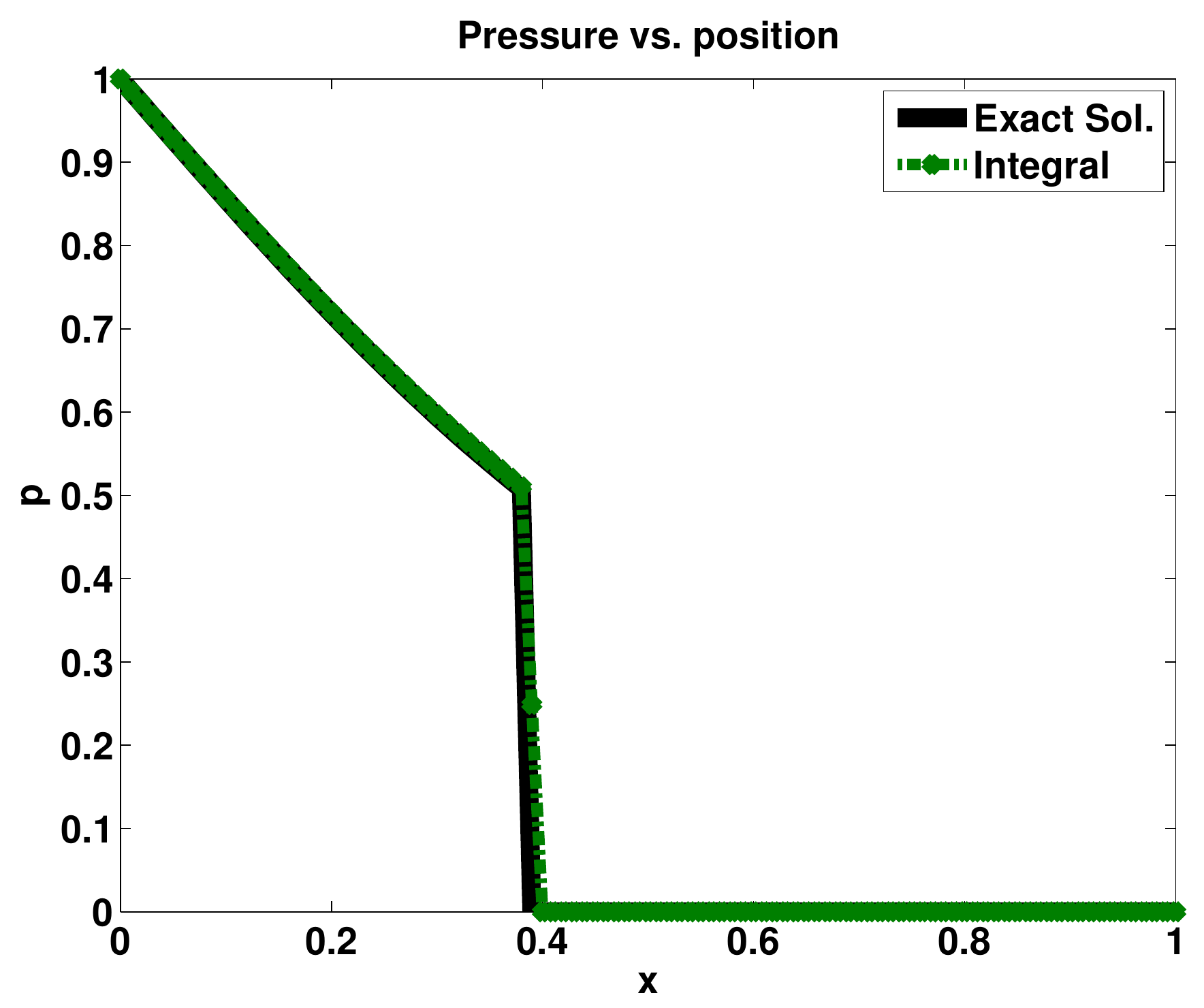}
			\caption{$N = 100$ grid points}
		\end{subfigure}
		\caption{Spatial profiles at time $t = 0.05$ with $\dt  = \dx^2/32$.}
		\label{fig:int_space}
\end{figure}

The numerical solution with the integral average is provided in Figures \ref{fig:int_time} and \ref{fig:int_space}.  The improvement with the integral average is clear.  It does not introduce high frequency oscillations near $p^*$.  
Although the low frequency, grid-dependent oscillations remain, the amplitude of these oscillations is smaller than those with arithmetic averaging.  The numerical solution is now monotonically increasing in time \cite{vandermeer2016}, matching the behavior of the true solution.  
The diffusive shock profile with the integral average is also illustrated in the spatial results in Figure \ref{fig:int_space}.

To explain why the integral average outperforms the arithmetic and harmonic averages, we write it in an alternate form than presented in \cite{vandermeer2016}.  The integral average $k^I_{i+1/2}$ in Eqn. \eqref{eq:int_avg} can be broken up at $p^*$ into two separate integrals to obtain 
\begin{equation}
k^I_{i+1/2} =
		\frac{k_{\min}(p_{i+1}-p^*) + k_{\max}(p^*-p_i)}{p_{i+1}-p_i}, 
	\label{eq:int_avg_val} 
\end{equation}
in the shock interval.  
It then contains some information about the shock, as encoded in the bounds of the integral. 
Unlike the harmonic and arithmetic averages, the integral average monotonically increases as the shock advances through the interval $[x_i, x_{i+1}]$.  

The corresponding continuous flux is given by 
\begin{equation}
	{F^+_{i}}^I = {F^-_{i+1}}^I = -\frac{k_{\min}(p_{i+1}-p^*) + k_{\max}(p^*-p_i)}{\dx}.
	\label{eq:int_flux}
\end{equation} 
Using this description, we can see that by utilizing $p^*$, the integral flux avoids computing the undefined gradient of $p$ across the jump, as done in the schemes with arithmetic and harmonic averaging.  
It is this incorporation of $p^*$ into the scheme that prevents the high frequency oscillations seen with arithmetic averaging, which explains the improved behavior.

As seen, the integral average does not remove the low frequency oscillations.  
This is because although $p^*$ is incorporated in the flux in Eqn. \eqref{eq:int_flux}, $x^*(t)$ is not.  In the case, where $k_{\min} = 0$, ${F^+_{i}}^I = {F^-_{i+1}}^I$ in Eqn. \eqref{eq:int_flux} appears to be similar to $F_{i}^+$ in Eqn. \eqref{eq:fluxi}.  The difference occurs in the denominator, where the relative shock position $\dx^*(t)$ detected in SAM is replaced with $\dx$ in the integral flux.  The integral flux is then assuming that the distance to the shock is fixed of size $\dx$.

 \begin{figure}[H]
		\center
		\begin{subfigure}[H]{0.49\textwidth}  
			\includegraphics[width =\textwidth]{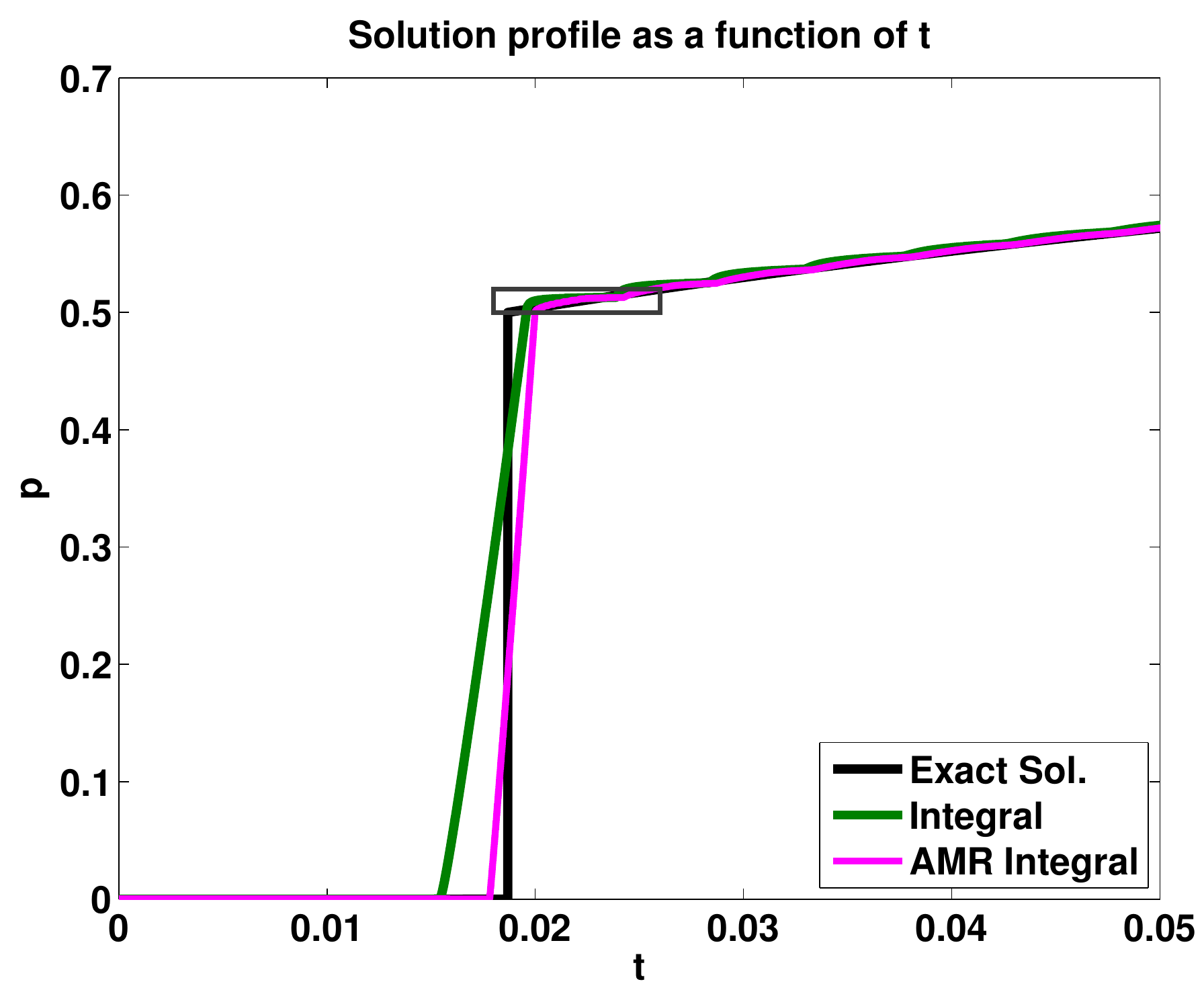}
		\end{subfigure}
		\begin{subfigure}[H]{0.49\textwidth}  
			\includegraphics[width =\textwidth]{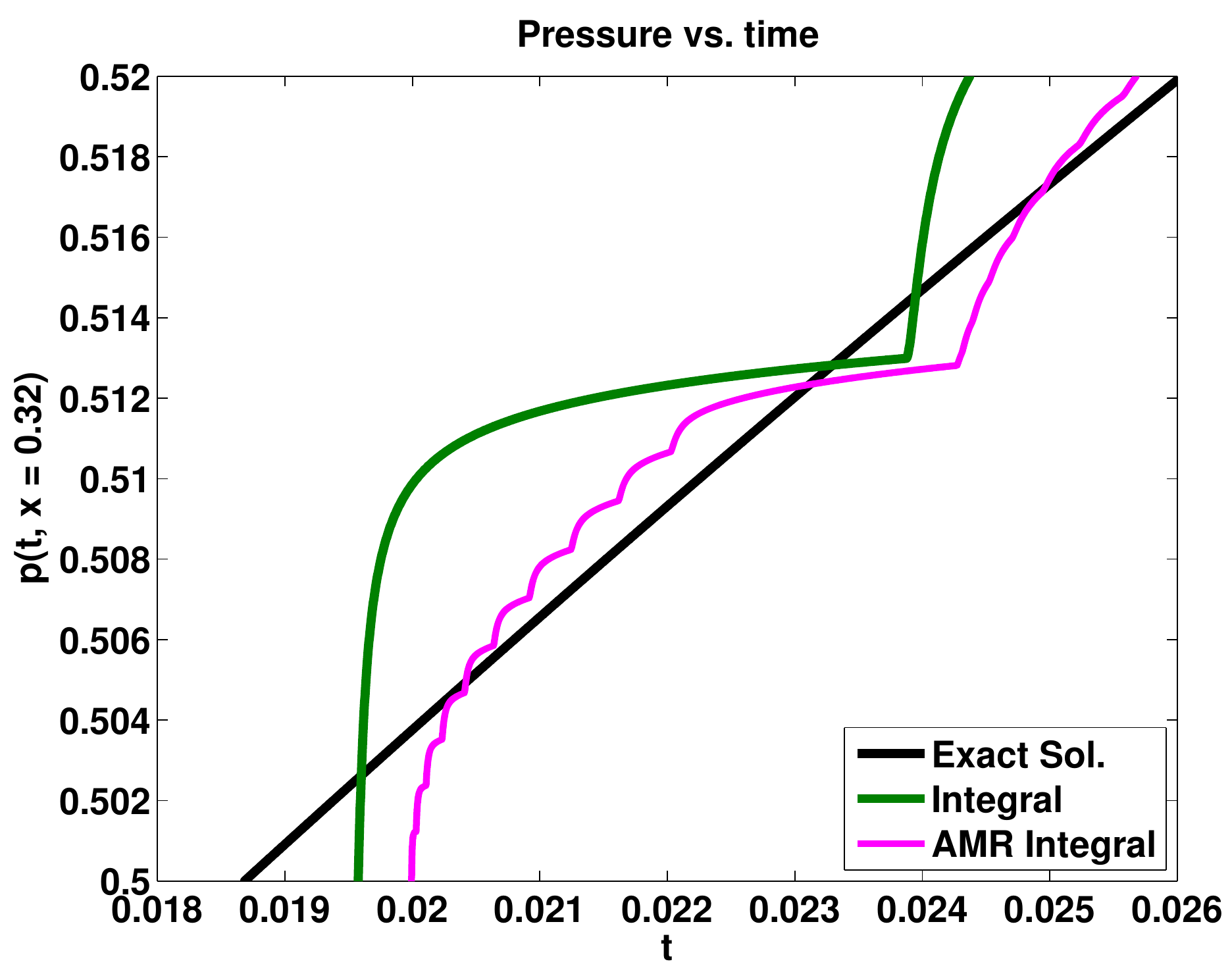}
		\end{subfigure}
		\caption{Temporal profiles at position $x = 0.32$ with $\dt  = \dx^2/32$ and $N = 100$ coarse grid points and 10 inner grid points.}
		\label{fig:AMR}
\end{figure}
To see what happens with the temporal oscillations under grid refinement, we also implement the scheme with integral averaging and Adaptive Mesh Refinement (AMR) \cite{berger84}.  We define a coarse mesh away from the shock and a fine mesh near the shock.  The fine mesh travels with the shock as it moves.  In Figure \ref{fig:AMR}, a coarse mesh size of $\dx = 1/100$ is utilized.  Within the coarse mesh cells surrounding the shock, a fine inner mesh size of $\dx_{\text{inner}} = 1/10$ is defined. 
Figure \ref{fig:AMR} displays that oscillations are present with a smaller period and damped amplitude.  
The zoomed in results in Figure \ref{fig:AMR} illustrate that there are $N_{\text{inner}} \equiv 1 / \dx_{\text{inner}}$ oscillations in between the coarse grid cells, as expected as the front crosses the $N_{\text{inner}} = 10$ inner grid cells for $\dx_{\text{inner}} = 1/10$.  The temporal oscillations are spatially dependent, and applying AMR does not remove them.

\subsection{SAM in Finite Volume Average Form and its Connection to the Integral Average}
\label{sam_discont}
 Writing SAM as an averaged-based method on a uniform Cartesian grid of spatial step size $\dx$ provides additional insight.  In particular, we will see another explanation of the temporal oscillations.  For $j \ne i$, we have defined $k^{SAM}_{j+1/2}$ by Eqn. \eqref{eq:gen_coeff} in Section \ref{sam_exact}.  The Cartesian and auxiliary grids are the same in the parabolic regions away from the shock.  
Using Eqns. \eqref{eq:fluxi} - \eqref{eq:fluxiplus}, we can write the out-flux $F_{i}^+$ of $CV_i$ and the in-flux $F_{i+1}^-$ of $CV_{i+1}$ for SAM in finite volume average form as
\begin{equation}
	\begin{aligned}
	\Huge F_{i}^+ &=   -\Bigg(\underbrace{k_{\max} \frac{\dx}{\dx^*}\frac{p^*-p_i}{p_{i+1}-p_i}}_{\mbox{\normalsize ${k^{SAM}_{i+1/2}}^+$ }}\Bigg)\frac{p_{i+1}-p_i}{\dx}, 
	\label{eq:sam_avgplus}
	\end{aligned}
\end{equation}
and
\begin{equation}
	\begin{aligned}
	F_{i+1}^- &= -\Bigg(\underbrace{k_{\min}\frac{\dx}{\dx - \dx^*}\frac{p_{i+1}-p^*}{p_{i+1}-p_i}}_{\mbox{\normalsize ${k^{SAM}_{i+1/2}}^-$ }}\Bigg)\frac{p_{i+1}-p_i}{\dx},
	\label{eq:sam_avgminus}
	\end{aligned}
\end{equation} 
respectively.  As opposed to the other averaging approaches, now $F_i^+ \ne F_{i+1}^-$: conservation is honored on the auxiliary grid, defined in Section \ref{aux_grid}.

We can express these SAM averages as a weighted linear combination  
\[
	c_1\frac{k_{\min}(p_{i+1}-p^*)}{p_{i+1}-p_i} + c_2\frac{k_{\max}(p^*-p_i)}{p_{i+1}-p_i},
\]  
of the integral average $k^{I}_{i+1/2} $ in Eqn. \eqref{eq:int_avg_val}.  
For  $c_1 = 0$ and $c_2 = \dx/\dx^*$, ${k^{SAM}_{i+1/2}}^+$ in Eqn. \eqref{eq:sam_avgplus} is recovered, while for $c_1 = \dx/(\dx - \dx^*)$ and $c_2 = 0$, ${k^{SAM}_{i+1/2}}^-$ in Eqn. \eqref{eq:sam_avgminus} is recovered.  Conversely, we can also express the integral flux in terms of the SAM fluxes surrounding the shock cell. 
Rewriting Eqn. \eqref{eq:int_flux}, we have
\begin{equation}
\begin{aligned}
	{F^+_{i}}^I &=  - \Big[k_{\min}\frac{p_{i+1}-p^*}{\dx-\dx^*}\frac{\dx-\dx^*}{\dx} + k_{\max}\frac{p^*-p_i}{\dx^*}\frac{\dx^*}{\dx}\Big] \\
	&=  F_{i+1}^-(1-y) + F_i^+y, 
	\end{aligned}
	\label{eq:int_anal}
\end{equation}
where $y \equiv \dx^*/\dx$.  The shock moves from left to right in the interval $[x_i,x_{i+1}]$, and $0 < y < 1$.  

As discussed in Section \ref{aux_grid}, the SAM flux is a first order approximation of the analytical flux, and the jump condition is satisfied.  The integral flux, on the other hand, does not converge to the analytical flux as $\dx \rightarrow 0$, and does not satisfy the jump condition.  Without loss of generality, we let the shock be to the right of the cell face $x_{i+1/2}$ by $\eps\dx$ for some small $\eps > 0$.  Then $y = 1/2 + \eps$.  Using Taylor series expansion in the parabolic region to the left of the shock gives 
\begin{equation}
	F_i^+ = F^L[p(x_{i+1/2+\eps})] + \mathcal{O}(\dx),
	\label{eq:int_analplus}
\end{equation}
where $F^L[p(x)]= -k_{\max}\nabla p(x)$.
Similarly,  Taylor series expansion in the parabolic region to the right of the shock gives
\begin{equation}
	F_{i+1}^- = F^R[p(x_{i+1/2 + \eps})] + \mathcal{O}(\dx),
	\label{eq:int_analminus}
\end{equation}
where $F^R[p(x)] = -k_{\min} \nabla p(x)$.
Substituting Eqns. \eqref{eq:int_analplus}-\eqref{eq:int_analminus} into Eqn. \eqref{eq:int_anal} gives the following 
\[
	{F^+_{i}}^I = F^R[p(x_{i+1/2+\eps})](0.5-\eps) +  F^L[p(x_{i+1/2+\eps})](0.5+\eps) + \mathcal{O}(\dx).
\]
Taking the limit as $\eps \rightarrow 0$, we obtain 
\[
	{F^+_{i}}^I = \frac{F^R[p(x_{i+1/2})] +  F^L[p(x_{i+1/2})]}{2} + \mathcal{O}(\dx),
\]
as the integral flux at the interface $x_{i+1/2}$.  Since $x_{i+1/2}$ is to the left of the shock, the analytical flux at $x_{i+1/2}$ is given by the left flux $F^{\text{L}}[p(x_{i+1/2})]$.  The integral flux is then averaging the left and right fluxes across the jump.  The right flux $F^R[p(x_{i+1/2})]$ should not have any impact on the value across the jump.  This jump condition violation is exactly what results in the temporal oscillations.  SAM does satisfy the jump condition and as the results in Sections \ref{exact_shock} and \ref{shock_speed} show, it does remove the artifacts entirely.  

\label{flux}

\section{Conclusions and Future Work} 
This paper explains what causes the mysterious artifacts in discretizations of the discontinuous Generalized Porous Medium Equation (GPME), and suggests an alternative method that results in numerical solutions without these artifacts.  The FTCS scheme with integral averaging performs better than the schemes with harmonic and arithmetic averaging because it has information about the shock value $p^*$.  By rewriting the integral average in the shock cell, it can be seen that full removal of the numerical artifacts requires more than $p^*$, and that the shock location $x^*(t)$ must also be included in the numerical scheme.  
The Shock-Based Averaging Method (SAM) incorporates both $x^*(t)$ and $p^*$, and satisfies the jump condition for integral conservation laws to result in numerical solutions with accurate and smooth temporal profiles.  Casting SAM in the finite-volume framework helps provide understanding on why the numerical artifacts were occurring with finite volume averaged-based methods.

Future work includes extensions to higher dimensions and to the Porous Medium Equation (PME) subclass.  Since the velocity for the GPME is given for arbitrary dimensions, it can be used as input to a level set implementation for extensions to higher dimensions, as discussed in Section \ref{shock_speed}.  
Another future direction would be to combine the approach from this paper with our previous paper \cite{maddix_pme}, discussing the PME subclass of the GPME with continuous coefficients.  
For the PME, there is no known specified $p$ value at the shock, as given by $p^*$ in this paper.  
With the additional control volume around the shock in the auxiliary finite volume grid from Section \ref{sam_exact}, $p^*(t)$ can be solved for as an additional degree of freedom, and an approach similar to SAM can be taken.





\appendix
\section{Exact Solution Derivation for arbitrary parameters}
\label{exact}
For the test problem with the initial condition in Figure \ref{IC}, the exact solution is known and is used for testing the numerical methods.  We give the full derivation, which is similar to \cite{Kraaijevanger2011}, and is generalized for arbitrary $k_{\max}$.  The domain is partitioned into $\Omega_1 = (0, x^*(t)]$ and $\Omega_2 = [x^*(t), \infty)$.  The problem can be subdivided into two constant coefficient heat equations  \cite{Kraaijevanger2011, crank1984}.  The solution is monotonically non-increasing and so in $\Omega_1, p_1(x,t) > p^*$ and in $\Omega_2, p_2(x,t) \le p^*$.  It can be verified that $p_1(x,t) = 1 - c_1\Phi\Big(x/(2\sqrt{k_{\max}t})\Big)$ and $p_2(x,t) = c_2\Big(1 - \Phi\Big(x/(2\sqrt{k_{\min}t})\Big)\Big)$ for some constants $c_1, c_2$ to be determined.  $\Phi(x) = \text{erf}(z) = \int_0^z \phi(y) dy$ is the standard Gaussian error function, where $\phi(y) = \frac{2}{\sqrt{\pi}} \exp(-y^2)$. 

The unknown shock location ($x^*(t)$) needs to be computed.  To do so, two additional boundary conditions are required at the shock $x = x^*(t)$: 
\begin{enumerate}
	\item $p_1(x^*(t), t) = p_2(x^*(t),t) = p^*, \hspace{.25cm} \forall t$ (Continuity)
	\item  $k_{\max}(\partial{p_1(x^*(t), t)}/\partial{x}) = k_{\min}(\partial{p_2(x^*(t), t)}/\partial{x}), \hspace{.25cm}  \forall t$ (Flux Continuity).
\end{enumerate}

Since Condition (1) must hold for all $t$, $x^*(t) = \alpha \sqrt{t}$, 
\[
	c_1 = \frac{1-p^*}{\Phi\Big(\alpha/(2\sqrt{k_{\max}})\Big)}, \hspace{.5cm} c_2 = \frac{p^*}{1-\Phi\Big(\alpha / (2\sqrt{k_{\min}})\Big)}.
\] 

Condition (2) is used to derive a nonlinear solve for the remaining unknown, $\alpha$.  Substituting in the expressions for the derivatives and simplifying leads to
\begin{equation}
	c_1\sqrt{k_{\max}}\phi\Big(\frac{\alpha}{2\sqrt{k_{\max}}}\Big) = c_2\sqrt{k_{\min}}\phi\Big(\frac{\alpha}{2\sqrt{k_{\min}}}\Big).
	\label{eq:flux_eq_bc}
\end{equation}
Let $z_1 = \alpha/(2\sqrt{k_{\max}})$ and $z_2 = \alpha/(2\sqrt{k_{\min}})$.  We substitute the expressions for $c_1$ and $c_2$ into Eqn. \eqref{eq:flux_eq_bc}, multiply both sides by  $\frac{\alpha}{2}$ and simplify to obtain
\begin{equation}
	\begin{aligned}
		\frac{1-p^*}{\Phi(z_1)} \sqrt{k_{\max}}\phi(z_1) &=  \frac{p^*}{1-\Phi(z_2)}\sqrt{k_{\min}}\phi(z_2) \iff \\
								(1-p^*)(1-\Phi(z_2))\exp(z_2^2)\frac{\alpha}{2\sqrt{k_{\min}}} &=  p^*\frac{\alpha}{2\sqrt{k_{\max}}}\Phi(z_1)\exp(z_1^2) \iff \\
								(1-p^*)(1-\Phi(z_2))\exp(z_2^2)z_2 &=  p^*\Phi(z_1)z_1\exp(z_1^2) .
		\label{eq:flux_eq_long}
	\end{aligned}
\end{equation}
By a series expansion from integration by parts, 

\[
	1-\Phi(z_2) = \frac{\exp(-z_2^2)}{z_2\sqrt{\pi}}\Bigg(1 - \frac{1}{2z_2^2} + \frac{3}{4z_2^4} - \dots\Bigg).
\] Using the above expression, Eqn. \eqref{eq:flux_eq_long} simplifies to
					\begin{equation}
					\frac{1 - p^*}{\sqrt{\pi}}\Bigg(1 - \frac{1}{2z_2^2} + \frac{3}{4z_2^4} - \dots\Bigg) =  p^*\Phi(z_1)z_1\exp(z_1^2).
					\label{eq:limit_jbp}
					\end{equation}

In our application and most of the numerical tests, $k_{\min} = 0$.  This implies that $z_2 \rightarrow \infty$.  The limit as $z_2 \rightarrow \infty$ of the left hand side of Eqn. \eqref{eq:limit_jbp} is simply $(1-p^*)/\sqrt{\pi}$.  This is an advantage of computing the series form, since we can easily take this limit as $k_{\min} \rightarrow 0$.   This form is also preferred numerically when $k_{\min}$ is small to avoid multiplication of the large $\exp(z_2^2)z_2$ terms.  For the $k_{\min} = 0$ case, the unknown ($\alpha$) only appears in $z_1$ on the right hand side of the equation, given by 
\begin{equation}
	\frac{1 - p^*}{\sqrt{\pi}}=  p^*\Phi(z_1)z_1\exp(z_1^2).
\label{eq:alpha}
\end{equation}
A simple one-dimensional nonlinear equation solver can be used to solve this equation for $z_1$, where $\alpha = 2\sqrt{k_{\max}}z_1$.


\section{Error Tables and Convergence Study}
\label{app:conv}
\begin{figure}[H]
		\center
		\begin{subfigure}[H]{.49\textwidth}  
			\includegraphics[width =\textwidth]{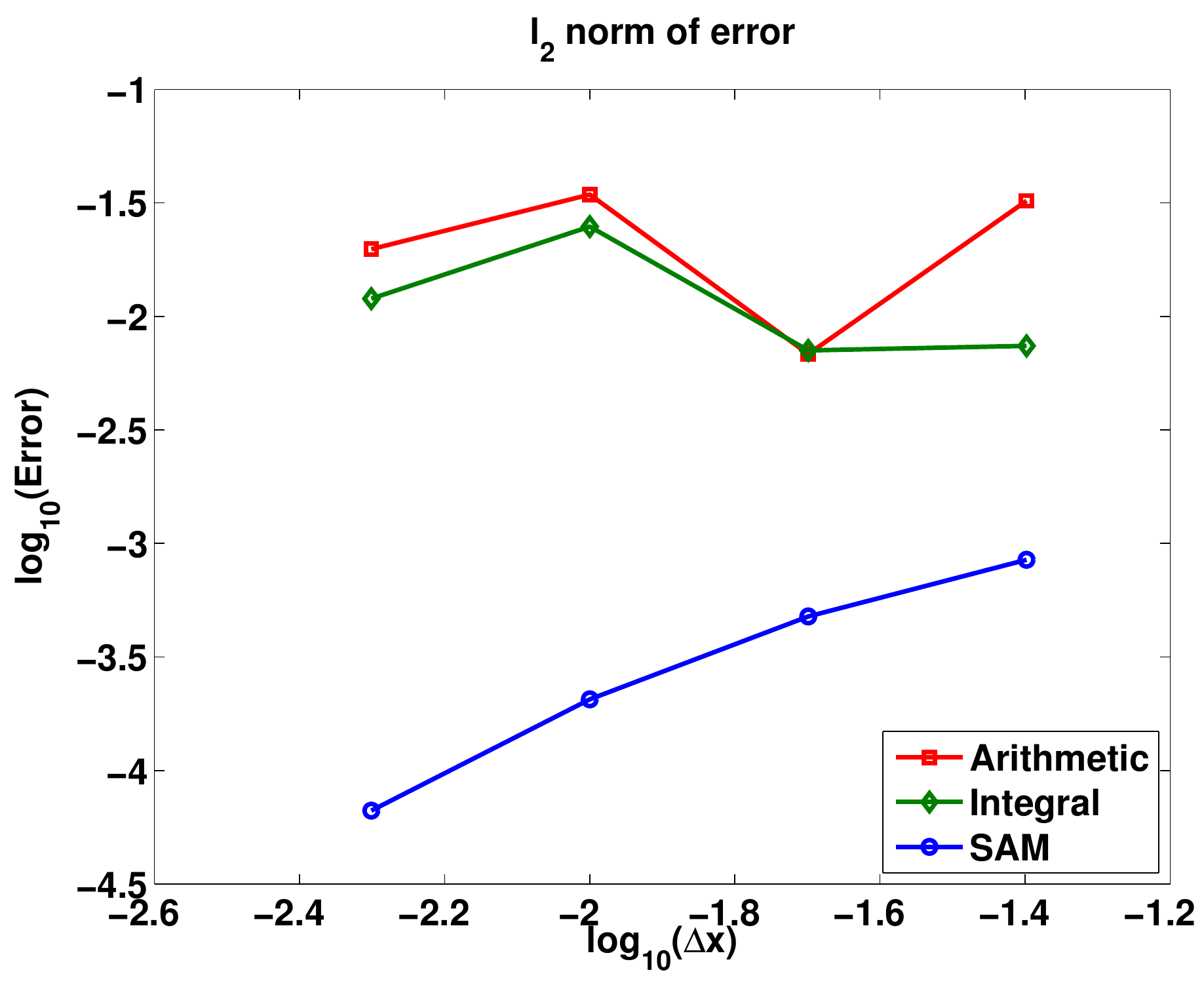}
			\caption {$l_2$ norm}
			\label{fig:log1}
			\end{subfigure}
			\begin{subfigure}[H]{.49\textwidth}  
			\includegraphics[width =\textwidth]{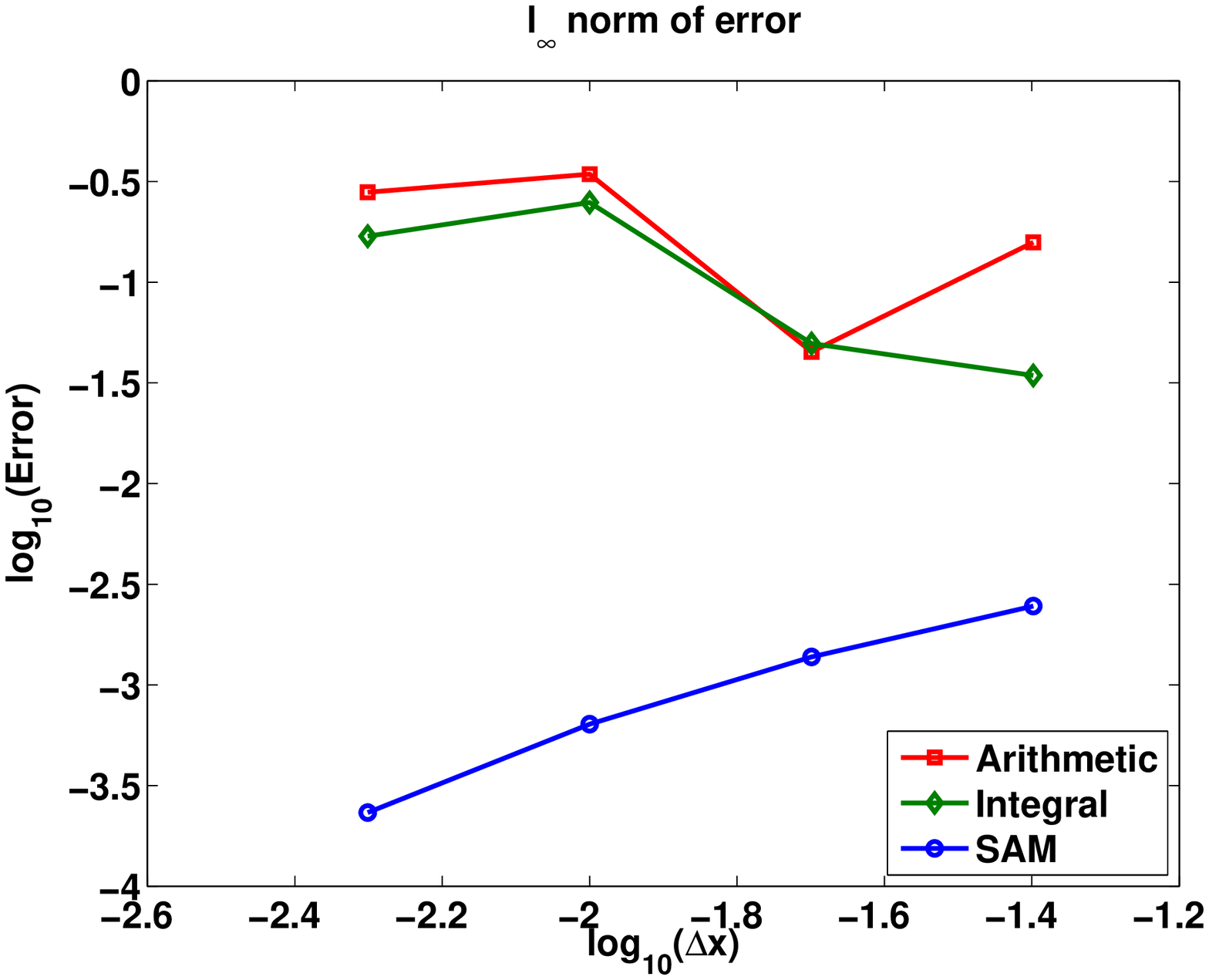}
			\caption {$l_\infty$ norm}
			\label{fig:log2}
			\end{subfigure}
		\caption{Loglog convergence plots for the Stefan problem.}
\end{figure}

This appendix contains the convergence results for the model problem tested with the same time step size $\dt = \dx^2/32$ as in the prior sections.  The $l_2$ and $l_\infty$ error norms are calculated with respect to the exact solution on the corresponding grid.  The error is measured globally, and includes the points surrounding the shock.  We compare the errors of the FTCS with arithmetic and integral averaging to those of SAM.  The errors of the FTCS with harmonic averaging are not shown, since the solution locks and does not converge.  For the methods with arithmetic and integral averaging, even on finer grids, the asymptotic region is not yet reached.  Tables \ref{table:timing}-\ref{table:timing2} show that SAM has approximately first order convergence and errors that are orders of magnitude lower than those for the averaged-based methods.  The first order convergence of SAM is also shown in Figures \ref{fig:log1}-\ref{fig:log2}.  We do not consider higher order methods because of the dominance of the shock.
\begin{table}[H]
\begin{center}
\caption{$l_2$ norm errors}
\small
\begin{tabular}{|c|c|c|c|c|c|c|}  
\hline 
& $N = 25$ & $N = 50$ & $N = 100$ &$N = 200$  &Order \\  \hline
Arithmetic			&3.2335e-02         &6.8351e-03& 3.4431e-02 & 1.9773e-02 &N/A  \\ \hline
Integral       &7.4250e-03        & 7.0866e-03& 2.4898e-02  &1.1974e-02  & N/A   \\  \hline
SAM  		&8.4813e-04 &4.7762e-04 &2.0583e-04 &6.6594e-05 & 1.2229
   \\ \hline
\end{tabular}
\label{table:timing}
\end{center}
\end{table} 

\begin{table}[H]
\begin{center}
\caption{$l_{\infty}$ norm errors}
\small
\begin{tabular}{|c|c|c|c|c|c|c|}  
\hline 
& $N = 25$ & $N = 50$ & $N = 100$ &$N = 200$  &Order \\  \hline
Arithmetic		&4.8237e-02         &1.5223e-02 &3.5742e-02 & 2.0394e-02&N/A \\ \hline
Integral	      &1.4104e-02          &9.1485e-03 &2.9971e-02 &1.6252e-02 &N/A     \\  \hline
SAM	  		&2.1859e-03  &1.7763e-03 &1.0586e-03  & 4.5134e-04	
 & 1.1336   \\ \hline
\end{tabular}
\label{table:timing2}
\end{center}
\end{table} 

\section*{Acknowledgements}
This material is based upon work supported by the National Science Foundation Graduate Research Fellowship under Grant No. DGE - 114747.  We would also like to thank Giuseppe Domenico Cervelli for his guidance.
\bibliographystyle{elsarticle-num-names}
 \bibliography{Foamreferences}
\end{document}